\documentclass[11pt,oneside,english,reqno]{amsart}
\usepackage[latin9]{inputenc}
\usepackage{geometry}
\usepackage{amstext}
\usepackage{amssymb}
\usepackage{setspace}
\makeatletter
\usepackage{amsthm}
\usepackage{amstext}
\usepackage{amsfonts}

\numberwithin{figure}{section}
\@ifundefined{date}{}{\date{}}
\numberwithin{figure}{section}
\newtheorem{thm}{Theorem}

\newtheorem{lem}{Lemma}

\newtheorem*{asmM}{Assumption M}
\newtheorem*{asmM3'}{Assumption M (iii')}
\newtheorem*{asmS}{Assumption S}

\newtheorem*{lemC}{Lemma C}
\newtheorem*{lemM}{Lemma M}
\newtheorem*{lemM'}{Lemma M'}
\newtheorem*{lemM1}{Lemma M1}
\newtheorem*{lemMS}{Lemma MS}
\newtheorem*{lemMS'}{Lemma MS'}
\newtheorem*{def-cube}{Definition (Generalized cube root class)}
\newtheorem*{def-cube-set}{Definition (Partially identified cube root class)}
\newtheorem*{asmD}{Assumption D}
\theoremstyle{definition}

\makeatother

\usepackage{babel}
\begin{document}
	
	\title{Local M-estimation with Discontinuous Criterion for Dependent and
		Limited Observations}
	
	\author{Myung Hwan Seo}
	
	\address{Department of Economics, Seoul National University, Seoul, Korea.}
	
	\email{myunghseo@snu.ac.kr}
	
	\author{Taisuke Otsu}
	
	\address{Department of Economics, London School of Economics, Houghton Street,
		London, WC2A 2AE, UK.}
	
	\email{t.otsu@lse.ac.uk}
	
	\thanks{The authors would like to thank Aureo de Paula, Marine Carrasco,
		Javier Hidalgo, Dennis Kristensen, Benedikt Pötscher, Peter Robinson,
		Kyungchul Song, Yoon-Jae Whang, and seminar and conference participants
		at Cambridge, CIREQ in Montreal, CORE in Louvain, CREATES in Aarhus,
		IHS in Vienna, ISNPS in Cádiz, LSE, Surrey, UCL, Vienna, and York
		for helpful comments. The authors also acknowledge helpful comments
		from an associate editor and anonymous referees. This research was
		partly supported by Promising-Pioneering Researcher Program through
		Seoul National University (Seo) and the ERC Consolidator Grant (SNP
		615882) (Otsu).}
	
	\begin{abstract}
		This paper examines asymptotic properties of local M-estimators under
		three sets of high-level conditions. These conditions are sufficiently
		general to cover the minimum volume predictive region, conditional
		maximum score estimator for a panel data discrete choice model, and
		many other widely used estimators in statistics and econometrics.
		Specifically, they allow for discontinuous criterion functions of
		weakly dependent observations, which may be localized by kernel smoothing
		and contain nuisance parameters whose dimension may grow to infinity.
		Furthermore, the localization can occur around parameter values rather
		than around a fixed point and the observation may take limited values,
		which leads to set estimators. Our theory produces three different
		nonparametric cube root rates and enables valid inference for the
		local M-estimators, building on novel maximal inequalities for weakly
		dependent data. Our results include the standard cube root asymptotics
		as a special case. To illustrate the usefulness of our results, we
		verify our conditions for various examples such as the Hough transform
		estimator with diminishing bandwidth, maximum score-type set estimator,
		and many others.
	\end{abstract}
	
	\maketitle

\section{Introduction}

There is a class of estimation problems in statistics where a point (or set-valued) estimator is obtained by maximizing a discontinuous and possibly localized criterion function. As a prototype, consider estimation of a simplified version of the minimum volume predictive region for $y$ at $x=c$ (Polonik and Yao, 2000). Let $\mathbb{I}\{\cdot\}$ be the indicator function, $K(\cdot)$ be a kernel function, and $h_{n}$ be a bandwidth. At a significance level $\alpha$, the estimator $[\hat{\theta}\pm\hat{\nu}]$ is obtained by the M-estimation
\begin{equation}
\max_{\theta\in\Theta}\sum_{t=1}^{n}\mathbb{I}\{|y_{t}-\theta|\leq\hat{\nu}\}K\left(\frac{x_{t}-c}{h_{n}}\right),\label{eq:PY}
\end{equation}
where $\Theta$ is some parameter space and
\[
\hat{\nu}=\inf\left\{ \nu\in\mathbb{R}:\max_{\theta\in\Theta}\frac{\sum_{t=1}^{n}\mathbb{I}\{|y_{t}-\theta|\leq\nu\}K\left(\frac{x_{t}-c}{h_{n}}\right)}{\sum_{t=1}^{n}K\left(\frac{x_{t}-c}{h_{n}}\right)}\geq\alpha\right\} .
\]
This problem exhibits several distinguishing features such as discontinuity of the criterion function, localization by kernel smoothing, and serial dependence in time series data which have prevented a full-blown asymptotic analysis of the M-estimator $\hat{\theta}$. Only consistency is reported in the literature.

This type of M-estimation has numerous applications. Since Chernoff's (1964) study on estimation of the mode, many papers have raised such estimation problems, for example the shorth (Andrews \emph{et al.,} 1972), least median of squares (Rousseeuw, 1984), nonparametric monotone density estimation (Prakasa Rao, 1969), and maximum score estimation (Manski, 1975). These classical examples are studied in a seminal work by Kim and Pollard (1990), which explained elegantly how this type of estimation problem induces so-called cube root asymptotics in a unified framework by means of empirical process theory. See also van der Vaart and Wellner (1996) and Kosorok (2008) for a general theory of M-estimation via empirical processes. However, these works do not cover the estimation problem in (\ref{eq:PY}) due to their focus on cross-sectional data among other things. It should be emphasized that this is not a pathological example. We provide various relevant examples in Section \ref{sec:ex} and Supplement (Section B) including the well-known Honor\'{e} and Kyriazidou's (2000) estimator for a dynamic panel discrete choice model and a localized maximum score estimator for a new binary choice model with random coefficients.

This paper covers a broader class of M-estimators than the above examples suggest. The baseline scenario (called \emph{local M-estimation}) is generalized in two directions. First, we accommodate not only variables taking limited values (e.g., interval-valued data) which typically lead to estimation of a set rather than a point, but also nuisance parameters with growing dimension. Set estimation problems due to limited observations are also known as partial identification problems in econometrics (e.g., Manski and Tamer, 2002). It is also novel to accommodate high-dimensional nuisance parameters in M-estimation with discontinuous criterion functions. Second, we allow for localization to be dependent on parameter values instead of prespecified values. For instance, the criterion function may take the form of $\sum_{t=1}^{n}\mathbb{I}\{|y_{t}-\theta|\leq h_{n}\}$ with $h_{n}\rightarrow 0$. Relevant examples include mode estimation (Chernoff, 1964, and Lee, 1989) and the Hough transform estimator in image analysis (Goldenshluger and Zeevi, 2004). Henceforth we call this case \emph{parameter-dependent local M-estimation}. Parameter-dependence brings some new features into our asymptotic analysis but in a different way from a classical example of parameter-dependency on the support such as the maximum likelihood estimator for $\mathrm{Uniform}[0,\theta]$.

The main contribution of this paper is to develop a general asymptotic theory for such M-estimation problems. Our theoretical results cover all the examples above and can be used to establish limit laws for
point estimators and convergence rates for set estimators. To this end, we develop suitable maximal inequalities which enable us to obtain nonparametric cube root rates of $(nh_{n})^{1/3}$, $\{nh_{n}/\log(nh_{n})\}^{1/3}$, and $(nh_{n}^{2})^{1/3}$ for the cases of local M-estimation, limited observations, and parameter-dependent localization, respectively. These inequalities are extended to establish stochastic asymptotic equicontinuity of normalized processes of the criterion functions so that an argmax theorem delivers limit laws of the M-estimators. It is worth noting that all the conditions are characterized through moment conditions and can be easily verified as illustrated in the examples. Thus, our results can be applied without prior knowledge of empirical process theory. It is often not trivial to verify entropy conditions such as uniform manageability in Kim and Pollard (1990). Particularly for dependent data, the covering, or bracketing, numbers often need to be calculated using a norm that hinges on the mixing coefficients and distribution of the data (e.g., the $L_{2,\beta}$-norm in Doukhan, Massart and Rio, 1995).

Another contribution is that we allow for weakly dependent data which are associated with absolutely
regular processes with exponentially decaying mixing coefficients. In some applications, the cube root asymptotic theory has been extended to time series data, for example Anevski and H\"{o}ssjer (2006) for monotone density estimation, Zinde-Walsh (2002) for least median of squares, de Jong and Woutersen (2011) for maximum score, and Koo and Seo (2015) for break estimation under misspecification. However, it is not clear whether they are able to handle the general class of estimation problems in this paper. 

The paper is organized as follows. Section 2 develops an asymptotic theory for local M-estimation and Section 3 provides several examples. In Section 4, we generalize the asymptotic theory to the cases of limited observations (Section 4.1) and parameter-dependent localization (Section 4.2). Section 5 concludes. All proofs, details for illustrations, and additional examples are contained in the Supplement.

\section{Local M-estimation\label{sec:general}}

This section studies the M-estimator $\hat{\theta}$ that maximizes 
\[
\mathbb{P}_{n}f_{n,\theta}=\frac{1}{n}\sum_{t=1}^{n}f_{n,\theta}(z_{t}),
\]
where $\{f_{n,\theta}\}$ is a sequence of criterion functions indexed by the parameters $\theta\in\Theta\subseteq\mathbb{R}^{d}$ and $\{z_{t}\}$ is a strictly stationary sequence of random variables with marginal $P$. We introduce a set of conditions for $f_{n,\theta}$ that induces a possibly localized counterpart of Kim and Pollard's (1990) cube root asymptotics. Their cube root asymptotics can be viewed as a special case of ours, where $f_{n,\theta}$ does not vary with $n$. Let $Pf=\int f dP$ for a function $f$, $|\cdot|$ be the Euclidean norm, and $\left\Vert \cdot\right\Vert _{2}$ be the $L_{2}(P)$-norm of a random variable. The class of criterion functions of interest is characterized as follows.

\begin{asmM}
For a sequence $\{h_{n}\}$ of positive numbers with $nh_{n}\rightarrow\infty$, $f_{n,\theta}$ satisfies the following conditions.\smallskip

\noindent(i) $h_{n}f_{n,\theta}$ is uniformly bounded, $\lim_{n\rightarrow\infty}Pf_{n,\theta}$ is uniquely
maximized at $\theta_{0}$, $Pf_{n,\theta}$ is twice continuously differentiable at $\theta_{0}$ for all $n$ large enough, and satisfies
\begin{equation}
P(f_{n,\theta}-f_{n,\theta_{0}})=\frac{1}{2}(\theta-\theta_{0})^{\prime}V(\theta-\theta_{0})+o(|\theta-\theta_{0}|^{2})+o((nh_{n})^{-2/3}),\label{eq:exp}
\end{equation}
for a negative definite matrix $V.$\smallskip

\noindent(ii) There exist positive constants $C$ and $C^{\prime}$ such that 
\[
|\theta_{1}-\theta_{2}|\leq Ch_{n}^{1/2}\left\Vert f_{n,\theta_{1}}-f_{n,\theta_{2}}\right\Vert _{2},
\]
for all $n$ large enough and $\theta_{1},\theta_{2}\in\{\Theta:|\theta-\theta_{0}|\leq C^{\prime}\}$.\smallskip

\noindent(iii) There exists a positive constant $C^{\prime\prime}$ such that 
\begin{equation}
P\sup_{\theta\in\Theta:|\theta-\theta^{\prime}|<\varepsilon}h_{n}|f_{n,\theta}-f_{n,\theta^{\prime}}|^{2}\leq C^{\prime\prime}\varepsilon,\label{eq:env}
\end{equation}
for all $n$ large enough, $\varepsilon>0$ small enough, and $\theta^{\prime}$ in a neighborhood of $\theta_{0}$.
\end{asmM}

$\{h_{n}\}$ is usually a sequence of bandwidths for localization. Although we are primarily interested in the case of $h_{n}\rightarrow0$, we do not exclude the case of $h_{n}=1$ which corresponds to the conventional cube root asymptotics in Kim and Pollard (1990). Also, we note that our conditions appear somewhat different from Kim and Pollard (1990). In fact, our conditions consist of directly verifiable moment conditions without resorting to the notion of empirical process theory such as uniform manageability. 

Assumption M (i) requires boundedness, point identification, and quadratic approximation of $Pf_{n,\theta}$. Boundedness of $h_{n}f_{n,\theta}$ is a major requirement but is satisfied for all examples in this paper and for Kim and Pollard (1990). In Section \ref{sec:set}, we relax the assumption of point identification. When the criterion function involves kernel smoothing for localization, it typically takes the form of a product of a bounded function and smoothing term $\frac{1}{h_{n}}K\left(\frac{x-c}{h_{n}}\right)$ (see (\ref{eq:PY}) and examples in Section \ref{sec:ex}).

Despite discontinuity of $f_{n,\theta}$, its population counterpart $Pf_{n,\theta}$ is smooth and approximated by a quadratic function as in (\ref{eq:exp}). This distinguishes our estimation problem from that of a change-point in a regression model, which also involves a discontinuous criterion function but the change-point estimator is super-consistent (e.g., Chan, 1993) unless the estimating equation is misspecified as in the split point estimator for decision trees (B\"{u}hlmann and Yu, 2002, and Banerjee and McKeague, 2007). 

Assumption M (ii) is used to relate the $L_{2}(P)$-norm for the criterion functions to the Euclidean norm for the parameters. This condition is implicit in Kim and Pollard (1990, Condition (v)) under independent observations and is often verified in the course of checking the expansion in (\ref{eq:exp}).

Assumption M (iii), an envelope condition for the class $\mathcal{F}_{n}=\{f_{n,\theta}-f_{n,\theta^{\prime}}:|\theta-\theta^{\prime}|\le\varepsilon\}$, plays a key role for cube root asymptotics. It should be noted that for the familiar square root asymptotics, the upper bound in (\ref{eq:env}) is of order $\varepsilon^{2}$ instead of $\varepsilon$. It is often the case that verifying the envelope condition for arbitrary $\theta^{\prime}$ in a neighborhood of $\theta_{0}$ is no more demanding than that for $\theta_{0}$.

In particular, Assumption M (iii) is used to guarantee an integrability condition on the metric entropy with bracketing for $\mathcal{F}_{n}$ in the $L_{2,\beta}$-norm so that the maximal inequality in Doukhan, Massart and Rio (1995, Theorem 3) can be applied to establish Lemma M below. On the other hand, Kim and Pollard (1990) used the concept of uniform manageability (Pollard, 1989) to control the size of $\mathcal{F}_{n}$ which is defined by the $\epsilon$-capacity, or metric entropy, by covering numbers. Generally the bracketing and covering numbers approaches are not directly comparable (see Section 2.5 of van der Vaart and Wellner, 1996, for example). It would be interesting to explore how the symmetrization argument combined with the suitable manageability concept can be applied in our setup.

We now study the asymptotic properties of the M-estimator which is precisely defined as a random variable $\hat{\theta}$ satisfying
\begin{equation}
\mathbb{P}_{n}f_{n,\hat{\theta}}\geq\sup_{\theta\in\Theta}\mathbb{P}_{n}f_{n,\theta}-o_{p}((nh_{n})^{-2/3}).\label{eq:M-est}
\end{equation}
The first step is to establish weak consistency $\hat{\theta}\overset{p}{\to}\theta_{0}$, which is rather standard and usually shown by establishing the uniform convergence of $\mathbb{P}_{n}f_{n,\theta}$. In this section we simply assume the consistency of $\hat{\theta}$. See the Supplement for some illustrations to show consistency.

The next step is to derive the convergence rate of $\hat{\theta}$. A key ingredient for this is to obtain the modulus of continuity of the empirical process $\{\mathbb{G}_{n}h_{n}^{1/2}(f_{n,\theta}-f_{n,\theta_{0}})\}$ by some maximum inequality, where $\mathbb{G}_{n}f=\sqrt{n}(\mathbb{P}_{n}f-Pf)$ for a function $f$. If $f_{n,\theta}$ does not vary with $n$ and $\{z_{t}\}$ is independent, several maximal inequalities are available in the literature (e.g., page 199 of Kim and Pollard, 1990). If $f_{n,\theta}$ varies with $n$ and $\{z_{t}\}$ is dependent, to the best of our knowledge, there is no maximal inequality which can be applied to the class of functions satisfying Assumption M. Our first task is to establish such an inequality.

To proceed, we now characterize the dependence structure of the data. Among several notions of dependence, this paper focuses on an absolutely regular process. See Doukhan, Massart and Rio (1995) for a discussion on empirical process theory of absolutely regular processes. Let $\mathcal{F}_{-\infty}^{0}$ and $\mathcal{F}_{m}^{\infty}$ be $\sigma$-fields of $\{\ldots,z_{t-1},z_{0}\}$ and $\{z_{m},z_{m+1},\ldots\}$, respectively. Define the $\beta$-mixing coefficient as $\beta_{m}=\frac{1}{2}\sup\sum_{(i,j)\in I\times J}|P\{A_{i}\cap B_{j}\}-P\{A_{i}\}P\{B_{j}\}|$, where the supremum is taken over all finite partitions $\{A_{i}\}_{i\in I}$ and $\{B_{j}\}_{j\in J}$, respectively $\mathcal{F}_{-\infty}^{0}$ and $\mathcal{F}_{m}^{\infty}$ measurable. Throughout the paper, we maintain the following assumption on $\{z_{t}\}$.

\begin{asmD} 
$\{z_{t}\}$ is a strictly stationary and absolutely regular process with $\beta$-mixing coefficients $\{\beta_{m}\}$ such that $\beta_{m}=O(\rho^{m})$ for some $0<\rho<1$.
\end{asmD}

This assumption obviously covers the case of independent observations. It also says the mixing coefficient $\beta_{m}$ should decay at an exponential rate.\footnote{Polynomial decays of $\beta_{m}$ are often associated with strong dependence and long memory type behaviors in sample statistics. See Chen, Hansen and Carrasco (2010) and references therein. In this case, asymptotic analysis for the M-estimator will be very different.} For example, various Markov, GARCH, and stochastic volatility models satisfy this assumption (Carrasco and Chen, 2002). See Section \ref{sub:D} below for further discussions.

Under this assumption, we obtain the following maximal inequality.

\begin{lemM} Under Assumptions M and D, there exist positive constants
$C$ and $C^{\prime}$ such that 
\[
P\sup_{\theta\in\Theta:|\theta-\theta_{0}|<\delta}|\mathbb{G}_{n}h_{n}^{1/2}(f_{n,\theta}-f_{n,\theta_{0}})|\leq C\delta^{1/2},
\]
for all $n$ large enough and $\delta\in[(nh_{n})^{-1/2},C^{\prime}]$.
\end{lemM}

This lemma provides a preliminary lemma to derive the convergence rate.

\begin{lem} \label{lem:KP4.1}
Under Assumptions M and D, for each $\varepsilon>0$, there exist random variables $\{R_{n}\}$ of order $O_{p}(1)$ and a positive constant $C$ such that 
\[
|\mathbb{P}_{n}(f_{n,\theta}-f_{n,\theta_{0}})-P(f_{n,\theta}-f_{n,\theta_{0}})|\leq\varepsilon|\theta-\theta_{0}|^{2}+(nh_{n})^{-2/3}R_{n}^{2},
\]
for all $\theta\in\{\Theta:(nh_{n})^{-1/3}\leq|\theta-\theta_{0}|\leq C\}$.
\end{lem}

We now derive the convergence rate of $\hat{\theta}$. Suppose $|\hat{\theta}-\theta_{0}|\geq(nh_{n})^{-1/3}$. Then by (\ref{eq:M-est}), Lemma \ref{lem:KP4.1}, and Assumption M (i), we can take a positive constant $c$ such that 
\begin{eqnarray*}
o_{p}((nh_{n})^{-2/3}) & \leq & \mathbb{P}_{n}(f_{n,\hat{\theta}}-f_{n,\theta_{0}})\\
 & \leq & P(f_{n,\hat{\theta}}-f_{n,\theta_{0}})+\varepsilon|\hat{\theta}-\theta_{0}|^{2}+(nh_{n})^{-2/3}R_{n}^{2}\\
 & \leq & (-c+\varepsilon)|\hat{\theta}-\theta_{0}|^{2}+o(|\hat{\theta}-\theta_{0}|^{2})+O_{p}((nh_{n})^{-2/3}),
\end{eqnarray*}
for each $\varepsilon>0$. Taking $\varepsilon$ small enough to satisfy $\varepsilon<c$ yields the convergence rate $\hat{\theta}-\theta_{0}=O_{p}((nh_{n})^{-1/3})$.

Given this, the final step is to establish the limiting distribution of $\hat{\theta}$. To this end, we apply a continuous mapping theorem of an argmax element (e.g., Theorem 2.7 of Kim and Pollard, 1990) and it is enough to show weak convergence of the normalized empirical process 
\[
Z_{n}(s)=n^{1/6}h_{n}^{2/3}\mathbb{G}_{n}(f_{n,\theta_{0}+s(nh_{n})^{-1/3}}-f_{n,\theta_{0}}),
\]
for $|s|\leq K$ with any $K>0$. Weak convergence of $Z_{n}$ may be characterized by its finite dimensional convergence and stochastic asymptotic equicontinuity (or tightness). If $f_{n,\theta}$ does not vary with $n$ and $\{z_{t}\}$ is independent as in Kim and Pollard (1990), a classical central limit theorem combined with the Cram\'{e}r-Wold device implies finite dimensional convergence, and a maximal inequality on a suitable class of functions guarantees stochastic asymptotic equicontinuity of the normalized empirical process. We adapt this approach to our local M-estimation problem with dependent observations.

Consider a function $\beta(\cdot)$ such that $\beta(t)=\beta_{[t]}$ if $t\geq1$ and $\beta(t)=1$ otherwise and denote its c\`{a}dl\`{a}g inverse by $\beta^{-1}(\cdot)$. Let $Q_{g}(u)$ be the inverse function of the tail probability function $x\mapsto P\{|g(z_{t})|>x\}$. For finite dimensional convergence, we employ Rio's (1997, Corollary 1) central limit theorem for $\alpha$-mixing arrays to our setup.

\begin{lemC}
Suppose Assumption D holds true, $Pg_{n}=0$, and 
\begin{equation}
\sup_{n\in\mathbb{N}}\int_{0}^{1}\beta^{-1}(u)Q_{g_{n}}(u)^{2}du<\infty.\label{eq:Linde}
\end{equation}
Then $\Sigma=\lim_{n\rightarrow\infty}\mathrm{Var}(\mathbb{G}_{n}g_{n})$
exists and $\mathbb{G}_{n}g_{n}\overset{d}{\to}N(0,\Sigma)$.
\end{lemC}

The finite dimensional convergence of $Z_{n}$ follows from Lemma C by setting $g_{n}$ as any finite dimensional projection of the process $\{g_{n,s}-Pg_{n,s}\}$ with
\begin{equation}
g_{n,s}=n^{1/6}h_{n}^{2/3}(f_{n,\theta_{0}+s(nh_{n})^{-1/3}}-f_{n,\theta_{0}}).\label{eq:gns}
\end{equation}
The requirement in (\ref{eq:Linde}) is the Lindeberg-type condition in Rio (1997, Corollary 1) and excludes polynomial decay of $\beta_{m}$. Note that for criterion functions satisfying Assumption M, the $(2+\delta)$-th moments $P|g_{n,s}|^{2+\delta}$ typically diverge because $g_{n,s}$ usually involves indicator functions. To verify (\ref{eq:Linde}), the following lemma is often useful.

\begin{lem} \label{lem:linde} Suppose Assumptions M and D hold true
and there is a positive constant $c$ such that
\begin{equation}
P\{|g_{n,s}|\geq c\}\leq c(nh_{n}^{-2})^{-1/3},\label{eq:suffL}
\end{equation}
for all $s$ and $n$ large enough. Then (\ref{eq:Linde}) holds true.
\end{lem}

In our examples, $g_{n,s}$ is zero or close to zero with high probability so that (\ref{eq:suffL}) is easily satisfied. See Section \ref{sec:ex} for illustrations.

We provide another maximal inequality that is useful to establish stochastic asymptotic equicontinuity of the process $Z_{n}$.

\begin{lemM'} Suppose Assumption D holds true. Consider a class of functions $\mathcal{G}_{n}=\{g_{n,s}:|s|\le K\}$ for some $K>0$ with envelope $G_{n}$. Suppose there is a positive constant $C$ such that
\begin{equation}
P\sup_{s:|s-s^{\prime}|<\varepsilon}|g_{n,s}-g_{n,s^{\prime}}|^{2}\leq C\varepsilon,\label{eq:Lpg}
\end{equation}
for all $n$ large enough, $|s^{\prime}|\le K$, and $\varepsilon>0$ small enough. Also, assume that there exist $0\leq\kappa<1/2$ and $C^{\prime}>0$ such that $G_{n}\leq C^{\prime}n^{\kappa}$ and $\left\Vert G_{n}\right\Vert _{2}\leq C^{\prime}$ for all $n$ large enough. Then for any $\sigma>0$, there exist $\delta>0$ and a positive integer $N_{\delta}$ such that 
\[
P\sup_{(s,s^{\prime}):|s-s^{\prime}|<\delta}|\mathbb{G}_{n}(g_{n,s}-g_{n,s^{\prime}})|\leq\sigma,
\]
for all $n\geq N_{\delta}$. \end{lemM'}

Stochastic asymptotic equicontinuity of $Z_{n}$ is implied from this lemma by setting $g_{n,s}$ as in (\ref{eq:gns}). Note that (\ref{eq:Lpg}) is satisfied by Assumption M (iii).\footnote{The upper bound in (\ref{eq:Lpg}) can be relaxed to $\varepsilon^{1/p}$ for $1\le p<\infty$. However, it is typically satisfied with $p=1$ for the examples we consider.} Compared to Lemma M used to derive the convergence rate of $\hat{\theta}$, Lemma M' is applied only to establish stochastic asymptotic equicontinuity of $Z_{n}$. Therefore, we do not need an exact decay rate on the right hand side of the maximal inequality.\footnote{In particular, $Z_{n}$ itself does not satisfy Assumption M (ii).}

By finite dimensional convergence and stochastic asymptotic equicontinuity of $Z_{n}$, its weak convergence is implied and the continuous mapping theorem for an argmax element (Theorem 2.7 of Kim and Pollard, 1990) yields the limiting distribution of $\hat{\theta}$. Define the covariance kernel
\[
H(s_{1},s_{2})=\lim_{n\to\infty}\sum_{t=-n}^{n}\mathrm{Cov}(g_{n,s_{1}}(z_{0}),g_{n,s_{2}}(z_{t})),
\]
if it exists. Throughout the paper, we use this notation for different choices of $g_{n,s}$. The main theorem of this section is presented as follows.

\begin{thm} \label{thm:cube1}
Suppose that Assumptions M and D hold, $\hat{\theta}$ defined in (\ref{eq:M-est}) converges in probability to $\theta_{0}\in\mathrm{int}\Theta$, and (\ref{eq:Linde}) holds with $g_{n,s}-Pg_{n,s}$ defined in (\ref{eq:gns}) for each $s$. Then 
\begin{equation}
(nh_{n})^{1/3}(\hat{\theta}-\theta_{0})\overset{d}{\rightarrow}\arg\max_{s\in\mathbb{R}^{d}}Z(s),\label{eq:dist}
\end{equation}
where $Z(s)$ is a Gaussian process with continuous sample paths, expected value $s^{\prime}Vs/2$, and covariance kernel $H(s_{1},s_{2})$.
\end{thm}

This theorem can be considered as an extension of the main theorem of Kim and Pollard (1990) to the cases where the criterion function can vary with the sample size and the observations can obey a dependent process. To the best of our knowledge, the (nonparametric) cube root convergence rate $(nh_{n})^{1/3}$ is new in the literature. It is interesting to note that similar to standard nonparametric estimation, $nh_{n}$ still plays the role of the ``effective sample size.''

\subsection{Nuisance parameters}

It is often the case that the criterion function contains some nuisance parameters, which can be estimated with rates faster than $(nh_{n})^{1/3}$. For the rest of this section, let $\hat{\theta}$ and $\tilde{\theta}$ satisfy
\begin{eqnarray*}
\mathbb{P}_{n}f_{n,\hat{\theta},\hat{\nu}} & \ge & \sup_{\theta\in\Theta}\mathbb{P}_{n}f_{n,\theta,\hat{\nu}}+o_{p}((nh_{n})^{-2/3}),\\
\mathbb{P}_{n}f_{n,\tilde{\theta},\nu_{0}} & \ge & \sup_{\theta\in\Theta}\mathbb{P}_{n}f_{n,\theta,\nu_{0}}+o_{p}((nh_{n})^{-2/3}),
\end{eqnarray*}
respectively, where $\nu_{0}$ is a vector of nuisance parameters
and $\hat{\nu}$ is its estimator satisfying $\hat{\nu}-\nu_{0}=o_{p}((nh_{n})^{-1/3})$.
Theorem \ref{thm:cube1} is extended as follows.

\begin{thm} \label{thm:cube1n} Suppose Assumption D holds true. Let $\{f_{n,\theta,\nu_{0}}:\theta\in\Theta\}$ satisfy Assumption M and $\{f_{n,\theta,\nu}:\theta\in\Theta,\nu\in\Lambda\}$ satisfy Assumption M (iii). Also assume that there exists a negative definite matrix $V_{1}$ such that 
\begin{eqnarray}
 & & P(f_{n,\theta,\nu}-f_{n,\theta_{0},\nu_{0}})  \label{eq:exp2n} \\
& = & \frac{1}{2}(\theta-\theta_{0})^{\prime}V_{1}(\theta-\theta_{0})+o(|\theta-\theta_{0}|^{2})+O(|\nu-\nu_{0}|^{2})+o((nh_{n})^{-2/3}), \nonumber
\end{eqnarray}
for all $\theta$ and $\nu$ in neighborhoods of $\theta_{0}$ and $\nu_{0}$, respectively. Then $\hat{\theta}=\tilde{\theta}+o_{p}((nh_{n})^{-1/3})$. Additionally, if (\ref{eq:Linde}) holds with $(g_{n,s}-Pg_{n,s})$ for each $s$ with $g_{n,s}$ being $n^{1/6}h_{n}^{2/3}(f_{n,\theta_{0}+s(nh_{n})^{-1/3},\nu_{0}}-f_{n,\theta_{0},\nu_{0}})$,
then 
\[
(nh_{n})^{1/3}(\hat{\theta}-\theta_{0})\overset{d}{\rightarrow}\arg\max_{s\in\mathbb{R}^{d}}Z(s),
\]
where $Z(s)$ is a Gaussian process with continuous sample paths, expected value $s^{\prime}V_{1}s/2$ and covariance kernel $H(s_{1},s_{2})$.
\end{thm}

A key step for the proof of this theorem is to confirm that the empirical process $\mathbb{G}_{n}f_{n,\theta,\nu_{0}+c(nh_{n})^{-1/3}}$ is well approximated by $\mathbb{G}_{n}f_{n,\theta,\nu_{0}}$ over $|\theta-\theta_{0}|\le\epsilon$ and $|c|\le\epsilon$ (see (A.10) in the Supplement). This is shown by applying Lemma M' with $g_{n,s}=n^{1/6}h_{n}^{2/3}(f_{n,\theta,\nu_{0}+c(nh_{n})^{-1/3}}-f_{n,\theta,\nu_{0}})$. Condition (\ref{eq:Lpg}) in Lemma M' demands more precise control on the size of the envelope for the class of $g_{n,s}$ than the comparable condition in Z-estimation with nuisance parameters (e.g., eq. (3) of van der Vaart and Wellner, 2007).

\subsection{Discussions}

\subsubsection{Inference}

Once we show that the M-estimator has a proper limiting distribution, Politis, Romano and Wolf (1999, Theorem 3.3.1) justify the use of subsampling to construct confidence intervals. Since Assumption D satisfies the requirement of their theorem, subsampling inference based on $s$ consecutive observations with $s/n\rightarrow\infty$ is asymptotically valid (in a pointwise sense explained below). See Politis, Romano and Wolf (1999, Section 3.6) for a discussion on data-dependent choices of $s$.

We note that this asymptotic validity of subsampling inference is in a pointwise sense rather than uniform. To be specific, suppose $\{z_{t}\}$ is an independent and identically distributed (iid) sample from the probability measure $P$ that belongs to a class of probability measures $\mathcal{P}$. Also denote the true parameters by $\theta_{0}(P)$ to make explicit the dependence on $P$. Based on Romano and Shaikh (2008), a confidence set $\mathcal{C}_{n}$ for $\theta_{0}(P)$ is called \emph{pointwise} valid in $(1-\alpha)$ level if
\[
\liminf_{n\rightarrow\infty}P\{\theta_{0}(P)\in\mathcal{C}_{n}\}\geq1-\alpha,
\]
for each $P\in\mathcal{P}$ and is called \emph{uniformly} valid in $(1-\alpha)$ level if
\[
\liminf_{n\rightarrow\infty}\inf_{P\in\mathcal{P}}P\{\theta_{0}(P)\in\mathcal{C}_{n}\}\geq1-\alpha.
\]
Our Theorems \ref{thm:cube1} and \ref{thm:cube1n} combined with Politis, Romano and Wolf (1999, Theorem 3.3.1) guarantee the pointwise validity of the subsampling confidence set based on quantiles of the subsample statistic $(sh_{s})^{1/3}(\hat{\theta}_{s}-\hat{\theta})$, where $\hat{\theta}_{s}$ and $\hat{\theta}$ are the M-estimators based on the subsample and full sample, respectively. Also a pointwise valid confidence interval for each element of $\theta_{0}(P)$ can be obtained in a similar manner. 

To investigate whether we can construct a uniformly valid confidence set in our setup, we assume that $\{z_{t}\}$ is iid and the distribution $J_{n}(\cdot,\theta,P)$ of $Q_{n}(\theta)=(nh_{n})^{2/3}\{\max_{\vartheta\in\Theta}\mathbb{P}_{n}f_{n,\vartheta}-\mathbb{P}_{n}f_{n,\theta}\}$ satisfies
\begin{equation}
\limsup_{n\rightarrow\infty}\sup_{\theta\in\Theta}\sup_{P\in\mathcal{P}:\theta=\theta_{0}(P)}\sup_{x\in\mathbb{R}}\{J_{s}(x,\theta,P)-J_{n}(x,\theta,P)\}\le0,\label{eq:J}
\end{equation}
Then, Romano and Shaikh (2008, Theorems 3.1 and 3.3) imply the uniform validity of the confidence set
\[
\mathcal{C}_{n}=\{\theta\in\Theta:Q_{n}(\theta)\le q_{s}(\theta,1-\alpha)\},
\]
over $\mathcal{P}$, where $q_{s}(\theta,1-\alpha)$ is the $(1-\alpha)$-th quantile of the distribution of the subsample statistic $Q_{s}(\theta)$. By inspection of Romano and Shaikh (2008), we can see that (\ref{eq:J}) is satisfied if $Q_{n}(\theta_{0}(P_{n}))$ converges in law to a unique continuous distribution for any sequence of $P_{n}\in\mathcal{P}$ yielding a row-wise iid triangular array. Our lemmas to obtain Theorem \ref{thm:cube1} can be readily extended to the array setting by restating Assumptions M and D and the additional conditions for Theorem \ref{thm:cube1} in the array setup. We note that computation of $\mathcal{C}_{n}$ may require an extensive numerical search over $\Theta$, where the quantile $q_{s}(\theta,1-\alpha)$ needs to be computed for each $\theta$.

The above uniformity result relies upon the general results in Romano and Shaikh (2008, Theorems 3.1 and 3.3) and there are at least three issues to be further considered. First, the iid assumption for the sample does not allow serial dependence as in Assumption D. To accommodate dependent data, the high level assumptions provided by Romano and Shaikh (2008, Theorems 3.1) for uniform validity should be modified. Second, it is not a trivial task to extend the results in Romano and Shaikh (2008) to inference on subvectors (or functions) of $\theta$ except for a conservative projection of $\mathcal{C}_{n}$ to a lower dimension. Third, a key result in Romano and Shaikh (2008, Theorems 3.1) holds for objects centered at the true parameter $\theta_{0}(P)$ instead of the estimator $\hat{\theta}$. Therefore, their result does not apply to the subsample statistic $(sh_{s})^{1/3}(\hat{\theta}_{s}-\hat{\theta})$. All of these issues require full length papers and are beyond the scope of this paper.

Another candidate to conduct inference based on the M-estimator is the bootstrap. However, even for independent observations, it is known that the naive nonparametric bootstrap is typically invalid under cube root asymptotics (Abrevaya and Huang, 2005, and Sen, Banerjee and Woodroofe, 2010).

\subsubsection{Generalization of Assumption D\label{sub:D}}

All the results in this section build upon Assumption D which requires $\{z_{t}\}$ to be strictly stationary and absolutely regular (or $\beta$-mixing) with exponentially decaying mixing coefficients. Assumption D is used for both the maximal inequality (Lemma M) and central limit theorem (Lemma C) which are building blocks to derive the asymptotic distribution of $\hat{\theta}$. It is of interest whether we can establish analogous results under more general setups, such as $\alpha$-mixing, by utilizing some recent developments in the empirical process theory for dependent data. For instance, Merlev\`{e}de, Peligrad and Rio (2009, 2011) obtained Bernstein type inequalities for $\alpha$-mixing processes and Baraud (2010) and Nickl and S\"{o}hl (2016, Section 3) explored the generic chaining argument by Talagrand (2005) for Markov chains. 

Since the central limit theorem in Rio (1997, Corollary 1) holds for $\alpha$-mixing arrays, we can modify Lemma C to accommodate $\alpha$-mixing processes. Thus, we focus on extending Lemma M, the maximal inequality. A crucial step for this extension is whether we can replace the key lemma in Doukhan, Massart and Rio (1995, Lemma 3), which leads to the maximal inequality for $\beta$-mixing processes (in eq. (A.6) of the Supplement) through a chaining argument. Specifically, consider a finite subclass $\mathcal{F}$ of bounded functions with cardinality $p\ge\exp(1)$. By a decoupling technique for $\beta$-mixing processes, Doukhan, Massart and Rio (1995, Lemma 3) showed that for positive constants $c$ and $c_{1}$, there exists a universal positive constant $C$ such that
\[
P\max_{f\in\mathcal{F}}|\mathbb{G}_{n}f|\le C\left(c\sqrt{\log p}+c_{1}q\frac{\log p}{\sqrt{n}}+c_{1}\beta_{q}\sqrt{n}\right),
\]
for all $q=1,\ldots,n$. Note that the above upper bound reduces to the first term $Cc\sqrt{\log p}$ for the iid case. By properly choosing $q$, the first term still dominates in the $\beta$-mixing case even if $\log p$ is close to $n$ so that Lemma M can be established. In contrast, the maximal inequality implied by Merlev\`{e}de, Peligrad and Rio (2009, (2.1) in Theorem 1) for $\alpha$-mixing would be written in the form of $C\left(c\sqrt{\log p}+c_{1}\log n\log\log n\frac{\log p}{\sqrt{n}}\right)$. Therefore, as $\log p$ becomes close to $n$, the second term will dominate. Since this order of cardinality $p$ (i.e., $\log p$ close to $n$) is required in the proof of Doukhan, Massart and Rio (1995, Theorem 2), the upper bound in Lemma M for $\alpha$-mixing processes would become larger.\footnote{Although a full investigation is beyond the scope of this paper, we conjecture that it is also the case for the generic chaining argument by Talagrand (2005). Indeed, eq. (1.9) on page 10 of Talagrand (2005) explains that generic chaining needs partitions of cardinality up to $2^{2^{n}}$. }

Another direction to extend our result is to accommodate general Markov chains that may not be covered by Assumption D. To this end, a chaining argument (see, Baraud, 2010, and Nickl and S\"{o}hl, 2016) based on Bernstein type inequalities for Markov chains (e.g., Adamczak, 2008, and Paulin, 2015) may yield an analog of Lemma M. Although this is an intriguing question, existing time series examples on cube root asymptotics mostly focus on mixing data (e.g., Polonik and Yao, 2000, and de Jong and Woutersen, 2011) and also typically involve additional conditioning or exogenous variables. Thus, we leave this extension for future work.

\section{Examples\label{sec:ex}}

We provide several examples to demonstrate the usefulness of the asymptotic theory in the last section. For the sake of space, we only sketch the arguments to verify the conditions to apply the theorems in Section \ref{sec:general}. Detailed verifications under primitive conditions are delegated to the Supplement.

\subsection{Dynamic panel discrete choice\label{subsec:HK}}

For a binary response $y_{it}$ and $k$-dimensional covariates $x_{it}$, consider a dynamic panel data model 
\begin{eqnarray*}
P\{y_{i0}=1|x_{i},\alpha_{i}\} & = & F_{0}(x_{i},\alpha_{i}),\\
P\{y_{it}=1|x_{i},\alpha_{i},y_{i0},\ldots,y_{it-1}\} & = & F(x_{it}^{\prime}\beta_{0}+\gamma_{0}y_{it-1}+\alpha_{i}),
\end{eqnarray*}
for $i=1,\ldots,n$ and $t=1,2,3$, where $\alpha_{i}$ is unobservable and both $F_{0}$ and $F$ are unknown. Honor\'{e} and Kyriazidou (2000) proposed the conditional maximum score estimator $(\hat{\beta},\hat{\gamma})$ that maximizes
\[
\sum_{i=1}^{n}K\left(\frac{x_{i2}-x_{i3}}{b_{n}}\right)(y_{i2}-y_{i1})\mathrm{sgn}\{(x_{i2}-x_{i1})^{\prime}\beta+(y_{i3}-y_{i0})\gamma\},
\]
where $K$ is a kernel function and $b_{n}$ is a bandwidth. Kernel smoothing is introduced to deal with the unknown link function $F$. Honor\'{e} and Kyriazidou (2000) obtained consistency of this estimator but the convergence rate and limiting distribution are unknown. Since the criterion function varies with the sample size through the bandwidth $b_{n}$, the cube root asymptotic theory of Kim and Pollard (1990) is not applicable here.

This open question can be addressed by Theorem \ref{thm:cube1}. Let $z=(z_{1}^{\prime},z_{2},z_{3}^{\prime})^{\prime}$ with $z_{1}=x_{2}-x_{3}$, $z_{2}=y_{2}-y_{1}$, and $z_{3}=((x_{2}-x_{1})^{\prime},y_{3}-y_{0})^{\prime}$. The above estimator for $\theta_{0}=(\beta_{0}^{\prime},\gamma_{0})^{\prime}$ can be written as an M-estimator using the criterion function
\begin{equation}
f_{n,\theta}(z)=e_{n}(z)(\mathbb{I}\{z_{3}^{\prime}\theta\geq0\}-\mathbb{I}\{z_{3}^{\prime}\theta_{0}\geq0\}),\label{eq:HK}
\end{equation}
where $e_{n}(z)=b_{n}^{-k}K(b_{n}^{-1}z_{1})z_{2}$. To apply Theorem \ref{thm:cube1}, it is enough to show that $f_{n,\theta}$ in (\ref{eq:HK}) satisfies Assumption M with $h_{n}=b_{n}^{k}$ and the condition in (\ref{eq:suffL}). Then the limiting distribution of Honor\'{e} and Kyriazidou's (2000) estimator is obtained as in (\ref{eq:dist}). 

Here we sketch the verification. See Section B.1 of the Supplement for detailed verifications and primitive conditions. For Assumption M (i), $\{h_{n}f_{n,\theta}\}$ is bounded for the bounded kernel $K$ and (\ref{eq:exp}) is obtained by a Taylor expansion combined with the argument in Kim and Pollard (1990, pp. 214-215). For Assumption M (ii), take any $\theta_{1}$ and $\theta_{2}$ and note that 
\begin{eqnarray*}
h_{n}^{1/2}\left\Vert f_{n,\theta_{1}}-f_{n,\theta_{2}}\right\Vert _{2} & = & \sqrt{ P\left\{ h_{n}E[e_{n}(z)^{2}|z_{3}]|\mathbb{I}\{z_{3}^{\prime}\theta_{1}\geq0\}-\mathbb{I}\{z_{3}^{\prime}\theta_{2}\geq0\}|\right\} }\\
 & \geq & c^{1/2}P\{z_{3}^{\prime}\theta_{1}\geq0>z_{3}^{\prime}\theta_{2}\text{ or }z_{3}^{\prime}\theta_{2}\geq0>z_{3}^{\prime}\theta_{1}\},
\end{eqnarray*}
for some $c>0$, where the inequality follows from $h_{n}E[e_{n}(z)^{2}|z_{3}]>c$ (by a change of variables and the condition on the density $z_{1}|(z_{2}\neq0,z_{3})$ being bounded away from zero) and Jensen's inequality. The right hand side is the probability for a pair of wedge shaped regions with an angle of order $|\theta_{1}-\theta_{2}|$. Thus, Assumption M (ii) is satisfied if the density of $z_{3}$ is bounded away from zero in a neighborhood of the origin. Assumption M (iii) can be verified in a similar way (by considering the upper bound instead). The Markov inequality and boundedness of the density imply (\ref{eq:suffL}).

\subsection{Random coefficient binary choice\label{subsec:rc}}

As a new statistical model which can be covered by our asymptotic theory, consider a regression model $y_{t}=x_{t}^{\prime}\theta(w_{t})+u_{t}$ with random coefficients. Suppose we observe $\{ \mathrm{sgn}(y_{t}), x_{t}, w_{t}\}$ and wish to estimate $\theta_{0}=\theta(c)$ at some given $c$.\footnote{Gautier and Kitamura (2013) studied identification and estimation of the random coefficient binary choice model, where $\theta_{t}=\theta(w_{t})$ is unobservable. Here we study the model where heterogeneity in the slope is caused by the observables $w_{t}$.} We propose a localized version of the maximum score estimator 
\begin{equation}
\hat{\theta}=\arg\max_{\theta\in S}\sum_{t=1}^{n}K\left(\frac{w_{t}-c}{b_{n}}\right)[\mathbb{I}\{y_{t}\ge0,x_{t}^{\prime}\theta\ge0\}+\mathbb{I}\{y_{t}<0,x_{t}^{\prime}\theta<0\}],\label{eq:rce}
\end{equation}
where $S$ is the surface of the unit sphere. Again, the cube root asymptotic theory of Kim and Pollard (1990) is not applicable due to the bandwidth.

Theorem \ref{thm:cube1} can be applied to obtain the limiting distribution of this estimator. Note that $\hat{\theta}$ in (\ref{eq:rce}) can be written as an M-estimator using the criterion function
\begin{equation}
f_{n,\theta}(x,w,u)=\frac{1}{h_{n}}K\left(\frac{w-c}{h_{n}^{1/k}}\right)h(x,u)[\mathbb{I}\{x^{\prime}\theta\ge0\}-\mathbb{I}\{x^{\prime}\theta_{0}\ge0\}],\label{eq:rcf}
\end{equation}
for $h_{n}=b_{n}^{k}$ and $h(x,u)=\mathbb{I}\{x^{\prime}\theta_{0}+u\ge0\}-\mathbb{I}\{x^{\prime}\theta_{0}+u<0\}$.
Once we check Assumption M and (\ref{eq:suffL}), Theorem \ref{thm:cube1} implies the limiting distribution. 

The verification is sketched as follows. See Section B.2 of the Supplement for detailed verifications and primitive conditions. Assumption M (i)-(ii) and (\ref{eq:suffL}) can be checked similarly as in Section \ref{subsec:HK}. Here we verify Assumption M (iii). By a change of variables and $h(x,u)^{2}=1$, there exists a positive constant $C^{\prime}$ such that
\begin{eqnarray*}
 &  & P\sup_{\theta\in\Theta:|\theta-\vartheta|<\varepsilon}h_{n}|f_{n,\theta}-f_{n,\vartheta}|^{2}\\
 & = & \int\int K(s)^{2}\sup_{\theta\in\Theta:|\theta-\vartheta|<\varepsilon}|[\mathbb{I}\{x^{\prime}\theta\ge0\}-\mathbb{I}\{x^{\prime}\vartheta\ge0\}]|^{2}p(x,c+sb_{n})dxds\\
 & \le & C^{\prime}E\left[\left.\sup_{\theta\in\Theta:|\theta-\vartheta|<\varepsilon}|[\mathbb{I}\{x^{\prime}\theta\ge0\}-\mathbb{I}\{x^{\prime}\vartheta\ge0\}]|^{2}\right|w=c\right],
\end{eqnarray*}
for all $\varepsilon>0$, $\vartheta$ in a neighborhood of $\theta_{0}$, and $n$ large enough, where $p$ is the joint density of $(x_{t},w_{t})$. Since the right hand side is the conditional probability for a pair of wedge shaped regions with an angle of order $\varepsilon$, Assumption M (iii) is guaranteed by some boundedness condition on the conditional density of $x_{t}$ given $w_{t}=c$.

\subsection{Minimum volume predictive region\label{subsec:mv}}

As an illustration of Theorem \ref{thm:cube1n}, we now consider the example in (\ref{eq:PY}), a simplified version of Polonik and Yao's (2000) minimum volume predictor. For notational convenience, assume $\theta_{0}=0$ and $\nu_{0}=1$.
By applying Lemma M', the convergence rate of the nuisance parameter estimator is obtained as $\hat{\nu}-1=O_{p}((nh_{n})^{-1/2}+h_{n}^{2})$ (see Section B.3 in the Supplement).

The criterion function for the maximization in (\ref{eq:PY}) can be written as 
\[
f_{n,\theta,\hat{\nu}}(y,x)=\frac{1}{h_{n}}K\left(\frac{x-c}{h_{n}}\right)[\mathbb{I}\{y\in[\theta-\hat{\nu},\theta+\hat{\nu}]\}-\mathbb{I}\{y\in[-\hat{\nu},\hat{\nu}]\}].
\]
We apply Theorem \ref{thm:cube1n} to obtain the convergence rate of $\hat{\theta}$. Details are provided in Section B.3 of the Supplement. Assumptions M for $f_{n,\theta,1}$ and M (iii) for $f_{n,\theta,\nu}$ are verified similarly as in Sections \ref{subsec:HK} and \ref{subsec:rc}. To check (\ref{eq:exp2n}), a Taylor expansion yields
\begin{eqnarray*}
 &  & P(f_{n,\theta,\nu}-f_{n,0,1})\\
& = & \frac{1}{2}V_{1}\theta^{2}+\{\dot{\gamma}_{y|x}(1|c)+\dot{\gamma}_{y|x}(-1|c)\}\gamma_{x}(c)\theta\nu+o(\theta^{2}+|\nu-1|^{2})+O(h_{n}^{2}),
\end{eqnarray*}
for $V_{1}=\{\dot{\gamma}_{y|x}(1|c)-\dot{\gamma}_{y|x}(-1|c)\}\gamma_{x}(c)$, where $\gamma$ and $\dot{\gamma}$ mean the density and its derivative, respectively. 

Therefore, Theorem \ref{thm:cube1n} implies $\hat{\theta}-\theta_{0}=O_{p}((nh_{n})^{-1/3}+h_{n})$, which confirms positively the conjecture of Polonik and Yao (2000, Remark 3b) on the exact convergence rate of $[\hat{\theta}\pm\hat{\nu}]$.

\subsection{Dynamic maximum score\label{subsec:ms}}

To illustrate the derivation of the covariance kernel $H$ in Theorem \ref{thm:cube1} for dependent data, we consider the maximum score estimator (Manski, 1975) for a regression model $y_{t}=x_{t}^{\prime}\theta_{0}+u_{t}$, that is 
\[
\hat{\theta}=\arg\max_{\theta\in S}\sum_{t=1}^{n}[\mathbb{I}\{y_{t}\ge0,x_{t}^{\prime}\theta\ge0\}+\mathbb{I}\{y_{t}<0,x_{t}^{\prime}\theta<0\}],
\]
where $S$ is the surface of the unit sphere. This estimator can be written as an M-estimator using the criterion function
\[
f_{\theta}(x,u)=h(x,u)[\mathbb{I}\{x^{\prime}\theta\ge0\}-\mathbb{I}\{x^{\prime}\theta_{0}\ge0\}],
\]
where $h(x,u)=\mathbb{I}\{x^{\prime}\theta_{0}+u\ge0\}-\mathbb{I}\{x^{\prime}\theta_{0}+u<0\}$. The conditions to apply Theorem \ref{thm:cube1} can be verified similarly as in the above examples (see Section B.4 of the Supplement). Here we focus on the derivation of the covariance kernel for the limiting distribution under Assumption D.

Let $q_{n,t}=f_{\theta_{0}+n^{-1/3}s_{1}}(x_{t},u_{t})-f_{\theta_{0}+n^{-1/3}s_{2}}(x_{t},u_{t})$. The covariance kernel is written as $H(s_{1},s_{2})=\frac{1}{2}\{L(s_{1},0)+L(0,s_{2})-L(s_{1},s_{2})\}$, where 
\[
L(s_{1},s_{2})=\lim_{n\to\infty}n^{1/3}\left\{ \mathrm{Var}(q_{n,t})+\sum_{m=1}^{\infty}\mathrm{Cov}(q_{n,t},q_{n,t+m})\right\} .
\]
The limit of $n^{1/3}\mathrm{Var}(q_{n,t})$ is given in Kim and Pollard (1990, p. 215). For the covariance we note that $q_{n,t}$ takes only three values: $-1$, $0$, or $1$. The definition of $\beta_{m}$ and Assumption D imply 
\[
|P\{q_{n,t}=j,q_{n,t+m}=k\}-P\{q_{n,t}=j\}P\{q_{n,t+m}=k\}|\leq n^{-2/3}\beta_{m},
\]
for all $n,m\ge1$ and $j,k=-1,0,1$. Thus, $\{q_{n,t}\}$ is a $\beta$-mixing array with mixing coefficients bounded by $n^{-2/3}\beta_{m}$. This in turn implies that $\{q_{n,t}\}$ is an $\alpha$-mixing array with mixing coefficients bounded by $2n^{-2/3}\beta_{m}$. By applying the $\alpha$-mixing inequality, the covariance is bounded as
\[
\mathrm{Cov}(q_{n,t},q_{n,t+m})\le Cn^{-2/3}\beta_{m}\left\Vert q_{n,t}\right\Vert _{p}^{2},
\]
for some $C>0$ and $p>2$. Note that 
\begin{eqnarray*}
 \left\Vert q_{n,t}\right\Vert _{p}^{2} & \leq & [P|\mathbb{I}\{x^{\prime}(\theta_{0}+s_{1}n^{-1/3})>0\}-\mathbb{I}\{x^{\prime}(\theta_{0}+s_{2}n^{-1/3})>0\}|]^{2/p}\\
 & = & O(n^{-2/(3p)}).
\end{eqnarray*}
Combining these results, we get $n^{1/3}\sum_{m=1}^{\infty}\mathrm{Cov}(q_{n,t},q_{n,t+m})\to0$ as $n\to\infty$. Therefore, the covariance kernel $H$ is the same as the independent case in Kim and Pollard (1990, p. 215).

\subsection{Other examples\label{subsec:other}}

In the Supplement, we present additional examples on the dynamic least median of squares estimator (Section B.5) and the monotone density estimator (Section B.6).

\section{Generalizations\label{sec:set}}

In this section, we consider two generalizations of the asymptotic theory in Section \ref{sec:general}. The first concerns data taking limited values such as interval-valued regressors and the second is to allow for localization to depend on the parameter values.

\subsection{Limited observations}

We consider the case where some of the variables take limited values. In particular, we relax the assumption of point identification of $\theta_{0}$ and study the case where the limiting criterion function is maximized at any element of a set $\Theta_{I}\subset\Theta$. The set $\Theta_{I}$ is called the identified set. In order to estimate $\Theta_{I}$, we consider a collection of approximate maximizers of the sample criterion function
\[
\hat{\Theta}=\{\theta\in\Theta:\max_{\theta\in\Theta}\mathbb{P}_{n}f_{n,\theta}-\mathbb{P}_{n}f_{n,\theta}\leq\hat{c}(nh_{n})^{-1/2}\},
\]
i.e., the level set based on $\mathbb{P}_{n}f_{n,\theta}$ from its maximum by a cutoff value $\hat{c}(nh_{n})^{-1/2}$. This section studies the convergence rate of $\hat{\Theta}$ to $\Theta_{I}$ under the Hausdorff distance defined below. We assume that $\Theta_{I}$ is convex. Then the projection $\pi_{\theta}=\arg\min_{\theta^{\prime}\in\Theta_{I}}|\theta^{\prime}-\theta|$ of $\theta$ on $\Theta_{I}$ is uniquely defined. To deal with the partially identified case, we modify Assumption M as follows.

\begin{asmS} For a sequence $\{h_{n}\}$ of positive numbers satisfying $nh_{n}\rightarrow\infty$, $f_{n,\theta}$ satisfies the following conditions.\smallskip

\noindent(i) $h_{n}f_{n,\theta}$ is uniformly bounded, $\lim_{n\rightarrow\infty}Pf_{n,\theta}$ is maximized at any $\theta$ in a bounded convex set $\Theta_{I}$, and there exist positive constants $c$ and $c^{\prime}$ such that 
\begin{equation}
P(f_{n,\pi_{\theta}}-f_{n,\theta})\ge c|\theta-\pi_{\theta}|^{2}+o(|\theta-\pi_{\theta}|^{2})+o((nh_{n})^{-2/3}),\label{eq:exps}
\end{equation}
for all $n$ large enough and all $\theta\in\{\Theta:0<|\theta-\pi_{\theta}|\leq c^{\prime}\}$.\smallskip

\noindent(ii) There exist positive constants $C$ and $C^{\prime}$ such that 
\[
|\theta-\pi_{\theta}|\leq Ch_{n}^{1/2}\left\Vert f_{n,\theta}-f_{n,\pi_{\theta}}\right\Vert _{2},
\]
for all $n$ large enough and all $\theta\in\{\Theta:0<|\theta-\pi_{\theta}|\leq C^{\prime}\}$.\smallskip

\noindent(iii) There exists a positive constant $C^{\prime\prime}$ such that 
\[
P\sup_{\theta\in\Theta:0<|\theta-\pi_{\theta}|<\varepsilon}h_{n}|f_{n,\theta}-f_{n,\pi_{\theta}}|^{2}\leq C^{\prime\prime}\varepsilon,
\]
for all $n$ large enough and all $\varepsilon>0$ small enough.
\end{asmS}

We allow $h_{n}=1$ for the case without a bandwidth in the criterion function. Similar comments to Assumption M apply. The main difference is that the conditions are imposed on the contrast $f_{n,\theta}-f_{n,\pi_{\theta}}$ using the projection $\pi_{\theta}$. Assumption S (i) contains boundedness and expansion conditions. The inequality in (\ref{eq:exps}) can be checked by a one-sided Taylor expansion using the directional derivative. Assumption S (ii) and (iii) play similar roles as Assumption M (ii) and (iii) and can be verified in a similar way.

We first establish the maximal inequality for the criterion functions satisfying Assumption S. Let $r_{n}=nh_{n}/\log(nh_{n})$.

\begin{lemMS} Under Assumptions D and S, there exist positive constants
$C$ and $C^{\prime}<1$ such that 
\[
P\sup_{\theta\in\Theta:0<|\theta-\pi_{\theta}|<\delta}|\mathbb{G}_{n}h_{n}^{1/2}(f_{n,\theta}-f_{n,\pi_{\theta}})|\leq C(\delta\log(1/\delta))^{1/2},
\]
for all $n$ large enough and $\delta\in[r_{n}^{-1/2},C^{\prime}]$.
\end{lemMS}

Compared to Lemma M, the additional log term on the right hand side is due to the fact that the supremum is taken over the $\delta$-tube (or manifold) instead of the $\delta$-ball, which increases the entropy. This maximal inequality is applied to obtain an analog of Lemma \ref{lem:KP4.1}.

\begin{lem} \label{lem:KP4.1-S} Under Assumptions D and S, for each $\varepsilon>0$, there exist random variables $\{R_{n}\}$ of order $O_{p}(1)$ and a positive constant $C$ such that
\[
|\mathbb{P}_{n}(f_{\theta}-f_{\pi_{\theta}})-P(f_{\theta}-f_{\pi_{\theta}})|\leq\varepsilon|\theta-\pi_{\theta}|^{2}+r_{n}^{-2/3}R_{n}^{2},
\]
for all $\theta\in\{\Theta:r_{n}^{-1/3}\leq|\theta-\pi_{\theta}|\leq C\}$.
\end{lem}

Let $\rho(A,B)=\sup_{a\in A}\inf_{b\in B}|a-b|$ and $H(A,B)=\max\{\rho(A,B),\rho(B,A)\}$ be the Hausdorff distance of sets $A,B\subset\mathbb{R}^{d}$. Based on these lemmas, the convergence rate of the set estimator $\hat{\Theta}$ is obtained as follows.

\begin{thm} \label{thm:cube1s}
Suppose that Assumptions D and S hold, $H(\hat{\Theta},\Theta_{I})\overset{p}{\to}0$, $\{h_{n}^{1/2}f_{n,\theta}:\theta\in\Theta_{I}\}$ is $P$-Donsker, and $\hat{c}=o_{p}((nh_{n})^{1/2})$. Then 
\[
\rho(\hat{\Theta},\Theta_{I})=O_{p}(\hat{c}^{1/2}(nh_{n})^{-1/4}+r_{n}^{-1/3}).
\]
Furthermore, if $\hat{c}\to\infty$, then $P\{\Theta_{I}\subset\hat{\Theta}\}\rightarrow1$
and
\[
H(\hat{\Theta},\Theta_{I})=O_{p}(\hat{c}^{1/2}(nh_{n})^{-1/4}).
\]
 \end{thm}

Note that $\rho$ is asymmetric in its arguments. In contrast to the convergence rate of $\rho(\hat{\Theta},\Theta_{I})$ obtained in the first part of this theorem, the second part says $P\{\Theta_{I}\subset\hat{\Theta}\}\rightarrow1$ (i.e., $\rho(\Theta_{I},\hat{\Theta})$ can converge to zero at an arbitrary rate) as far as $\hat{c}\to\infty$. For example, we may set $\hat{c}=\log(nh_{n})$. These results are combined to imply the convergence rate $H(\hat{\Theta},\Theta_{I})=O_{p}(\hat{c}^{1/2}(nh_{n})^{-1/4})$ under the Hausdorff distance. When $\hat{c}\to\infty$, the cube root term of order $r_{n}^{-1/3}$ in the rate of $\rho(\hat{\Theta},\Theta_{I})$ is dominated by the term of order $\hat{c}^{1/2}(nh_{n})^{-1/4}$.

We next consider the case where the criterion function contains nuisance parameters. In particular, we allow the dimension $k_{n}$ of the nuisance parameters $\nu$ to grow as the sample size increases. For instance the nuisance parameters might be coefficients in a sieve estimation procedure. It is important to allow the growing dimension of $\nu$ to cover Manski and Tamer's (2002) set estimator, where the criterion function contains some nonparametric estimate and its transform by the indicator. The rest of this subsection considers the set estimator
\[
\hat{\Theta}=\{\theta\in\Theta:\max_{\theta\in\Theta}\mathbb{P}_{n}f_{n,\theta,\hat{\nu}}-\mathbb{P}_{n}f_{n,\theta,\hat{\nu}}\leq\hat{c}(nh_{n})^{-1/2}\},
\]
with some preliminary estimator $\hat{\nu}$ and cutoff value $\hat{c}$.

Let $g_{n,s}=h_{n}^{1/2}(f_{n,\theta,\nu}-f_{n,\theta,\nu_{0}})$ with $s=(\theta^{\prime},\nu^{\prime})^{\prime}$ and consider $\mathcal{G}_{n}=\{g_{n,s}:|\theta-\pi_{\theta}|\leq K_{1},|\nu-\nu_{0}|\leq a_{n}K_{2}\}$ for some $K_{1},K_{2}>0$ with the envelope function $G_{n}=\sup_{\mathcal{G}_{n}}|g_{n,s}|$. The maximal inequality in Lemma MS is modified as follows.

\begin{lemMS'}
Suppose Assumption D holds true and there exists a positive constant $C$ such that 
\begin{eqnarray}
P\sup_{s:\theta\in\Theta,|\nu-\nu_{0}|\leq\varepsilon}|g_{n,s}|^{2} & \leq & C\sqrt{k_{n}}\varepsilon,\label{eq:Lpgs}\\
\sup_{s:\theta\in\Theta,|\nu-\nu_{0}|\leq\varepsilon}\{|\nu-\nu_{0}|-C\left\Vert g_{n,s}\right\Vert _{2}\} & \leq & 0,\label{eq:Lpgs1}
\end{eqnarray}
for all $n$ large enough and $\varepsilon$ small enough. Also assume that there exist $0\leq\kappa<1/4$ and $C^{\prime}>0$ such that $G_{n}\leq C^{\prime}n^{\kappa}$ and $\left\Vert G_{n}\right\Vert _{2}\leq C^{\prime}$ for all $n$ large enough. Then there exists $K_{3}>0$ such that
\[
P\sup_{g_{n,s}\in\mathcal{G}_{n}}|\mathbb{G}_{n}g_{n,s}|\leq K_{3}a_{n}^{1/2}k_{n}^{3/4}\sqrt{\log k_{n}a_{n}^{-1}},
\]
for all $n$ large enough. \end{lemMS'}

The increasing dimension $k_{n}$ of $\nu$ affects the upper bound via two routes. First, it increases the size of envelope by a factor of $\sqrt{k_{n}}$, which in turn increases the entropy of the space. Second, it also demands us to consider an inflated class of functions to apply the more fundamental maximal inequality by Doukhan, Massart and Rio (1995), which relies on the $\left\Vert \cdot\right\Vert _{2,\beta}$ norm. Note that the envelope condition in (\ref{eq:Lpgs}) allows for step functions containing some nonparametric estimates.

Based on this lemma, the convergence rate of the set estimator $\hat{\Theta}$ is characterized as follows.

\begin{thm} \label{thm:cube1sn}
Suppose Assumption D holds true. Let $\{f_{n,\theta,\nu_{0}}:\theta\in\Theta\}$ satisfy Assumption S and $\{h_{n}^{1/2}f_{n,\theta,\nu_{0}}:\theta\in\Theta_{I}\}$ be a $P$-Donsker class. Assume $\rho(\hat{\Theta},\Theta_{I})\overset{p}{\to}0$, $\hat{c}=o_{p}((nh_{n})^{1/2})$, $k_{n}\rightarrow\infty$, and $|\hat{\nu}-\nu_{0}|=o_{p}(a_{n})$ for some $\{a_{n}\}$ such that $h_{n}/a_{n}\rightarrow\infty$. Furthermore, there exist some $\varepsilon>0$ and neighborhoods $\{\theta\in\Theta:|\theta-\pi_{\theta}|<\varepsilon\}$ and $\{\nu:|\nu-\nu_{0}|\leq\varepsilon\}$, where $h_{n}^{1/2}(f_{n,\theta,\nu}-f_{n,\theta,\nu_{0}})$ satisfies (\ref{eq:Lpgs}) and (\ref{eq:Lpgs1}) and
\begin{equation}
P(f_{n,\theta,\nu}-f_{n,\pi_{\theta},\nu}-f_{n,\theta,\nu_{0}}+f_{n,\pi_{\theta},\nu_{0}}) =  o(|\theta-\pi_{\theta}|^{2})+O(|\nu-\nu_{0}|^{2}+r_{n}^{-2/3}).\label{eq:exp4}
\end{equation}
Then 
\begin{eqnarray}
\rho(\hat{\Theta},\Theta_{I}) & = & O_{p}(\hat{c}^{1/2}(nh_{n})^{-1/4}+r_{n}^{-1/3}) \label{eq:set-rate1} \\
& & +O_{p}((nh_{n}a_{n}^{-1})^{-1/4}(\log k_{n})^{1/2}) +o(a_{n}).\nonumber
\end{eqnarray}
Furthermore, if $\hat{c}\to\infty$, then $P\{\Theta_{I}\subset\hat{\Theta}\}\rightarrow1$ and
\begin{equation}
H(\hat{\Theta},\Theta_{I})=O_{p}(\hat{c}^{1/2}(nh_{n})^{-1/4}+(nh_{n})^{-1/4}a_{n}^{1/4}k_{n}^{3/8}\log^{1/4}n)+o(a_{n}).\label{eq:set-rate2}
\end{equation}
 \end{thm}

Compared to Theorem \ref{thm:cube1s}, we have two extra terms in (\ref{eq:set-rate2}) due to the (nonparametric) estimation of $\nu_{0}$. However, they can be shown to be dominated by the first term under standard conditions. Suppose that $k_{n}^{4}\log k_{n}/n\to0$ and the preliminary estimator $\hat{\nu}$ satisfies $|\hat{\nu}-\nu_{0}|=O_{p}(n^{-1/2}(k_{n}\log k_{n})^{1/2})$, which is often the case as in sieve estimation (see, e.g., Chen, 2007).\footnote{Alternatively $\nu_{0}$ can be estimated by some high-dimensional method (e.g. Belloni, Chen, Chernozhukov and Hansen, 2012) which also typically guarantees $a_{n}=o(n^{-1/4})$.} Then we can set $a_{n}=n^{-1/2}(k_{n}\log k_{n})^{1/2}$ so that $a_{n}^{1/4}k_{n}^{3/8}\to0$. Now by choosing $\hat{c}=\log n$, the first term in (\ref{eq:set-rate2}) dominates the other terms.

\subsubsection{Example: Binary choice with interval regressor}

As an illustration of partially identified models, we consider a binary choice model with an interval-valued regressor studied by Manski and Tamer (2002). Let $y=\mathbb{I}\{x^{\prime}\theta_{0}+w+u\ge0\}$ where $x$ is a vector of observable regressors, $w$ is an unobservable regressor, and $u$ is an unobservable error term satisfying $P\{u\le0|x,w\}=\alpha$ (we set $\alpha=.5$ to simplify the notation). Instead of $w$, we observe an interval $[w_{l},w_{u}]$ such that $P\{w_{l}\le w\le w_{u}\}=1$. Here we normalize the coefficient of $w$ to be one. In this setup, the parameter $\theta_{0}$ is partially identified and its identified set is written as (Manski and Tamer 2002, Proposition 2)
\[
\Theta_{I}=\{\theta\in\Theta:P\{x^{\prime}\theta+w_{u}\le0<x^{\prime}\theta_{0}+w_{l}\mbox{ or }x^{\prime}\theta_{0}+w_{u}\le0<x^{\prime}\theta+w_{l}\}=0\}.
\]
Let $\tilde{x}=(x^{\prime},w_{l},w_{u})^{\prime}$ and $q_{\hat{\nu}}(\tilde{x})$ be an estimator of $q_{\nu_{0}}(\tilde{x})=P\{y=1|\tilde{x}\}$ with the estimated parameters $\hat{\nu}$. By exploring the maximum score approach, Manski and Tamer (2002) developed the set estimator for $\Theta_{I}$
\begin{equation}
\hat{\Theta}=\{\theta\in\Theta:\max_{\theta\in\Theta}S_{n}(\theta)-S_{n}(\theta)\leq\epsilon_{n}\},\label{eq:MTE}
\end{equation}
where
\[
S_{n}(\theta)=\mathbb{P}_{n}(y-.5)[\mathbb{I}\{q_{\hat{\nu}}(\tilde{x})>.5\}\mathrm{sgn}(x^{\prime}\theta+w_{u})+\mathbb{I}\{q_{\hat{\nu}}(\tilde{x})\le.5\}\mathrm{sgn}(x^{\prime}\theta+w_{l})].
\]
Manski and Tamer (2002) established $H(\hat{\Theta},\Theta_{I})\overset{p}{\to}0$ by assuming that the cutoff value $\epsilon_{n}$ is bounded from below
by the (almost sure) uniform convergence rate of $S_{n}(\theta)$ to the limiting object. As Manski and Tamer (2002, Footnote 3) argued, characterization of this rate is a complex task because $S_{n}(\theta)$ is a step function and $\mathbb{I}\{q_{\hat{\nu}}(\tilde{x})>.5\}$ is a step function transform of the nonparametric estimate of $P\{y=1|\tilde{x}\}$. As such, it has been an open question. Obtaining the lower bound of $\epsilon_{n}$ is important because we wish to minimize the volume of the estimator $\hat{\Theta}$ without losing the asymptotic validity. By applying Theorem \ref{thm:cube1sn}, we can explicitly characterize the lower bound of $\epsilon_{n}$ and establish the convergence rate of $\hat{\Theta}$.

A little algebra shows that the set estimator in (\ref{eq:MTE}) is written as
\[
\hat{\Theta}=\{\theta\in\Theta:\max_{\theta\in\Theta}\mathbb{P}_{n}f_{\theta,\hat{\nu}}-\mathbb{P}_{n}f_{\theta,\hat{\nu}}\leq\hat{c}n^{-1/2}\},
\]
where $z=(x^{\prime},w,w_{l},w_{u},u)^{\prime}$, $h(x,w,u)=\mathbb{I}\{x^{\prime}\theta_{0}+w+u\ge0\}-\mathbb{I}\{x^{\prime}\theta_{0}+w+u<0\}$,
and 
\[
f_{\theta,\nu}(z)=h(x,w,u)[\mathbb{I}\{x^{\prime}\theta+w_{u}\ge0,q_{\nu}(\tilde{x})>.5\}-\mathbb{I}\{x^{\prime}\theta+w_{l}<0,q_{\nu}(\tilde{x})\le.5\}].
\]
To apply Theorem \ref{thm:cube1sn}, we check Assumption S with $h_{n}=1$. See Section B.7 of the Supplement for details. Here we illustrate the verifications of (\ref{eq:Lpgs}) and (\ref{eq:exp4}). Let $I_{\nu}(\tilde{x})=\mathbb{I}\{q_{\nu}(\tilde{x})>.5\ge q_{\nu_{0}}(\tilde{x})\mbox{ or }q_{\nu}(\tilde{x})\le.5<q_{\nu_{0}}(\tilde{x})\}$. Note that
\[
|f_{\theta,\nu}-f_{\theta,\nu_{0}}|^{2}\leq\mathbb{I}\{x^{\prime}\theta\ge-w_{u}\ge x^{\prime}\pi_{\theta}\mbox{ or }x^{\prime}\theta<-w_{u}<x^{\prime}\pi_{\theta}\}I_{\nu}(\tilde{x})\leq I_{\nu}(\tilde{x}).
\]
Also, we have
\[
P\sup_{\nu\in\Lambda:|\nu-\nu_{0}|<\varepsilon}I_{\nu}(\tilde{x}) \leq CP\sup_{\nu\in\Lambda:|\nu-\nu_{0}|<\varepsilon}|q_{\nu}(\tilde{x})-q_{\nu_{0}}(\tilde{x})|\le C\sqrt{k_{n}}\varepsilon,
\]
for some $C>0$, where the first inequality holds under boundedness of the conditional density of $q_{\nu_{0}}(\tilde{x})$ and the second under smoothness of $q_{\nu}$. This verifies (\ref{eq:Lpgs}). Also, (\ref{eq:Lpgs1}) is verified in the same manner as Assumption S (ii).

For (\ref{eq:exp4}) note that 
\begin{eqnarray}
 &  & |P(f_{\theta,\nu}-f_{\theta,\nu_{0}})-P(f_{\pi_{\theta},\nu}-f_{\pi_{\theta},\nu_{0}})| \label{eq:MT1} \\
 & \le & P\mathbb{I}\{x^{\prime}\theta\ge-w_{u}\ge x^{\prime}\pi_{\theta}\mbox{ or }x^{\prime}\theta<-w_{u}<x^{\prime}\pi_{\theta}\}I_{\nu}(\tilde{x})\nonumber \\
 &  & +P\mathbb{I}\{x^{\prime}\theta\ge-w_{l}\ge x^{\prime}\pi_{\theta}\mbox{ or }x^{\prime}\theta<-w_{l}<x^{\prime}\pi_{\theta}\}I_{\nu}(\tilde{x}),\nonumber
\end{eqnarray}
for each $\theta\in\{\Theta:|\theta-\pi_{\theta}|<\varepsilon\}$ and $\nu$ in a neighborhood of $\nu_{0}$. For the first term of (\ref{eq:MT1}), the law of iterated expectation and an expansion of $q_{\nu}(\tilde{x})$ around $\nu_{0}$ imply 
\begin{eqnarray*}
 &  & P\mathbb{I}\{x^{\prime}\theta\ge-w_{u}\ge x^{\prime}\pi_{\theta}\mbox{ or }x^{\prime}\theta<-w_{u}<x^{\prime}\pi_{\theta}\}I_{\nu}(\tilde{x})\\
 & \le & P\mathbb{I}\{x^{\prime}\theta\ge-w_{u}\ge x^{\prime}\pi_{\theta}\mbox{ or }x^{\prime}\theta<-w_{u}<x^{\prime}\pi_{\theta}\}A(w_{u},x)|v-\nu_{0}|,
\end{eqnarray*}
for some bounded function $A$. The second term of (\ref{eq:MT1}) is bounded in the same manner. Therefore, $|P(f_{\theta,\nu}-f_{\theta,\nu_{0}})-P(f_{\pi_{\theta},\nu}-f_{\theta,\nu_{0}})|=O(|\theta-\pi_{\theta}||v-\nu_{0}|)$ and (\ref{eq:exp4}) is verified. By applying Theorem \ref{thm:cube1sn}, we can conclude that the convergence rate of Manski and Tamer's (2002) set estimator $\hat{\Theta}$ in (\ref{eq:MTE}) is characterized by (\ref{eq:set-rate1}) and (\ref{eq:set-rate2}).

Compared to Manski and Tamer (2002), we provide a sharper lower bound on the the tuning parameter $\epsilon_{n},$ which is $\hat{c}n^{-1/2}$ with $\hat{c}\rightarrow\infty$. For example, if we set $\hat{c}=\log n$, the convergence rate becomes $H(\hat{\Theta},\Theta_{I})=O_{p}(n^{-1/4}(\log n)^{1/2})$. We basically verify the high level assumption of Chernozhukov, Hong and Tamer (2007, Condition C.2) in the cube root context. However, we mention that in the above setup, the criterion function contains nuisance parameters with increasing dimension and the result in Chernozhukov, Hong and Tamer (2007) does not apply directly.

Furthermore, our result enables us to construct a confidence set by subsampling as described in Chernozhukov, Hong and Tamer (2007). Specifically, the maximal inequality in Lemma MS' and the assumption that $\{h_{n}^{1/2}f_{n,\theta,\nu_{0}}:\theta\in\Theta_{I}\}$ is $P$-Donsker are sufficient to satisfy their Conditions C.4 and C.5.

\subsection{Parameter-dependent localization}

We now consider a setup where localization of the criterion function depends on the parameter values. A leading example is mode estimation. Chernoff (1964) studied asymptotic properties of the mode estimator that maximizes $(nh)^{-1}\sum_{t=1}^{n}\mathbb{I}\{|y_{t}-\beta|\le h\}$ with respect to $\beta$ for some fixed $h$. Lee (1989) extended this estimator to regression models, established its consistency, and conjectured the cube root convergence rate. To estimate $\beta$ consistently for a broader family of distributions, however, we need to treat $h$ as a bandwidth parameter and let $h\to0$ as in Yao, Lindsay and Li (2012) for example. 

This parameter-dependent localization alters Assumption M (iii) because it increases the size (in terms of the $L_{2}(P)$-norm) of the envelope of the class $\{h^{-1}(\mathbb{I}\{|y_{t}-\beta|\le h\}-\mathbb{I}\{|y_{t}-\beta_{0}|\le h\}):\mbox{ }|\beta-\beta_{0}|\leq\varepsilon\}$. More precisely, we replace Assumption M (iii) with the following.

\begin{asmM} (iii') There exists a positive constant $C^{\prime\prime}$ such that 
\[
P\sup_{\theta\in\Theta:|\theta-\theta^{\prime}|<\varepsilon}h_{n}^{2}|f_{n,\theta}-f_{n,\theta^{\prime}}|^{2}\leq C^{\prime\prime}\varepsilon,
\]
for all $n$ large enough, $\varepsilon>0$ small enough, and $\theta^{\prime}$ in a neighborhood of $\theta_{0}$. \end{asmM}

Under this assumption, Lemma M in Section \ref{sec:general} is modified
as follows.

\begin{lemM1} Under Assumption M (i), (ii), and (iii'), there exist positive constants $C$ and $C^{\prime}$ such that 
\[
P\sup_{\theta\in\Theta:|\theta-\theta_{0}|<\delta}|\mathbb{G}_{n}h_{n}^{1/2}(f_{n,\theta}-f_{n,\theta_{0}})|\leq Ch_{n}^{-1/2}\delta^{1/2},
\]
for all $n$ large enough and $\delta\in[(nh_{n}^{2})^{-1/2},C^{\prime}]$.
\end{lemM1}

Parameter dependency arises in different contexts and may lead to different types of non-standard distributions. For instance, the maximum likelihood estimator for $\mathrm{Uniform}[0,\theta]$ yields super consistency (see Hirano and Porter, 2003, for a general discussion). This contrast is similar to the difference between estimation of a change point in regression analysis and mode regression. 

Once we have obtained this lemma, the remaining steps are similar to those in Section \ref{sec:general} by replacing ``$h_{n}$'' with ``$h_{n}^{2}$''. Here we present the result without nuisance parameters $\nu$ for the sake of expositional simplicity. 

\begin{thm} \label{thm:cube1n-1}
Suppose that Assumptions D and M (i), (ii), and (iii') hold. Also suppose (\ref{eq:Linde}) holds with $(g_{n,s}-Pg_{n,s})$ for each $s$, where $g_{n,s}=n^{1/6}h_{n}^{4/3}(f_{n,\theta_{0}+s(nh_{n}^{2})^{-1/3}}-f_{n,\theta_{0}})$. Then 
\begin{equation}
(nh_{n}^{2})^{1/3}(\hat{\theta}-\theta_{0})\overset{d}{\rightarrow}\arg\max_{s\in\mathbb{R}^{d}}Z(s),\label{eq:dist1n-1}
\end{equation}
where $Z(s)$ is a Gaussian process with continuous sample paths, expected value $s^{\prime}Vs/2$ and covariance kernel $H(s_{1},s_{2})$.
\end{thm}

\subsubsection{Example: Hough transform estimator\label{subsec:mode}}

In the statistics literature on the computer vision algorithm, Goldenshluger and Zeevi (2004) investigated the so-called Hough transform estimator for regression models
\begin{equation}
\hat{\beta}=\arg\max_{\beta\in B}\sum_{t=1}^{n}\mathbb{I}\{|y_{t}-x_{t}^{\prime}\beta|\le h|x_{t}|\},\label{eq:hough}
\end{equation}
where $B$ is some parameter space, $x_{t}=(1,\tilde{x}_{t})^{\prime}$ for a scalar $\tilde{x}_{t}$, and $h$ is a fixed tuning constant. Goldenshluger and Zeevi (2004) derived the cube root asymptotics for $\hat{\beta}$ with fixed $h$ and discussed carefully the practical choice of $h$. However, for this estimator, $h$ plays the role of the bandwidth and the analysis for the case of $h_{n}\to0$ is a substantial open question (see pp. 1915-6 of Goldenshluger and Zeevi, 2004). Here we study the asymptotic property of $\hat{\beta}$ in (\ref{eq:hough}) with $h=h_{n}\to0$. The estimators by Chernoff (1964) and Lee (1989) with varying $h$ can be analyzed in the same manner.

Let $z=(x',u)'$. Note that $\hat{\theta}=\hat{\beta}-\beta_{0}$ is written as an M-estimator using the criterion function
\[
f_{n,\theta}(z)=h_{n}^{-1}\mathbb{I}\{|u-x^{\prime}\theta|\le h_{n}|x|\}.
\]
To apply Theorem \ref{thm:cube1n-1}, we need to verify that $f_{n,\theta}$ satisfies Assumption M (i), (ii), and (iii'). Here we focus on showing M (iii') while other details are found in Section B.8 of the Supplement. Observe that
\begin{eqnarray*}
 & & P\sup_{\theta\in\Theta:|\theta-\vartheta|<\varepsilon}h_{n}^{2}|f_{n,\theta}-f_{n,\vartheta}|^{2}\\
 & \le & P\sup_{\theta\in\Theta:|\theta-\vartheta|<\varepsilon}\mathbb{I}\{|u-x^{\prime}\vartheta|\le h_{n}|x|,\mbox{ }|u-x^{\prime}\theta|>h_{n}|x|\}\\
 &  & +P\sup_{\theta\in\Theta:|\theta-\vartheta|<\varepsilon}\mathbb{I}\{|u-x^{\prime}\theta|\le h_{n}|x|,\mbox{ }|u-x^{\prime}\vartheta|>h_{n}|x|\},
\end{eqnarray*}
for all $\vartheta$ in a neighborhood of $0$. Since the same argument applies to the second term, we focus on the first term (say, $T$). If $\varepsilon\le2h_{n}$, an expansion around $\varepsilon=0$
implies 
\[
T\le P\{(h_{n}-\varepsilon)|x|\le u\le h_{n}|x|\}=P\gamma(h_{n}|x|)|x|\varepsilon+o(\varepsilon),
\]
assuming independence between $u$ and $x$. Also, if $\varepsilon>2h_{n}$, an expansion around $h_{n}=0$ implies 
\[
T\le P\{-h_{n}|x|\le u\le h_{n}|x|\}\le P\gamma(0)|x|\varepsilon+o(h_{n}).
\]
Therefore, Assumption M (iii') is satisfied.

Finally, the covariance kernel is obtained in a similar way as in Section \ref{subsec:ms}. Let $r_{n}=(nh_{n}^{2})^{1/3}$. The covariance kernel is written by $H(s_{1},s_{2})=\frac{1}{2}\{L(s_{1},0)+L(0,s_{2})-L(s_{1},s_{2})\}$, where $L(s_{1},s_{2})=\lim_{n\to\infty}\mathrm{Var}(r_{n}^{2}\mathbb{P}_{n}g_{n,t})$ with $g_{n,t}=f_{n,s_{1}/r_{n}}-f_{n,s_{2}/r_{n}}$. An expansion implies $n^{-1}\mathrm{Var}(r_{n}^{2}g_{n,t})\rightarrow2\gamma(0)P|x^{\prime}(s_{1}-s_{2})|$, where $\gamma$ is the density of $u$. We can also see that the covariance term is negligible. Therefore, by Theorem \ref{thm:cube1n-1}, the limiting distribution of the Hough transform estimator with diminishing bandwidth is obtained as in (\ref{eq:dist1n-1}) with $V=\ddot{\gamma}(0)P(|x|xx^{\prime})$ and $H(s_{1},s_{2})=2\gamma(0)P|x^{\prime}(s_{1}-s_{2})|$.

\section{Conclusion}

This paper develops general asymptotic theory, which encompasses a wide class of non-regular M-estimation problems. Many of these problems have been left without a proper inference method for a long time. It is worthwhile to emphasize that our theory validates inference based on subsampling for this important class of estimators, including construction of confidence sets for set-valued parameters in Manski and Tamer's (2002) binary choice model with an interval regressor. An interesting line of future research is to develop valid bootstrap methods for these estimators. Naive applications of standard bootstrap resampling lead to inconsistent inference as shown by Abrevaya and Huang (2005) and Sen, Banerjee and Woodroofe (2010) among others.

\newpage
\appendix

\section*{Supplement to ``Local M-estimation with Discontinuous Criterion for
	Dependent and Incomplete Observations''}

	Section A presents the proofs of Lemmas and Theorems in the paper.
	In Section B, we provide primitive conditions and detailed verifications
	for the examples in Sections 3 and 4. Also Section B contains additional
	examples, which are omitted from the paper for brevity.

\section{Proofs of theorems and lemmas}

\subsection{Notation}

We employ the same notation in the paper. Recall that $Q_{g}(u)$
is the inverse function of the tail probability function $x\mapsto P\{|g(z_{t})|>x\}$
and that $\{\beta_{m}\}$ is the $\beta$-mixing coefficients used
in Assumption D. Let $\beta(\cdot)$ be a function such that $\beta(t)=\beta_{[t]}$
if $t\geq1$ and $\beta(t)=1$ otherwise and $\beta^{-1}(\cdot)$
be the càdlàg inverse of $\beta(\cdot)$. The $L_{2,\beta}(P)$-norm
is defined as 
\begin{equation}
	\left\Vert g\right\Vert _{2,\beta}=\sqrt{\int_{0}^{1}\beta^{-1}(u)Q_{g}(u)^{2}du}.\label{pf:L2b}
\end{equation}

\subsection{Proof of Lemma M}

Pick any $C^{\prime}>0$ and then pick any $n$ satisfying $(nh_{n})^{-1/2}\leq C^{\prime}$
and any $\delta\in[(nh_{n})^{-1/2},C^{\prime}]$. Throughout the proof,
positive constants $C_{j}$ ($j=1,2,\ldots$) are independent of $n$
and $\delta$.

First, we introduce some notation. Consider the sets defined by different
norms: 
\begin{eqnarray*}
	\mathcal{G}_{n,\delta}^{1} & = & \left\{ h_{n}^{1/2}(f_{n,\theta}-f_{n,\theta_{0}}):|\theta-\theta_{0}|<\delta\mbox{ for }\theta\in\Theta\right\} ,\\
	\mathcal{G}_{n,\delta}^{2} & = & \left\{ h_{n}^{1/2}(f_{n,\theta}-f_{n,\theta_{0}}):\left\Vert h_{n}^{1/2}(f_{n,\theta}-f_{n,\theta_{0}})\right\Vert _{2}<\delta\mbox{ for }\theta\in\Theta\right\} ,\\
	\mathcal{G}_{n,\delta}^{\beta} & = & \left\{ h_{n}^{1/2}(f_{n,\theta}-f_{n,\theta_{0}}):\left\Vert h_{n}^{1/2}(f_{n,\theta}-f_{n,\theta_{0}})\right\Vert _{2,\beta}<\delta\mbox{ for }\theta\in\Theta\right\} .
\end{eqnarray*}
For any $g\in\mathcal{G}_{n,\delta}^{1}$, $g$ is bounded (by Assumption
M (i)) and so is $Q_{g}$. Thus we can always find a function $\hat{g}$
such that $\left\Vert g\right\Vert _{2}^{2}\leq\left\Vert \hat{g}\right\Vert _{2}^{2}\leq2\left\Vert g\right\Vert _{2}^{2}$
and 
\begin{equation}
	Q_{g}(u)\le Q_{\hat{g}}(u)=\sum_{j=1}^{m}a_{j}\mathbb{I}\{(j-1)/m\leq u<j/m\},\label{pf:Qghat}
\end{equation}
for some positive integer $m$ and sequence of positive constants
$\{a_{j}\}$.

Next, we derive the set inclusion relationships 
\begin{equation}
	\mathcal{G}_{n,\delta}^{\beta}\subset\mathcal{G}_{n,\delta}^{2}\subset\mathcal{G}_{n,C_{1}\delta}^{1},\qquad\mathcal{G}_{n,\delta}^{1}\subset\mathcal{G}_{n,C_{2}\delta^{1/2}}^{\beta},\label{pf:set}
\end{equation}
for some positive constants $C_{1}$ and $C_{2}$. The relation $\mathcal{G}_{n,\delta}^{\beta}\subset\mathcal{G}_{n,\delta}^{2}$
follows from $\left\Vert \cdot\right\Vert _{2}\le\left\Vert \cdot\right\Vert _{2,\beta}$
(Doukhan, Massart and Rio, 1995, Lemma 1). The relation $\mathcal{G}_{n,\delta}^{2}\subset\mathcal{G}_{n,C_{1}\delta}^{1}$
follows from Assumption M (ii). Pick any $g\in\mathcal{G}_{\delta}^{1}$.
The relation $\mathcal{G}_{n,\delta}^{1}\subset\mathcal{G}_{n,C_{2}\delta^{1/2}}^{\beta}$
is obtained by 
\begin{eqnarray}
	\left\Vert g\right\Vert _{2,\beta}^{2} & \leq & \sum_{j=1}^{m}a_{j}^{2}\left\{ \int_{(j-1)/m}^{j/m}\beta^{-1}(u)du\right\} \le\left\{ m\int_{0}^{1/m}\beta^{-1}(u)du\right\} \int_{0}^{1}Q_{\hat{g}}(u)^{2}du\nonumber \\
	& \le & \left\{ \sup_{0<a\leq1}a\int_{0}^{1/a}\beta^{-1}(u)du\right\} 2\left\Vert g\right\Vert _{2}^{2}\le C_{2}^{2}\delta,\label{pf:g}
\end{eqnarray}
for some positive constant $C_{2}$, where the first inequality follows
from $Q_{g}\leq Q_{\hat{g}}$, the second inequality follows from
monotonicity of $\beta^{-1}(u)$ and $\int_{0}^{1}Q_{\hat{g}}(u)^{2}du=\frac{1}{m}\sum_{j=1}^{m}a_{j}^{2}$,
the third inequality follows by $\int_{0}^{1}Q_{\hat{g}}(u)^{2}du=\left\Vert \hat{g}\right\Vert _{2}^{2}\leq2\left\Vert g\right\Vert _{2}^{2}$,
and the last inequality follows from $\sup_{0<a\leq1}a\int_{0}^{1/a}\beta^{-1}(u)du<\infty$
(by Assumption D) and Assumption M (iii).

Third, based on (\ref{pf:set}), we derive some relationships for
the bracketing numbers. Let $N_{[]}(\nu,\mathcal{G},\left\Vert \cdot\right\Vert )$
be the bracketing number for a class of functions $\mathcal{G}$ with
radius $\nu>0$ and norm $\left\Vert \cdot\right\Vert $. Note that
\[
N_{[]}(\nu,\mathcal{G}_{n,\delta}^{\beta},\left\Vert \cdot\right\Vert _{2,\beta})\leq N_{[]}(\nu,\mathcal{G}_{n,C_{1}\delta}^{1},\left\Vert \cdot\right\Vert _{2})\leq C_{3}\left(\frac{\delta}{\nu}\right)^{2d},
\]
for some positive constant $C_{3}$, where the first inequality follows
from $\mathcal{G}_{n,\delta}^{\beta}\subset\mathcal{G}_{n,C_{1}\delta}^{1}$
(by (\ref{pf:set})) and $\left\Vert \cdot\right\Vert _{2}\le\left\Vert \cdot\right\Vert _{2,\beta}$,
and the second inequality follows from the argument to derive Andrews
(1993, eq. (4.7)) based on Assumption M (iii). Therefore, by the indefinite
integral formula $\int\log xdx=\mbox{const.}+x(\log x-1)$, there
exists a positive constant $C_{4}$ such that 
\begin{equation}
	\varphi_{n}(\delta)=\int_{0}^{\delta}\sqrt{\log N_{[]}(\nu,\mathcal{G}_{n,\delta}^{\beta},\left\Vert \cdot\right\Vert _{2,\beta})}d\nu\leq C_{4}\delta.\label{pf:phi}
\end{equation}

Finally, based on the entropy condition (\ref{pf:phi}), we apply
the maximal inequality of Doukhan, Massart and Rio (1995, Theorem
3), i.e., there exists a positive constant $C_{5}$ such that 
\begin{equation}
	P\sup_{g\in\mathcal{G}_{n,\delta}^{\beta}}|\mathbb{G}_{n}g|\leq C_{5}[1+\delta^{-1}q_{G_{n,\delta}}(\min\{1,v_{n}(\delta)\})]\varphi_{n}(\delta),\label{pf:max}
\end{equation}
where $q_{G_{n,\delta}}(v)=\sup_{u\leq v}Q_{G_{n,\delta}}(u)\sqrt{\int_{0}^{u}\beta^{-1}(\tilde{u})d\tilde{u}}$
with the envelope function $G_{n,\delta}$ of $\mathcal{G}_{n,\delta}^{\beta}$,
and $v_{n}(\delta)$ is the unique solution of 
\[
\frac{v_{n}(\delta)^{2}}{\int_{0}^{v_{n}(\delta)}\beta^{-1}(\tilde{u})d\tilde{u}}=\frac{\varphi_{n}(\delta)^{2}}{n\delta^{2}}.
\]
Since $\varphi_{n}(\delta)\leq C_{4}\delta$ from (\ref{pf:phi}),
it holds $v_{n}(\delta)\leq C_{5}n^{-1}$ for some positive constant
$C_{5}$. Now take some $n_{0}$ such that $v_{n_{0}}(\delta)\leq1$
and then pick again any $n\geq n_{0}$ and $\delta\in[(nh_{n})^{-1/2},C^{\prime}]$.
We have 
\begin{equation}
	q_{G_{n,\delta}}(\min\{1,v_{n}(\delta)\})\leq C_{6}Q_{G_{n,\delta}}(v_{n}(\delta))\sqrt{v_{n}(\delta)}\leq C_{7}(nh_{n})^{-1/2},\label{pf:qG}
\end{equation}
for some positive constants $C_{6}$ and $C_{7}$. Therefore, combining
(\ref{pf:phi})-(\ref{pf:qG}), we obtain 
\begin{equation}
	P\sup_{g\in\mathcal{G}_{n,C_{2}\delta^{1/2}}^{\beta}}|\mathbb{G}_{n}g|\leq C_{8}\delta^{1/2},\label{pf:delta}
\end{equation}
for some positive constant $C_{8}$. The conclusion follows from the
second relation in (\ref{pf:set}).

\subsection{Proof of Lemma 1}

Pick any $C>0$ and $\varepsilon>0$. Define 
\begin{eqnarray*}
	A_{n} & = & \{\theta\in\Theta:(nh_{n})^{-1/3}\leq|\theta-\theta_{0}|\leq C\},\\
	R_{n}^{2} & = & (nh_{n})^{2/3}\sup_{\theta\in A_{n}}\{|\mathbb{P}_{n}(f_{n,\theta}-f_{n,\theta_{0}})-P(f_{n,\theta}-f_{n,\theta_{0}})|-\varepsilon|\theta-\theta_{0}|^{2}\}.
\end{eqnarray*}
It is enough to show $R_{n}=O_{p}(1)$. Let
\[
A_{n,j}=\{\theta\in\Theta:(j-1)(nh_{n})^{-1/3}\leq|\theta-\theta_{0}|<j(nh_{n})^{-1/3}\}.
\]
There exists a positive constant $C^{\prime}$ such that 
\begin{eqnarray*}
	&  & P\{R_{n}>m\}\\
	& \le & P\left\{ |\mathbb{P}_{n}(f_{n,\theta}-f_{n,\theta_{0}})-P(f_{n,\theta}-f_{n,\theta_{0}})|>\varepsilon|\theta-\theta_{0}|^{2}+(nh_{n})^{-2/3}m^{2}\quad\text{for some }\theta\in A_{n}\right\} \\
	& \leq & \sum_{j=1}^{\infty}P\left\{ (nh_{n})^{2/3}|\mathbb{P}_{n}(f_{n,\theta}-f_{n,\theta_{0}})-P(f_{n,\theta}-f_{n,\theta_{0}})|>\varepsilon(j-1)^{2}+m^{2}\quad\text{for some }\theta\in A_{n,j}\right\} \\
	& \leq & \sum_{j=1}^{\infty}\frac{C^{\prime}\sqrt{j}}{\varepsilon(j-1)^{2}+m^{2}},
\end{eqnarray*}
for all $m>0$, where the last inequality is due to the Markov inequality
and Lemma M. Since the above sum is finite for all $m>0$, the conclusion
follows.

\subsection{Proof of Lemma C}

First of all, any $\beta$-mixing process is $\alpha$-mixing with
the mixing coefficient $\alpha_{m}\le\beta_{m}/2$. Thus it is sufficient
to check Conditions (a) and (b) of Rio (1997, Corollary 1). Under
eq. (5) of the paper, Condition (a) is verified by Rio (1997, Proposition
1), which guarantees $\mathrm{Var}(\mathbb{G}_{n}g_{n})\le\int_{0}^{1}\beta^{-1}(u)Q_{g_{n}}(u)^{2}du$
for all $n$. Since $\mathrm{Var}(\mathbb{G}_{n}g_{n})$ is bounded
(by eq. (5)) and $\{z_{t}\}$ is strictly stationary under Assumption
D, Condition (b) of Rio (1997, Corollary 1) can be written as 
\[
\int_{0}^{1}\beta^{-1}(u)Q_{g_{n}}(u)^{2}\inf_{n}\{n^{-1/2}\beta^{-1}(u)Q_{g_{n}}(u),1\}du\rightarrow0,
\]
as $n\to\infty$. Pick any $u\in(0,1)$. Since $\beta^{-1}(u)Q_{g_{n}}(u)^{2}$
is non-increasing in $u\in(0,1)$, the condition in eq. (5) implies
$\beta^{-1}(u)Q_{g_{n}}(u)^{2}<C<\infty$ for all $n$. Therefore,
for each $u\in(0,1)$, it holds $n^{-1/2}\beta^{-1}(u)Q_{g_{n}}(u)\to0$
as $n\to\infty$. Then the dominated convergence theorem based on
eq. (5) implies Condition (b).

\subsection{Proof of Lemma 2}

By Assumption M (i), it holds $|g_{n,s}|\le2C(nh_{n}^{-2})^{1/6}$
for all $n$ and $s$, which implies $Q_{g_{n,s}-Pg_{n,s}}(u)^{2}\le16C^{2}(nh_{n}^{-2})^{1/3}$
for all $n$, $s$, and $u\in(0,1)$. By the condition of this lemma,
it holds $Q_{g_{n,s}}(u)\le c$ for all $n$ large enough and $u>c(nh_{n}^{-2})^{-1/3}$.
By the triangle inequality and the definition of $Q_{g}$, 
\[
P\{|g_{n,s}-Pg_{n,s}|\ge Q_{g_{n,s}}(u)+|Pg_{n,s}|\}\le P\{|g_{n,s}|\ge Q_{g_{n,s}}(u)\}=P\{|g_{n,s}-Pg_{n,s}|>Q_{g_{n,s}-Pg_{n,s}}(u)\},
\]
which implies $Q_{g_{n,s}-Pg_{n,s}}(u)\le Q_{g_{n,s}}(u)+|Pg_{n,s}|$.
Thus, for all $n$ large enough, $s$, and $u>c(nh_{n}^{-2})^{-1/3}$,
it holds 
\[
Q_{g_{n,s}-Pg_{n,s}}(u)^{2}\le c^{2}+|Pg_{n,s}|^{2}+2c|Pg_{n,s}|.
\]
Combining these bounds, eq. (5) of the paper is verified as 
\begin{eqnarray*}
	&  & \int_{0}^{1}\beta^{-1}(u)Q_{g_{n,s}-Pg_{n,s}}(u)^{2}du\\
	& \leq & 16C^{2}(nh_{n}^{-2})^{1/3}\int_{0}^{c(nh_{n}^{-2})^{-1/3}}\beta^{-1}(u)du+\{c^{2}+(Pg_{n,s})^{2}+2c|Pg_{n,s}|\}\int_{c(nh_{n}^{-2})^{-1/3}}^{1}\beta^{-1}(u)du\\
	& < & \infty,
\end{eqnarray*}
for all $n$ large enough, where the second inequality follows by
Assumptions M (i) and D (which guarantees $\sup_{n}(nh_{n}^{-2})^{1/3}\int_{0}^{c(nh_{n}^{-2})^{-1/3}}\beta^{-1}(u)du<\infty$
and $\int_{0}^{1}\beta^{-1}(u)du<\infty$).

\subsection{Proof of Lemma M'}

Pick any $K>0$ and $\sigma>0$. Let $g_{n,s,s^{\prime}}=g_{n,s}-g_{n,s^{\prime}}$,
\begin{eqnarray*}
	\mathcal{G}_{n}^{K} & = & \{g_{n,s,s^{\prime}}:|s|\leq K,|s^{\prime}|\leq K\},\\
	\mathcal{G}_{n,\delta}^{1} & = & \{g_{n,s,s^{\prime}}\in\mathcal{G}_{n}^{K}:|s-s^{\prime}|<\delta\},\\
	\mathcal{G}_{n,\delta}^{\beta} & = & \{g_{n,s,s^{\prime}}\in\mathcal{G}_{n}^{K}:\left\Vert g_{n,s,s^{\prime}}\right\Vert _{2,\beta}<\delta\}.
\end{eqnarray*}
Since $g_{n,s}$ satisfies the condition in eq. (8) of the paper,
there exists a positive constant $C_{1}$ such that $\mathcal{G}_{n,\delta}^{1}\subset\{g_{n,s,s^{\prime}}\in\mathcal{G}_{n}^{K}:\left\Vert g_{n,s,s^{\prime}}\right\Vert _{2}<C_{1}\delta^{1/2}\}$
for all $n$ large enough and all $\delta>0$ small enough. Also,
by the same argument to derive (\ref{pf:g}), there exists a positive
constant $C_{2}$ such that $\left\Vert g_{n,s,s^{\prime}}\right\Vert _{2,\beta}\leq C_{2}\left\Vert g_{n,s,s^{\prime}}\right\Vert _{2}$
for all $n$ large enough, $|s|\leq K$, and $|s^{\prime}|\leq K$.
The constant $C_{2}$ depends only on the mixing sequence $\{\beta_{m}\}$.
Combining these results, we obtain
\begin{equation}
	\mathcal{G}_{n,\delta}^{1}\subset\mathcal{G}_{n,C_{1}C_{2}\delta^{1/2}}^{\beta},\label{pf:g1}
\end{equation}
for all $n$ large enough and all $\delta>0$ small enough.

Also note that the bracketing numbers satisfy 
\[
N_{[]}(\nu,\mathcal{G}_{n,\delta}^{\beta},\left\Vert \cdot\right\Vert _{2,\beta})\leq N_{[]}(\nu,\mathcal{G}_{n}^{K},\left\Vert \cdot\right\Vert _{2})\leq C_{3}\nu^{-d/2},
\]
where the first inequality follows from $\mathcal{G}_{n,\delta}^{\beta}\subset\mathcal{G}_{n}^{K}$
(by the definitions) and $\left\Vert \cdot\right\Vert _{2}\le\left\Vert \cdot\right\Vert _{2,\beta}$
(Doukhan, Massart and Rio, 1995, Lemma 1), and the second inequality
follows from the argument to derive Andrews (1993, eq. (4.7)) based
on eq. (8) of the paper. Thus, there is a function $\varphi(\eta)$
such that $\varphi(\eta)\rightarrow0$ as $\eta\rightarrow0$ and
\[
\varphi_{n}(\eta)=\int_{0}^{\eta}\sqrt{\log N_{[]}(\nu,\mathcal{G}_{n,\eta}^{\beta},\left\Vert \cdot\right\Vert _{2,\beta})}d\nu\leq\varphi(\eta),
\]
for all $n$ large enough and all $\eta>0$ small enough.

Based on the above entropy condition, we can apply the maximal inequality
of Doukhan, Massart and Rio (1995, Theorem 3), i.e., there exists
a positive constant $C_{3}$ depending only on the mixing sequence
$\{\beta_{m}\}$ such that 
\[
P\sup_{g\in\mathcal{G}_{n,\eta}^{\beta}}|\mathbb{G}_{n}g|\leq C_{4}[1+\eta^{-1}q_{G_{n}}(\min\{1,v_{n}(\eta)\})]\varphi(\eta),
\]
for all $n$ large enough and all $\eta>0$ small enough, where $q_{2G_{n}}(v)=\sup_{u\leq v}Q_{2G_{n}}(u)\sqrt{\int_{0}^{u}\beta^{-1}(\tilde{u})d\tilde{u}}$
with the envelope function $2G_{n}$ of $\mathcal{G}_{n,\eta}^{\beta}$
(note: by the definition of $\mathcal{G}_{n,\eta}^{\beta}$, the envelope
$2G_{n}$ does not depend on $\eta$) and $v_{n}(\eta)$ is the unique
solution of 
\[
\frac{v_{n}(\eta)^{2}}{\int_{0}^{v_{n}(\eta)}\beta^{-1}(\tilde{u})d\tilde{u}}=\frac{\varphi_{n}^{2}(\eta)}{n\eta^{2}}.
\]
Now pick any $\eta>0$ small enough so that $2C_{4}\varphi(\eta)<\sigma$.
Since $\varphi_{n}(\eta)\leq\varphi(\eta)$, there is a positive constant
$C_{5}$ such that $v_{n}(\eta)\leq C_{5}\frac{\varphi(\eta)}{n\eta^{2}}$
for all $n$ large enough and $\eta>0$ small enough. Since $G_{n}\leq C^{\prime}n^{\kappa}$
by the definition of $\mathcal{G}_{n,\eta}^{\beta}$, there exist
$C_{6}>0$ and $0<\kappa<1/2$ such that
\[
q_{2G_{n}}(\min\{1,v_{n}(\eta)\})\leq C_{6}\sqrt{\varphi(\eta)}\eta^{-1}n^{\kappa-1/2},
\]
for all $n$ large enough. Therefore, by setting $\eta=C_{1}C_{2}\delta^{1/2}$,
we obtain 
\[
P\sup_{g\in\mathcal{G}_{n,C_{1}C_{2}\delta^{1/2}}^{\beta}}|\mathbb{G}_{n}g|\leq\sigma,
\]
for all $n$ large enough. The conclusion follows by (\ref{pf:g1}).

\subsection{Proof of Theorem 1}

As discussed in the paper, Lemma 1 yields the convergence rate of
the M-estimator $\hat{\theta}$. This enables us to consider the centered
and normalized process $Z_{n}(s)$, which can be defined on arbitrary
compact parameter space. Based on finite dimensional convergence and
stochastic asymptotic equicontinuity of $Z_{n}$ shown by Lemmas C
and M', respectively, we establish weak convergence of $Z_{n}$.
Then a continuous mapping theorem of an argmax element (Kim and Pollard,
1990, Theorem 2.7) yields the limiting distribution of $\hat{\theta}$.

\subsection{Proof of Theorem 2}

To ease notation, let $\theta_{0}=\nu_{0}=0$. First, we show that
$\hat{\theta}=O_{p}((nh_{n})^{-1/3})$. Since $\{f_{n,\theta,\nu}\}$
satisfies Assumption M (iii), we can apply Lemma M' with $g_{n,s}=n^{1/6}h_{n}^{2/3}(f_{n,\theta,c(nh_{n})^{-1/3}}-f_{n,\theta,0})$
for $s=(\theta^{\prime},c^{\prime})^{\prime}$, which implies 
\begin{equation}
	\sup_{|\theta|\leq\epsilon,|c|\leq\epsilon}n^{1/6}h_{n}^{2/3}\mathbb{G}_{n}(f_{n,\theta,c(nh_{n})^{-1/3}}-f_{n,\theta,0})=O_{p}(1),\label{eq:1n1}
\end{equation}
for all $\epsilon>0$. Also from eq. (10) of the paper and $\hat{\nu}=o_{p}((nh_{n})^{-1/3})$,
we have 
\begin{equation}
	P(f_{n,\theta,\hat{\nu}}-f_{n,\theta,0})-P(f_{n,0,\hat{\nu}}-f_{n,0,0})\leq2\epsilon|\theta|^{2}+O_{p}((nh_{n})^{-2/3}),\label{eq:1n2}
\end{equation}
for all $\theta$ in a neighborhood of $\theta_{0}$ and all $\epsilon>0$.
Combining (\ref{eq:1n1}), (\ref{eq:1n2}), and Lemma 1, 
\begin{eqnarray*}
	\mathbb{P}_{n}(f_{n,\theta,\hat{\nu}}-f_{n,0,\hat{\nu}}) & = & n^{-1/2}\{\mathbb{G}_{n}(f_{n,\theta,\hat{\nu}}-f_{n,\theta,0})+\mathbb{G}_{n}(f_{n,\theta,0}-f_{n,0,0})-\mathbb{G}_{n}(f_{n,0,\hat{\nu}}-f_{n,0,0})\}\\
	&  & +P(f_{n,\theta,\hat{\nu}}-f_{n,\theta,0})+P(f_{n,\theta,0}-f_{n,0,0})-P(f_{n,0,\hat{\nu}}-f_{n,0,0})\\
	& \leq & P(f_{n,\theta,0}-f_{n,0,0})+2\epsilon|\theta|^{2}+O_{p}((nh_{n})^{-2/3})\\
	& \leq & \frac{1}{2}\theta^{\prime}V_{1}\theta+3\epsilon|\theta|^{2}+O_{p}((nh_{n})^{-2/3}),
\end{eqnarray*}
for all $\theta$ in a neighborhood of $\theta_{0}$ and all $\epsilon>0$,
where the last inequality follows from eq. (10) of the paper. From
$\mathbb{P}_{n}(f_{n,\hat{\theta},\hat{\nu}}-f_{n,0,\hat{\nu}})\geq o_{p}((nh_{n})^{-2/3})$,
negative definiteness of $V_{1}$, and $\hat{\nu}=o_{p}((nh_{n})^{-1/3})$,
we can find $c>0$ such that 
\[
o_{p}((nh_{n})^{-2/3})\leq-c|\hat{\theta}|^{2}+|\hat{\theta}|o_{p}((nh_{n})^{-1/3})+O_{p}((nh_{n})^{-2/3}),
\]
which implies $|\hat{\theta}|=O_{p}((nh_{n})^{-1/3})$.

Next, we show that $\hat{\theta}-\tilde{\theta}=o_{p}((nh_{n})^{-1/3})$.
By reparametrization, 
\[
(nh_{n})^{1/3}\hat{\theta}=\arg\max_{s}(nh_{n})^{2/3}[(\mathbb{P}_{n}-P)(f_{n,s(nh_{n})^{-1/3},\hat{\nu}}-f_{n,0,\hat{\nu}})+P(f_{n,s(nh_{n})^{-1/3},\hat{\nu}}-f_{n,0,\hat{\nu}})]+o_{p}(1).
\]
By Lemma M' (replacing $\theta$ with $(\theta,\nu)$) and $\hat{\nu}=o_{p}((nh_{n})^{-1/3})$,
\[
(\mathbb{P}_{n}-P)(f_{n,s(nh_{n})^{-1/3},\hat{\nu}}-f_{n,0,0})-(\mathbb{P}_{n}-P)(f_{n,s(nh_{n})^{-1/3},0}-f_{n,0,0})=o_{p}((nh_{n})^{-2/3}),
\]
uniformly in $s$. Also eq. (10) of the paper implies
\[
P(f_{n,s(nh_{n})^{-1/3},\hat{\nu}}-f_{n,0,\hat{\nu}})-P(f_{n,s(nh_{n})^{-1/3},0}-f_{n,0,0})=o_{p}((nh_{n})^{-2/3}),
\]
uniformly in $s$. Given $\hat{\theta}-\tilde{\theta}=o_{p}((nh_{n})^{-1/3})$,
an application of Theorem 1 to the class $\{f_{n,\theta,\nu_{0}}:\theta\in\Theta\}$
implies the limiting distribution of $\hat{\theta}$.

\subsection{Proof of Lemma MS}

First, we introduce some notation. Let 
\begin{align*}
	\mathcal{G}_{n,\delta}^{\beta} & =\left\{ h_{n}^{1/2}(f_{n,\theta}-f_{n,\pi_{\theta}}):\left\Vert h_{n}^{1/2}(f_{n,\theta}-f_{n,\pi_{\theta}})\right\Vert _{2,\beta}<\delta\mbox{ for }\theta\in\Theta\right\} ,\\
	\mathcal{G}_{n,\delta}^{1} & =\left\{ h_{n}^{1/2}(f_{n,\theta}-f_{n,\pi_{\theta}}):|\theta-\pi_{\theta}|<\delta\mbox{ for }\theta\in\Theta\right\} ,\\
	\mathcal{G}_{n,\delta}^{2} & =\left\{ h_{n}^{1/2}(f_{n,\theta}-f_{n,\pi_{\theta}}):\left\Vert h_{n}^{1/2}(f_{n,\theta}-f_{n,\pi_{\theta}})\right\Vert _{2}<\delta\mbox{ for }\theta\in\Theta\right\} .
\end{align*}
For any $g\in\mathcal{G}_{n,\delta}^{1}$, $g$ is bounded (Assumption
S (i)) and so is $Q_{g}$. Thus we can always find a function $\hat{g}$
such that $\left\Vert g\right\Vert _{2}^{2}\leq\left\Vert \hat{g}\right\Vert _{2}^{2}\leq2\left\Vert g\right\Vert _{2}^{2}$
and 
\[
Q_{g}(u)\le Q_{\hat{g}}(u)=\sum_{j=1}^{m}a_{j}\mathbb{I}\{(j-1)/m\leq u<j/m\},
\]
for some positive integer $m$ and sequence of positive constants
$\{a_{j}\}$. Let $r_{n}=nh_{n}/\log(nh_{n})$. Pick any $C^{\prime}>0$
and then pick any $n$ satisfying $r_{n}^{-1/2}\leq C^{\prime}$ and
any $\delta\in[(r_{n}^{-1/2},C^{\prime}]$. Throughout the proof,
positive constants $C_{j}$ ($j=1,2,\ldots$) are independent of $n$
and $\delta$.

Next, we derive some set inclusion relationships. Let $M=\frac{1}{2}\sup_{0<x\leq1}x^{-1}\int_{0}^{x}\beta^{-1}(u)du$.
For any $g\in\mathcal{G}_{\delta}^{1}$, it holds 
\begin{align}
	\left\Vert g\right\Vert _{2}^{2} & \leq\int_{0}^{1}\beta^{-1}(u)Q_{g}(u)^{2}du\leq\frac{1}{m}\sum_{j=1}^{m}a_{j}^{2}\left\{ m\int_{(j-1)/m}^{j/m}\beta^{-1}(u)du\right\} \nonumber \\
	& \leq\left\{ m\int_{0}^{1/m}\beta^{-1}(u)du\right\} \int_{0}^{1}Q_{\hat{g}}(u)^{2}du\nonumber \\
	& \leq M\left\Vert g\right\Vert _{2}^{2},\label{pf:g-s}
\end{align}
where the first inequality is due to Doukhan, Massart and Rio (1995,
Lemma 1), the second inequality follows from $Q_{g}\leq Q_{\hat{g}}$,
the third inequality follows from monotonicity of $\beta^{-1}(u)$,
and the last inequality follows by $\left\Vert \hat{g}\right\Vert _{2}^{2}\leq2\left\Vert g\right\Vert _{2}^{2}$.
Therefore, 
\begin{equation}
	\left\Vert f_{n,\theta}-f_{n,\pi_{\theta}}\right\Vert _{2}\leq\left\Vert f_{n,\theta}-f_{n,\pi_{\theta}}\right\Vert _{2,\beta}\leq M^{1/2}\left\Vert f_{n,\theta}-f_{n,\pi_{\theta}}\right\Vert _{2},\label{pf:norm-s}
\end{equation}
for each $\theta\in\{\theta\in\Theta:|\theta-\pi_{\theta}|<\delta\}$,
where the first inequality follows from Doukhan, Massart and Rio (1995,
Lemma 1) and the second inequality follows from (\ref{pf:g-s}). Based
on this, we can deduce the inclusion relationships: there exist positive
constants $C_{1}$ and $C_{2}$ such that 
\begin{equation}
	\mathcal{G}_{n,\delta}^{1}\subset\mathcal{G}_{n,C_{1}\delta^{1/2}}^{2}\subset\mathcal{G}_{n,M^{1/2}C_{1}\delta^{1/2}}^{\beta},\qquad\mathcal{G}_{n,\delta}^{\beta}\subset\mathcal{G}_{n,\delta}^{2}\subset\mathcal{G}_{n,C_{2}\delta}^{1},\label{pf:set-s}
\end{equation}
where the relation $\mathcal{G}_{n,\delta}^{1}\subset\mathcal{G}_{n,C_{1}\delta^{1/2}}^{2}$
follows from Assumption S (iii) and the relation $\mathcal{G}_{n,\delta}^{2}\subset\mathcal{G}_{n,C_{2}\delta}^{1}$
follows from Assumption S (ii).

Third, based on the above set inclusion relationships, we derive some
relationships for the bracketing numbers. Let $N_{[]}(\nu,\mathcal{G},\left\Vert \cdot\right\Vert )$
be the bracketing number for a class of functions $\mathcal{G}$ with
radius $\nu>0$ and norm $\left\Vert \cdot\right\Vert $. By (\ref{pf:norm-s})
and the second relation in (\ref{pf:set-s}), 
\[
N_{[]}(\nu,\mathcal{G}_{n,\delta}^{\beta},\left\Vert \cdot\right\Vert _{2,\beta})\leq N_{[]}(\nu,\mathcal{G}_{n,C_{2}\delta}^{1},\left\Vert \cdot\right\Vert _{2})\leq C_{3}\frac{\delta}{\nu^{2d}},
\]
for some positive constant $C_{3}$. Note that the upper bound here
is different from the point identified case. Therefore, for some positive
constant $C_{4}$, it holds 
\begin{equation}
	\varphi_{n}(\delta)=\int_{0}^{\delta}\sqrt{\log N_{[]}(\nu,\mathcal{G}_{n,\delta}^{\beta},\left\Vert \cdot\right\Vert _{2,\beta})}d\nu\leq C_{4}\delta\log\delta^{-1}.\label{pf:phi-s}
\end{equation}

Finally, based on the above entropy condition, we apply the maximal
inequality of Doukhan, Massart and Rio (1995, Theorem 3), i.e., there
exists a positive constant $C_{5}$ depending only on the mixing sequence
$\{\beta_{m}\}$ such that 
\[
P\sup_{g\in\mathcal{G}_{n,\delta}^{\beta}}|\mathbb{G}_{n}g|\leq C_{5}[1+\delta^{-1}q_{G_{n,\delta}}(\min\{1,v_{n}(\delta)\})]\varphi_{n}(\delta),
\]
where $q_{G_{n,\delta}}(v)=\sup_{u\leq v}Q_{G_{n,\delta}}(u)\sqrt{\int_{0}^{u}\beta^{-1}(\tilde{u})d\tilde{u}}$
with the envelope function $G_{n,\delta}$ of $\mathcal{G}_{n,\delta}^{\beta}$
(note: $\mathcal{G}_{n,\delta}^{\beta}$ is a class of bounded functions)
and $v_{n}(\delta)$ is the unique solution of 
\[
\frac{v_{n}(\delta)^{2}}{\int_{0}^{v_{n}(\delta)}\beta^{-1}(\tilde{u})d\tilde{u}}=\frac{\varphi_{n}(\delta)^{2}}{n\delta^{2}}.
\]
Since $\varphi_{n}(\delta)\leq C_{4}\delta\log\delta^{-1}$ from (\ref{pf:phi-s}),
it holds
\[
v_{n}(\delta)\leq C_{5}n^{-1}(\log\delta^{-1})^{2}\leq C_{5}n^{-1}\{\log(nh_{n})^{1/2}\}^{2},
\]
for some positive constant $C_{5}$. Now take some $n_{0}$ such that
$v_{n_{0}}(\delta)\leq1$, and then pick again any $n\geq n_{0}$
and $\delta\in[r_{n}^{-1/2},C^{\prime}]$. We have 
\[
q_{G_{n,\delta}}(\min\{1,v_{n}(\delta)\})\leq C_{6}\sqrt{v_{n}(\delta)}Q_{G_{n,\delta}}(v_{n}(\delta))\leq C_{7}n^{-1/2}\log(nh_{n})^{1/2},
\]
for some positive constants $C_{6}$ and $C_{7}$. Therefore, combining
this with (\ref{pf:phi-s}), the conclusion follows by 
\[
P\sup_{g\in\mathcal{G}_{n,\delta}^{1}}|\mathbb{G}_{n}g|\leq P\sup_{g\in\mathcal{G}_{n,M^{1/2}C_{1}\delta^{1/2}}^{\beta}}|\mathbb{G}_{n}g|\leq C_{8}(\delta\log\delta^{-1})^{1/2},
\]
where the first inequality follows from the first relation in (\ref{pf:set-s}).

\subsection{Proof of Lemma 3}

Pick any $C>0$ and $\varepsilon>0$. Then define $A_{n}=\{\theta\in\Theta\setminus\Theta_{I}:r_{n}^{-1/3}\leq|\theta-\pi_{\theta}|\leq C\}$
and 
\[
R_{n}^{2}=r_{n}^{2/3}\sup_{\theta\in A_{n}}\{|\mathbb{P}_{n}(f_{n,\theta}-f_{n,\pi_{\theta}})-P(f_{n,\theta}-f_{n,\pi_{\theta}})|-\varepsilon|\theta-\pi_{\theta}|^{2}\}.
\]
It is enough to show $R_{n}=O_{p}(1)$. Letting $A_{n,j}=\{\theta\in\Theta:(j-1)r_{n}^{-1/3}\leq|\theta-\pi_{\theta}|<jr_{n}^{-1/3}\}$,
there exists a positive constant $C^{\prime}$ such that 
\begin{eqnarray*}
	&  & P\{R_{n}>m\}\\
	& \le & P\left\{ |\mathbb{P}_{n}(f_{n,\theta}-f_{n,\pi_{\theta}})-P(f_{n,\theta}-f_{n,\pi_{\theta}})|>\varepsilon|\theta-\pi_{\theta}|^{2}+r_{n}^{-2/3}m^{2}\quad\text{for some }\theta\in A_{n}\right\} \\
	& \leq & \sum_{j=1}^{\infty}P\left\{ r_{n}^{2/3}|\mathbb{P}_{n}(f_{n,\theta}-f_{n,\pi_{\theta}})-P(f_{n,\theta}-f_{n,\pi_{\theta}})|>\varepsilon(j-1)^{2}+m^{2}\quad\text{for some }\theta\in A_{n,j}\right\} \\
	& \leq & \sum_{j=1}^{\infty}\frac{C^{\prime}\sqrt{j}}{\varepsilon(j-1)^{2}+m^{2}},
\end{eqnarray*}
for all $m>0$, where the last inequality is due to the Markov inequality
and Lemma MS. Since the above sum is finite for all $m>0$, the conclusion
follows.

\subsection{Proof of Theorem 3}

Pick any $\vartheta\in\hat{\Theta}$. By the definition of $\hat{\Theta}$,
\[
\mathbb{P}_{n}(f_{n,\vartheta}-f_{n,\pi_{\vartheta}})\geq\max_{\theta\in\Theta}\mathbb{P}_{n}f_{n,\theta}-(nh_{n})^{-1/2}\hat{c}-\mathbb{P}_{n}f_{n,\pi_{\vartheta}}\geq-(nh_{n})^{-1/2}\hat{c}.
\]
Now, suppose $H(\vartheta,\Theta_{I})=|\vartheta-\pi_{\vartheta}|>r_{n}^{-1/3}$.
By Lemma 3 and Assumption S (i), 
\begin{align*}
	\mathbb{P}_{n}(f_{n,\vartheta}-f_{n,\pi_{\vartheta}}) & \leq P(f_{n,\vartheta}-f_{n,\pi_{\vartheta}})+\varepsilon|\vartheta-\pi_{\vartheta}|^{2}+r_{n}^{-2/3}R_{n}^{2}\\
	& \leq(-c+\varepsilon)|\vartheta-\pi_{\vartheta}|^{2}+o(|\vartheta-\pi_{\vartheta}|^{2})+O_{p}(r_{n}^{-2/3}),
\end{align*}
for any $\varepsilon>0$. Note that $c$, $\varepsilon$, and $R_{n}$
do not depend on $\vartheta$. By taking $\varepsilon$ small enough,
the convergence rate of $\rho(\hat{\Theta},\Theta_{I})$ is obtained
as 
\[
\rho(\hat{\Theta},\Theta_{I})=\sup_{\vartheta\in\hat{\Theta}}|\vartheta-\pi_{\vartheta}|\leq O_{p}(\hat{c}^{1/2}(nh_{n})^{-1/4}+r_{n}^{-1/3}).
\]
Furthermore, for the maximizer $\hat{\theta}$ of $\mathbb{P}_{n}f_{n,\theta}$,
it holds $\mathbb{P}_{n}(f_{n,\hat{\theta}}-f_{n,\pi_{\hat{\theta}}})\ge0$
and this implies $\hat{\theta}-\pi_{\hat{\theta}}=O_{p}(r_{n}^{-1/3})$.

For the convergence rate of $\rho(\Theta_{I},\hat{\Theta})$, we show
$P\{\Theta_{I}\subset\hat{\Theta}\}\rightarrow1$ for $\hat{c}\to\infty$,
which implies that $\rho(\Theta_{I},\hat{\Theta})$ can converge at
arbitrarily fast rate. To see this, note that 
\begin{eqnarray}
	&  & (nh_{n})^{1/2}\max_{\theta^{\prime}\in\Theta_{I}}|(\max_{\theta\in\Theta}\mathbb{P}_{n}f_{n,\theta}-\mathbb{P}_{n}f_{n,\theta^{\prime}})|\nonumber \\
	& \leq & |\mathbb{G}_{n}(f_{n,\hat{\theta}}-f_{n,\pi_{\hat{\theta}}})|+(nh_{n})^{1/2}|P(f_{n,\hat{\theta}}-f_{n,\pi_{\hat{\theta}}})|+2(nh_{n})^{1/2}|\max_{\theta^{\prime}\in\Theta_{I}}(\mathbb{P}_{n}f_{n,\theta^{\prime}}-Pf_{n,\theta^{\prime}})|\nonumber \\
	& = & 2h_{n}^{1/2}|\max_{\theta^{\prime}\in\Theta_{I}}\mathbb{G}_{n}f_{n,\theta^{\prime}}|+o_{p}(1),\label{eq:step2}
\end{eqnarray}
where the inequality follows from the triangle inequality and the
equality follows from Lemmas MS and 3, Assumption S (i), and the rate
$\hat{\theta}-\pi_{\hat{\theta}}=O_{p}(r_{n}^{-1/3})$ obtained above.
Since $\{h_{n}^{1/2}f_{n,\theta},\theta\in\Theta_{I}\}$ is $P$-Donsker
(Assumption S (i)), it follows $P\{\Theta_{I}\subset\hat{\Theta}\}\rightarrow1$
if $\hat{c}\to\infty$.

\subsection{Proof of Lemma MS'}

To ease notation, let $\nu_{0}=0$. First, we introduce some notation.
Let 
\begin{eqnarray*}
	\mathcal{G}_{n} & = & \{g_{n,s}=f_{n,\theta,\nu}-f_{n,\theta,0}:|\theta-\pi_{\theta}|\leq K_{1},|\nu|\leq a_{n}K_{2},s=(\theta^{\prime},\nu^{\prime})^{\prime}\},\\
	\mathcal{G}_{n,\delta}^{\beta} & = & \{g_{n,s}:\left\Vert g_{n,s}\right\Vert _{2,\beta}<\delta\},\\
	\mathcal{G}_{n,\delta}^{1} & = & \{g_{n,s}:|\theta-\pi_{\theta}|<K_{1}\mbox{ and }|\nu|\leq\delta\},\\
	\mathcal{G}_{n,\delta}^{2} & = & \{g_{n,s}:\left\Vert g_{n,s}\right\Vert _{2}<\delta\}.
\end{eqnarray*}
Since $g_{n,s}$ satisfies eq. (16) of the paper, there exists a positive
constant $C_{1}$ such that $\mathcal{G}_{n,\delta}^{1}\subset\{g_{n,s}\in\mathcal{G}_{n}:\left\Vert g_{n,s}\right\Vert _{2}<C_{1}k_{n}^{1/4}\delta^{1/2}\}$
for all $n$ large enough and all $\delta>0$ small enough. Also,
by the same argument to derive (\ref{pf:g}), there exists a positive
constant $C_{2}$ such that $\left\Vert g_{n,s}\right\Vert _{2,\beta}\leq C_{2}\left\Vert g_{n,s}\right\Vert _{2}$
for all $n$ large enough. The constant $C_{2}$ depends only on the
mixing sequence $\{\beta_{m}\}$. Combining these results, we obtain
\begin{equation}
	\mathcal{G}_{n,\delta}^{1}\subset\mathcal{G}_{n,C_{1}C_{2}k_{n}^{1/4}\delta^{1/2}}^{\beta},\label{pf:g1-1}
\end{equation}
for all $n$ large enough and all $\delta>0$ small enough. On the
other hand, for any $\delta'$ small enough, there exists some $C_{3}$
such that 
\[
\mathcal{G}_{n,\delta'}^{\beta}\subset\mathcal{G}_{n,\delta'}^{2}\subset\mathcal{G}_{n,C_{3}\delta'}^{1},
\]
due to the fact that $\left\Vert \cdot\right\Vert _{2}\leq\left\Vert \cdot\right\Vert _{2,\beta}$
(Doukhan, Massart and Rio, 1995, Lemma 1) and eq. (17) of the paper.
Then the bracketing numbers satisfy 
\[
N_{[]}(\nu,\mathcal{G}_{n,\delta'}^{\beta},\left\Vert \cdot\right\Vert _{2,\beta})\leq N_{[]}(\nu,\mathcal{G}_{n,C_{3}\delta'}^{1},\left\Vert \cdot\right\Vert _{2}).
\]
Furthermore, the bracketing number $N_{[]}(\nu,\mathcal{G}_{n,C_{3}\delta'}^{1},\left\Vert \cdot\right\Vert _{2})$
can be bounded by the covering number of the parameter space, say,
$N_{\Theta}(\nu,[-K_{1},K_{1}]^{d}\times[-C_{3}\delta',C_{3}\delta']^{k_{n}})$
following the argument in Andrews (1993, eq. (4.7)) based on eq. (16)
of the paper. 

Now we set $\delta=a_{n}K_{2}$ so that $\mathcal{G}_{n,\delta}^{1}=\mathcal{G}_{n}$.
Also set $\delta^{\prime}=C_{1}C_{2}k_{n}^{1/4}\delta^{1/2}$ and
compute the covering number $N_{\Theta}(\nu,[-K_{1},K_{1}]^{d}\times[-K_{2}^{\prime}a_{n}^{1/2}k_{n}^{1/4},K_{2}^{\prime}a_{n}^{1/2}k_{n}^{1/4}]^{k_{n}})$,
where $K_{2}^{\prime}=C_{1}C_{2}K_{2}^{1/2}$. By direct calculation,
this covering number is bounded by $(2K_{1})^{d}\left(\frac{\sqrt{d+k_{n}}}{2\nu}\right)^{d+k_{n}}(2K_{2}^{\prime}a_{n}^{1/2}k_{n}^{1/4})^{k_{n}}$.\footnote{The circumradius of the unit $s$-dimensional hypercube is $\sqrt{s}/2$.
	Or $\sqrt{\sum_{i=1}^{s}a_{i}^{2}}/2$ for the hypercube of side lengths
	$(a_{1},\ldots,a_{s})$.} Building on this, we compute the quantities in the maximal inequality
in (\ref{pf:max}). First, 
\begin{eqnarray*}
	\varphi_{n}(\delta^{\prime}) & = & \int_{0}^{\delta^{\prime}}\sqrt{\log N_{[]}(\nu,\mathcal{G}_{n,\delta^{\prime}}^{\beta},\left\Vert \cdot\right\Vert _{2,\beta})}d\nu\\
	& \leq & \int_{0}^{K_{2}^{\prime}a_{n}^{1/2}k_{n}^{1/4}}C_{3}\sqrt{k_{n}\left(\log k_{n}^{3}a_{n}^{2}-\log\nu\right)}d\nu\\
	& \leq & K_{3}a_{n}^{1/2}k_{n}^{3/4}\sqrt{\log k_{n}a_{n}^{-1}},
\end{eqnarray*}
for some $C_{3}$ and $K_{3}$, where the last inequality follows
from the indefinite integral formula $\int\log xdx=\mbox{const.}+x(\log x-1)$.
Second, as in the discussion following (\ref{pf:max}), we have
\[
v_{n}(\delta^{\prime})\leq\varphi_{n}(\delta^{\prime})^{2}/(n\delta^{\prime2})\leq k_{n}\log k_{n}a_{n}^{-1}/n,
\]
which can be made smaller than $1$ for large $n$. Then we obtain
$q_{G_{n,\delta^{\prime}}}(\min\{1,v_{n}(\delta^{\prime})\})\leq C_{4}n^{\kappa}\sqrt{v_{n}(\delta^{\prime})}$
and then $\delta^{\prime-1}q_{G_{n,\delta}}(\min\{1,v_{n}(\delta^{\prime})\})\leq C_{5}$
for $\delta^{\prime}=K_{2}^{\prime}a_{n}^{1/2}k_{n}^{1/4}.$ Putting
these together, we can bound the right hand side of (\ref{pf:max})
by $C_{6}a_{n}^{1/2}k_{n}^{3/4}\sqrt{\log k_{n}a_{n}^{-1}}$ for some
$C_{6}>0$.

\subsection{Proof of Theorem 4}

To ease notation, let $\nu_{0}=0$. From eq. (18) of the paper, we
have 
\begin{equation}
	P(f_{n,\theta,\hat{\nu}}-f_{n,\theta,0})-P(f_{n,\pi_{\theta},\hat{\nu}}-f_{n,\pi_{\theta},0})=o(|\theta-\pi_{\theta}|^{2})+O(|\hat{\nu}|^{2})+O_{p}(r_{n}^{-2/3}),\label{eq:S-1n2}
\end{equation}
for all $\theta$ in a neighborhood of $\Theta_{I}$ and all $\epsilon>0$.
Combining Lemma MS', eq. (18) of the paper, Assumption S (i), and
Lemma 3, 
\begin{eqnarray}
	\mathbb{P}_{n}(f_{n,\theta,\hat{\nu}}-f_{n,\pi_{\theta},\hat{\nu}}) & = & n^{-1/2}\{\mathbb{G}_{n}(f_{n,\theta,\hat{\nu}}-f_{n,\theta,0})-\mathbb{G}_{n}(f_{n,\pi_{\theta},\hat{\nu}}-f_{n,\pi_{\theta},0})+\mathbb{G}_{n}(f_{n,\theta,0}-f_{n,\pi_{\theta},0})\}\nonumber \\
	&  & +P(f_{n,\theta,\hat{\nu}}-f_{n,\theta,0})-P(f_{n,\pi_{\theta},\hat{\nu}}-f_{n,\pi_{\theta},0})+P(f_{n,\theta,0}-f_{n,\pi_{\theta},0})\nonumber \\
	& \leq & O_{p}((nh_{n}a_{n}^{-1})^{-1/2}k_{n}^{3/4}\log^{1/2}n)+\epsilon|\theta-\pi_{\theta}|^{2}+O_{p}(r_{n}^{-2/3})\nonumber \\
	&  & -c|\theta-\pi_{\theta}|^{2}+\epsilon|\theta-\pi_{\theta}|^{2}+O_{p}(|\hat{\nu}|^{2})+O_{p}(r_{n}^{-2/3}),\label{eq:S-1n3}
\end{eqnarray}
for all $\theta$ in a neighborhood of $\Theta_{I}$ and all $\epsilon>0$,
where the inequality follows from eq. (18) of the paper. Here, $\sqrt{\log k_{n}a_{n}^{-1}}$
in Lemma MS' is bounded by $\sqrt{\log n}$ up to a constant. 

Let $\hat{\theta}=\arg\max_{\theta\in\Theta}\mathbb{P}_{n}f_{n,\theta,\hat{v}}$.
If $|\hat{\theta}-\pi_{\hat{\theta}}|>a_{n}+r_{n}^{-1/3}$, then $\mathbb{P}_{n}(f_{n,\theta,\hat{\nu}}-f_{n,\pi_{\theta},\hat{\nu}})\geq0$
and thus by (\ref{eq:S-1n3}), 
\begin{equation}
	|\hat{\theta}-\pi_{\hat{\theta}}|\leq o(a_{n})+O_{p}(r_{n}^{-1/3})+O_{p}((nh_{n}a_{n}^{-1})^{-1/4}k_{n}^{3/8}\log^{1/4}n).\label{eq:rate2}
\end{equation}
Also for any $\theta^{\prime}\in\hat{\Theta}$, if $|\theta^{\prime}-\pi_{\theta^{\prime}}|>a_{n}+r_{n}^{-1/3}$,
it holds 
\[
-(nh_{n})^{-1/2}\hat{c}\leq\max_{\theta\in\Theta}\mathbb{P}_{n}f_{n,\theta,\hat{v}}-\mathbb{P}_{n}f_{n,\pi_{\theta^{\prime}},\hat{v}}-c_{n}^{-1}\hat{c}\leq\mathbb{P}_{n}f_{n,\theta^{\prime},\hat{v}}-\mathbb{P}_{n}f_{n,\pi_{\theta^{\prime}},\hat{v}},
\]
and by (\ref{eq:S-1n3}), 
\[
|\theta^{\prime}-\pi_{\theta^{\prime}}|\leq o(a_{n})+O_{p}(r_{n}^{-1/3})+O_{p}((nh_{n}a_{n}^{-1})^{-1/4}k_{n}^{3/8}\log^{1/4}n)+(nh_{n})^{-1/4}\hat{c}^{1/2}.
\]
It remains to show that $P\{\Theta_{I}\subset\hat{\Theta}\}\rightarrow1$
for $\hat{c}\to\infty$. Proceeding as in (\ref{eq:step2}), we get
\begin{eqnarray*}
	&  & (nh_{n})^{1/2}\max_{\theta^{\prime}\in\Theta_{I}}|(\max_{\theta\in\Theta}\mathbb{P}_{n}f_{n,\theta,\hat{v}}-\mathbb{P}_{n}f_{n,\theta^{\prime},\hat{v}})|\\
	& \leq & |h_{n}^{1/2}\mathbb{G}_{n}(f_{n,\hat{\theta},\hat{v}}-f_{n,\pi_{\hat{\theta}},\hat{v}})|+(nh_{n})^{1/2}|P(f_{n,\hat{\theta},\hat{v}}-f_{n,\pi_{\hat{\theta}},\hat{v}})|\\
	&  & +2(nh_{n})^{1/2}|\max_{\theta^{\prime}\in\Theta_{I}}(\mathbb{P}_{n}f_{n,\theta^{\prime},\hat{v}}-Pf_{n,\theta^{\prime},\hat{v}})|\\
	& = & 2|\max_{\theta^{\prime}\in\Theta_{I}}h_{n}^{1/2}\mathbb{G}_{n}f_{n,\theta^{\prime},\hat{v}}|+o_{p}(1),
\end{eqnarray*}
where the first term after the inequality being $o_{p}(1)$ is due
to Lemmas 3 and MS' and the second term is to (\ref{eq:S-1n2}) and
Assumption S (i) together with the rate for $\hat{\theta}$ in (\ref{eq:rate2}).
Finally, due to Lemma MS' and the class $\{h_{n}^{1/2}f_{n,\theta},\theta\in\Theta_{I}\}$
being a $P$-Donsker, we conclude $\Pr\{\Theta_{I}\subset\hat{\Theta}_{I}\}\rightarrow1$.

\subsection{Proof of Lemma M1}

The proof is similar to that of Lemma M except that for some positive
constant $C^{\prime\prime\prime}$, we have 
\[
\mathcal{G}_{\delta}^{1}\subset\mathcal{G}_{C^{\prime\prime}h_{n}^{-1/2}\delta^{1/2}}^{2}\subset\mathcal{G}_{C^{\prime\prime\prime}h_{n}^{-1/2}\delta^{1/2}}^{\beta},
\]
which reflects the component ``$h_{n}^{2}$'' in Assumption M (iii')
instead of ``$h_{n}$'' in Assumption M (iii). As a consequence of
this change, the upper bound in the maximal inequality becomes $Ch_{n}^{-1/2}\delta^{1/2}$
instead of $C\delta^{1/2}$. All the other parts remain the same.

\subsection{Proof of Theorem 5}

The proof is similar to that of Theorem 1 given Lemma M1.

\newpage{}

\section{Details on examples in Section 3 and 4}

\subsection{Dynamic panel discrete choice}

Consider a dynamic panel data model with a binary dependent variable
\begin{eqnarray*}
	P\{y_{i0}=1|x_{i},\alpha_{i}\} & = & F_{0}(x_{i},\alpha_{i}),\\
	P\{y_{it}=1|x_{i},\alpha_{i},y_{i0},\ldots,y_{it-1}\} & = & F(x_{it}^{\prime}\beta_{0}+\gamma_{0}y_{it-1}+\alpha_{i}),
\end{eqnarray*}
for $i=1,\ldots,n$ and $t=1,2,3$, where $y_{it}$ is binary, $x_{it}$
is a $k$-vector, and both $F_{0}$ and $F$ are unknown functions.
We observe $\{y_{it},x_{it}\}$ but do not observe $\alpha_{i}$.
Honoré and Kyriazidou (2000) proposed the conditional maximum score
estimator for $(\beta_{0},\gamma_{0})$, 
\begin{equation}
	(\hat{\beta},\hat{\gamma})=\arg\max_{(\beta,\gamma)\in\Theta}\sum_{i=1}^{n}K\left(\frac{x_{i2}-x_{i3}}{b_{n}}\right)(y_{i2}-y_{i1})\mathrm{sgn}\{(x_{i2}-x_{i1})^{\prime}\beta+(y_{i3}-y_{i0})\gamma\},\label{eq:HK}
\end{equation}
where $K$ is a kernel function and $b_{n}$ is a bandwidth. In this
case, nonparametric smoothing is introduced to deal with the unknown
link function $F$. Honoré and Kyriazidou (2000) obtained consistency
of this estimator but the convergence rate and limiting distribution
are unknown. Since the criterion function for the estimator varies
with the sample size due to the bandwidth $b_{n}$, the cube root
asymptotic theory of Kim and Pollard (1990) is not applicable. Here
we show that Theorem 1 can be applied to answer these open questions.

Let $z=(z_{1}^{\prime},z_{2},z_{3}^{\prime})^{\prime}$ with $z_{1}=x_{2}-x_{3}$,
$z_{2}=y_{2}-y_{1}$, and $z_{3}=((x_{2}-x_{1})^{\prime},y_{3}-y_{0})$.
Also define $x_{21}=x_{2}-x_{1}$. The criterion function of the estimator
$\hat{\theta}=(\hat{\beta}^{\prime},\hat{\gamma})^{\prime}$ in (\ref{eq:HK})
is written as 
\begin{align}
	f_{n,\theta}(z) & =b_{n}^{-k}K(b_{n}^{-1}z_{1})z_{2}\{\mathrm{sgn}(z_{3}^{\prime}\theta)-\mathrm{sgn}(z_{3}^{\prime}\theta_{0})\}\nonumber \\
	& =e_{n}(z)(\mathbb{I}\{z_{3}^{\prime}\theta\geq0\}-\mathbb{I}\{z_{3}^{\prime}\theta_{0}\geq0\}),\label{eq:HK-f}
\end{align}
and $e_{n}(z)=2b_{n}^{-k}K(b_{n}^{-1}z_{1})z_{2}$. Based on Honoré
and Kyriazidou (2000, Theorem 4), we impose the following assumptions.
\begin{description}
	\item [{(a)}] $\{z_{i}\}_{i=1}^{n}$ is an iid sample. $z_{1}$ has a bounded
	density which is continuously differentiable at zero. The conditional
	density of $z_{1}|z_{2}\neq0,z_{3}$ is positive in a neighborhood
	of zero, and $P\{z_{2}\neq0|z_{3}\}>0$ for almost every $z_{3}$.
	Support of $x_{21}$ conditional on $z_{1}$ in a neighborhood of
	zero is not contained in any proper linear subspace of $\mathbb{R}^{k}$.
	There exists at least one $j\in\{1,\ldots,k\}$ such that $\beta_{0}^{(j)}\neq0$
	and $x_{21}^{(j)}|x_{21}^{j-},z_{1}$, where $x_{21}^{j-}=(x_{21}^{(1)},\ldots,x_{21}^{(j-1)},x_{21}^{(j+1)},\ldots,x_{21}^{(k)})$,
	has everywhere positive conditional density for almost every $x_{21}^{j-}$
	and almost every $z_{1}$ in a neighborhood of zero. $E[z_{2}|z_{3},z_{1}=0]$
	is differentiable in $z_{3}$. $E[z_{2}\mathrm{sgn}((\beta_{0}^{\prime},\gamma_{0})^{\prime}z_{3})|z_{1}]$
	is continuously differentiable at $z_{1}=0$. $F$ is strictly increasing. 
	\item [{(b)}] $K$ is a bounded symmetric density function with bounded
	support and  $\int s_{j}s_{j^{\prime}}K(s)ds<\infty$ for any $j,j^{\prime}\in\{1,\ldots,k\}$.
	As $n\rightarrow\infty$, it holds $nb_{n}^{k}/\ln n\rightarrow\infty$
	and $nb_{n}^{k+3}\rightarrow0$. 
\end{description}
We verify that $\{f_{n,\theta}\}$ satisfies Assumption M with $h_{n}=b_{n}^{k}$.
We first check Assumption M (ii). By the definition of $z_{2}=y_{2}-y_{1}$
(which can take $-1$, $0$, or $1$) and change of variables $a=b_{n}^{-1}z_{1}$,
we obtain 
\[
E[e_{n}(z)^{2}|z_{3}]=4h_{n}^{-1}\int K(a)^{2}p_{1}(b_{n}a|z_{2}\neq0,z_{3})daP\{z_{2}\neq0|z_{3}\},
\]
almost surely for all $n$, where $p_{1}$ is the conditional density
of $z_{1}$ given $z_{2}\neq0$ and $z_{3}$. Thus under (a), $h_{n}E[e_{n}(z)^{2}|z_{3}]>c$
almost surely for some $c>0$. Pick any $\theta_{1}$ and $\theta_{2}$.
Note that 
\begin{eqnarray*}
	h_{n}^{1/2}\left\Vert f_{n,\theta_{1}}-f_{n,\theta_{2}}\right\Vert _{2} & = & \left(P\left\{ h_{n}E[e_{n}(z)^{2}|z_{3}]|\mathbb{I}\{z_{3}^{\prime}\theta_{1}\geq0\}-\mathbb{I}\{z_{3}^{\prime}\theta_{2}\geq0\}|\right\} \right)^{1/2}\\
	& \geq & c^{1/2}P|\mathbb{I}\{z_{3}^{\prime}\theta_{1}\geq0\}-\mathbb{I}\{z_{3}^{\prime}\theta_{2}\geq0\}|\\
	& \geq & c_{1}|\theta_{1}-\theta_{2}|,
\end{eqnarray*}
for some $c_{1}>0$, where the last inequality follows from the same
argument to the maximum score example in Section B.1 of the supplementary
material using (a). Similarly, Assumption M (iii) is verified as 
\[
h_{n}P\sup_{\theta\in\Theta:|\theta-\vartheta|<\varepsilon}|f_{n,\theta}-f_{n,\vartheta}|^{2}\leq C_{1}P\sup_{\theta\in\Theta:|\theta-\vartheta|<\varepsilon}|\mathbb{I}\{z_{3}^{\prime}\theta\geq0\}-\mathbb{I}\{z_{3}^{\prime}\vartheta\geq0\}|\leq C_{2}\varepsilon,
\]
for some positive constants $C_{1}$ and $C_{2}$ and all $\vartheta$
in a neighborhood of $\theta_{0}$ and $n$ large enough. We now verify
Assumption M (i). Since $h_{n}f_{n,\theta}$ is clearly bounded, it
is enough to verify eq. (2) in the paper. A change of variables $a=b_{n}^{-1}z_{1}$
and (b) imply 
\begin{eqnarray*}
	Pf_{n,\theta} & = & \int K(a)E[z_{2}\{\mathrm{sgn}(z_{3}^{\prime}\theta)-\mathrm{sgn}(z_{3}^{\prime}\theta_{0})\}|z_{1}=b_{n}a]p_{1}(b_{n}a)da\\
	& = & p_{1}(0)E[z_{2}\{\mathrm{sgn}(z_{3}^{\prime}\theta)-\mathrm{sgn}(z_{3}^{\prime}\theta_{0})\}|z_{1}=0]\\
	&  & +b_{n}^{2}\int K(a)a^{\prime}\left.\frac{\partial^{2}E[z_{2}\{\mathrm{sgn}(z_{3}^{\prime}\theta)-\mathrm{sgn}(z_{3}^{\prime}\theta_{0})\}|z_{1}=t]p_{1}(t)}{\partial t\partial t^{\prime}}\right|_{t=t_{a}}ada,
\end{eqnarray*}
where $t_{a}$ is a point on the line joining $a$ and $0$, and the
second equality follows from the dominated convergence and mean value
theorems. Since $b_{n}^{2}=o((nb_{n}^{k})^{-2/3})$ by (b), the second
term is negligible. Thus, for the condition in eq. (2), it is enough
to derive a second order expansion of $E[z_{2}\{\mathrm{sgn}(z_{3}^{\prime}\theta)-\mathrm{sgn}(z_{3}^{\prime}\theta_{0})\}|z_{1}=0]$.
Let $\mathcal{Z}_{\theta}=\{z_{3}:\mathbb{I}\{z_{3}^{\prime}\theta\geq0\}\neq\mathbb{I}\{z_{3}^{\prime}\theta_{0}\geq0\}\}$.
Honoré and Kyriazidou (2000, p. 872) showed that 
\[
-E[z_{2}\{\mathrm{sgn}(z_{3}^{\prime}\theta)-\mathrm{sgn}(z_{3}^{\prime}\theta_{0})\}|z_{1}=0]=2\int_{\mathcal{Z}_{\theta}}|E[z_{2}|z_{1}=0,z_{3}]|dF_{z_{3}|z_{1}=0}>0,
\]
for all $\theta\neq\theta_{0}$ on the unit sphere and that $\mathrm{sgn}(E[z_{2}|z_{3},z_{1}=0])=\mathrm{sgn}(z_{3}^{\prime}\theta_{0})$.
Therefore, by applying the same argument as Kim and Pollard (1990,
pp. 214-215), we obtain
\[
\left.\frac{\partial}{\partial\theta}E[z_{2}\mathrm{sgn}(z_{3}^{\prime}\theta)|z_{1}=0]\right\vert _{\theta=\theta_{0}}=0,
\]
and
\[
-\frac{\partial^{2}E[z_{2}\{\mathrm{sgn}(z_{3}^{\prime}\theta)-\mathrm{sgn}(z_{3}^{\prime}\theta_{0})\}|z_{1}=0]}{\partial\theta\partial\theta^{\prime}}=\int\mathbb{I}\{z_{3}^{\prime}\theta_{0}=0\}\dot{\kappa}(z_{3})^{\prime}\theta_{0}z_{3}z_{3}^{\prime}p_{3}(z_{3}|z_{1}=0)d\mu_{\theta_{0}},
\]
where $\dot{\kappa}(z_{3})=\frac{\partial}{\partial z_{3}}E[z_{2}|z_{3},z_{1}=0]$,
$p_{3}$ is the conditional density of $z_{3}$ given $z_{1}=0$,
and $\mu_{\theta_{0}}$ is the surface measure on the boundary of
$\{z_{3}:z_{3}^{\prime}\theta_{0}\geq0\}$. Combining these results,
the condition in eq. (2) is satisfied with the negative definite matrix
\begin{equation}
	V=-2p_{1}(0)\int\mathbb{I}\{z_{3}^{\prime}\theta_{0}=0\}\dot{\kappa}(z_{3})^{\prime}\theta_{0}z_{3}z_{3}^{\prime}p_{3}(z_{3}|z_{1}=0)d\mu_{\theta_{0}}.\label{eq:HK-V}
\end{equation}

We now verify eq. (7) of the paper to apply the central limit theorem
in Lemma C. In this example, the normalized criterion function is
written as
\[
g_{n,s}(z)=n^{1/6}b_{n}^{2k/3}e_{n}(z)A_{n,s}(z_{3}),
\]
where $e_{n}(z)=2b_{n}^{-k}K(b_{n}^{-1}z_{1})z_{2}$ and $A_{n,s}(z_{3})=\mathbb{I}\{z_{3}^{\prime}(\theta_{0}+sn^{-1/3}b_{n}^{-k/3})\geq0\}-\mathbb{I}\{z_{3}^{\prime}\theta_{0}\geq0\}$.
Since $|z_{2}|\le1$ and $|A_{n,s}(z_{3})|$ takes only $0$ or $1$,
it holds
\begin{eqnarray*}
	P\{|g_{n,s}|\geq c\} & \le & P\left\{ \left.|K(b_{n}^{-1}z_{1})|\geq2^{-1}cn^{-1/6}b_{n}^{k/3}\right||A_{n,s}(z_{3})|=1\right\} P\{|A_{n,s}(z_{3})|=1\}\\
	& \leq & P\left\{ \left.|b_{n}^{-1}z_{1}|\leq C\right||A_{n,s}(z_{3})|=1\right\} (nb_{n}^{k})^{-1/3}\\
	& \le & C^{\prime}(nb_{n}^{-2k})^{-1/3},
\end{eqnarray*}
for some $C,C^{\prime}>0$, where the second inequality follows from
the bonded support of $K$, boundedness of the conditional density
of $z_{1}$ given $|A_{n,s}(z_{3})|=1$ (by (a)), the fact that $P\{|A_{n,s}(z_{3})|=1\}$
is proportional to $(nb_{n}^{k})^{-1/3}$, and the last inequality
follows from boundedness of the conditional density of $z_{1}$ given
$|A_{n,s}(z_{3})|=1$ (by (a)). Since $h_{n}=b_{n}^{k}$ in this example,
we can apply Lemma 2 to conclude that the condition in eq. (5) holds
true.

Since the criterion function (\ref{eq:HK-f}) satisfies Assumption
M and the Lindeberg-type condition in eq. (5), Theorem 1 implies the
limiting distribution of Honoré and Kyriazidou's (2000) estimator
as in eq. (9). The matrix $V$ is given in (\ref{eq:HK-V}). The covariance
kernel $H$ is obtained in the same manner as Kim and Pollard (1990).
That is, decompose $z_{3}$ into $r^{\prime}\theta_{0}+\bar{z}_{3}$
with $\bar{z}_{3}$ orthogonal to $\theta_{0}.$ Then it holds $H(s_{1},s_{2})=L(s_{1})+L(s_{2})-L(s_{1}-s_{2})$,
where 
\[
L(s)=4p_{1}(0)\int|\bar{z}_{3}^{\prime}s|p_{3}(0,\bar{z}_{3}|z_{1}=0)d\bar{z}_{3}.
\]

\subsection{Random coefficient binary choice}

As a new statistical model which can be covered by our asymptotic
theory, let us consider the regression model with a random coefficient
$y_{t}=x_{t}^{\prime}\theta(w_{t})+u_{t}$. We observe $x_{t}\in\mathbb{R}^{d}$,
$w_{t}\in\mathbb{R}^{k}$, and the sign of $y_{t}$. We wish to estimate
$\theta_{0}=\theta(c)$ at some given $c\in\mathbb{R}^{k}$. In this
setup, we can consider a localized version of the maximum score estimator
\[
\hat{\theta}=\arg\max_{\theta\in S}\sum_{t=1}^{n}K\left(\frac{w_{t}-c}{b_{n}}\right)[\mathbb{I}\{y_{t}\ge0,x_{t}^{\prime}\theta\ge0\}+\mathbb{I}\{y_{t}<0,x_{t}^{\prime}\theta<0\}],
\]
where $S$ is the surface of the unit sphere in $\mathbb{R}^{d}$.
Let $h(x,u)=\mathbb{I}\{x^{\prime}\theta_{0}+u\ge0\}-\mathbb{I}\{x^{\prime}\theta_{0}+u<0\}$.
We impose the following assumptions.
\begin{description}
	\item [{(a)}] $\{x_{t},w_{t},u_{t}\}$ satisfies Assumption D. The density
	$p(x,w)$ of $(x_{t},w_{t})$ is continuous at all $x$ and $w=c$.
	The conditional distribution $x|w=c$ has compact support and continuously
	differentiable conditional density. The angular component of $x|w=c$,
	considered as a random variable on $S$, has a bounded and continuous
	density, and the density for the orthogonal angle to $\theta_{0}$
	is bounded away from zero. 
	\item [{(b)}] Assume that $|\theta_{0}|=1$, $\mathrm{median}(u|x,w=c)=0$,
	the function $\kappa(x,w)=E[h(x_{t},u_{t})|x_{t}=x,w_{t}=w]$ is continuous
	at all $x$ and $w=c$, $\kappa(x,c)$ is non-negative for $x^{\prime}\theta_{0}\ge0$
	and non-positive for $x^{\prime}\theta_{0}<0$ and is continuously
	differentiable in $x$, and
	\[
	P\left\{ x^{\prime}\theta_{0}=0,\left(\frac{\partial\kappa(x,w)}{\partial x}\right)^{\prime}\theta_{0}p(x,w)>0|w=c\right\} >0.
	\]
	\item [{(c)}] $K$ is a bounded symmetric density function with $\int s^{2}K(s)ds<\infty$.
	As $n\rightarrow\infty$, it holds $nb_{n}^{k^{\prime}}\to\infty$
	for some $k^{\prime}>k$. 
\end{description}
Note that the criterion function is written as 
\[
f_{n,\theta}(x,w,u)=\frac{1}{h_{n}}K\left(\frac{w-c}{h_{n}^{1/k}}\right)h(x,u)[\mathbb{I}\{x^{\prime}\theta\ge0\}-\mathbb{I}\{x^{\prime}\theta_{0}\ge0\}],
\]
where $h_{n}=b_{n}^{k}$. We can see that $\hat{\theta}=\arg\max_{\theta\in S}\mathbb{P}_{n}f_{n,\theta}$
and $\theta_{0}=\arg\max_{\theta\in S}\lim_{n\to\infty}Pf_{n,\theta}$.
Existence and uniqueness of $\theta_{0}$ are guaranteed by the change
of variables and (b) (see, Manski, 1985). Also the uniform law of
large numbers for an absolutely regular process by Nobel and Dembo
(1993, Theorem 1) implies $\sup_{\theta\in S}|\mathbb{P}_{n}f_{n,\theta}-Pf_{n,\theta}|\overset{p}{\to}0$.
Therefore, $\hat{\theta}$ is consistent for $\theta_{0}$.

We next compute the expected value and covariance kernel of the limit
process (i.e., $V$ and $H$ in Theorem 1). Due to strict stationarity
(in Assumption D), we can apply the same argument to Kim and Pollard
(1990, pp. 214-215) to obtain the second derivative 
\[
V=\lim_{n\to\infty}\left.\frac{\partial^{2}Pf_{n,\theta}}{\partial\theta\partial\theta^{\prime}}\right|_{\theta=\theta_{0}}=-\int\mathbb{I}\{x^{\prime}\theta_{0}=0\}\left(\frac{\partial\kappa(x,c)}{\partial x}\right)^{\prime}\theta_{0}p(x,c)xx^{\prime}d\sigma(x),
\]
where $\sigma$ is the surface measure on the boundary of the set
$\{x:x^{\prime}\theta_{0}\ge0\}$. The matrix $V$ is negative definite
under the last condition of (b). Now pick any $s_{1}$ and $s_{2}$,
and define $q_{n,t}=f_{n,\theta_{0}+(nh_{n})^{-1/3}s_{1}}(x_{t},w_{t},u_{t})-f_{n,\theta_{0}+(nh_{n})^{-1/3}s_{2}}(x_{t},w_{t},u_{t})$.
The covariance kernel is written as\\
$H(s_{1},s_{2})=\frac{1}{2}\{L(s_{1},0)+L(0,s_{2})-L(s_{1},s_{2})\}$,
where 
\[
L(s_{1},s_{2})=\lim_{n\to\infty}(nh_{n})^{4/3}\mathrm{Var}(\mathbb{P}_{n}q_{n,t})=\lim_{n\to\infty}(nh_{n})^{1/3}\{\mathrm{Var}(q_{n,t})+\sum_{m=1}^{\infty}\mathrm{Cov}(q_{n,t},q_{n,t+m})\}.
\]
The limit of $(nh_{n})^{1/3}\mathrm{Var}(q_{n,t})$ is obtained in
the same manner as Kim and Pollard (1990, p. 215). For the covariance,
the $\alpha$-mixing inequality implies 
\[
|\mathrm{Cov}(q_{n,t},q_{n,t+m})|\leq C\beta_{m}\left\Vert q_{n,t}\right\Vert _{p}^{2}=O(\rho^{m})O((nh_{n})^{-\frac{2}{3p}}h_{n}^{\frac{2(1-p)}{p}}),
\]
for some $C>0$ and $p>2$, where the equality follows from the change
of variables and Assumption D. Also, by the change of variables $|\mathrm{Cov}(q_{n,t},q_{n,t+m})|=|Pq_{n,t}q_{n,t+m}-(Pq_{n,t})^{2}|=O((nh_{n})^{-2/3})$.
By using these bounds (note: if $0<A\le\min\{B_{1},B_{2}\}$, then
$A\le B_{1}^{\ell}B_{2}^{1-\ell}$ for any $\ell\in[0,1]$), there
exists a positive constant $C^{\prime}$ such that 
\[
(nh_{n})^{1/3}\sum_{m=1}^{\infty}|\mathrm{Cov}(q_{n,t},q_{n,t+m})|\leq C^{\prime}(nh_{n})^{-\frac{1}{3}+\frac{2(p-1)\ell}{3}}h_{n}^{-\frac{2(p-1)\ell}{p}}\sum_{m=1}^{\infty}\rho^{\ell m},
\]
for any $\ell\in[0,1]$. Thus, by taking $\ell$ sufficiently small,
we obtain
\[
\lim_{n\rightarrow\infty}(nh_{n})^{1/3}\sum_{m=1}^{\infty}\mathrm{Cov}(q_{n,t},q_{n,t+m})=0,
\]
due to $nb_{n}^{k^{\prime}}\to\infty$.

We now verify that $\{f_{n,\theta}:\theta\in S\}$ satisfies Assumption
M with $h_{n}=b_{n}^{k}$. Assumption M (i) is already verified. By
the change of variables and Jensen's inequality (also note that $h(x,u)^{2}=1$
everywhere), there exists a positive constant $C$ such that 
\begin{eqnarray*}
	h_{n}^{1/2}\left\Vert f_{n,\theta_{1}}-f_{n,\theta_{2}}\right\Vert _{2} & = & \sqrt{\int\int K(s)^{2}|\mathbb{I}\{x^{\prime}\theta_{1}\ge0\}-\mathbb{I}\{x^{\prime}\theta_{2}\ge0\}|p(x,c+sb_{n})dxds}\\
	& \ge & CE\left[\left.|\mathbb{I}\{x^{\prime}\theta_{1}\ge0\}-\mathbb{I}\{x^{\prime}\theta_{2}\ge0\}|\right|w=c\right]\\
	& = & CP\{x^{\prime}\theta_{1}\ge0>x^{\prime}\theta_{2}\mbox{ or }x^{\prime}\theta_{2}\ge0>x^{\prime}\theta_{1}|w=c\},
\end{eqnarray*}
for all $\theta_{1},\theta_{2}\in S$ and all $n$ large enough. Since
the right hand side is the conditional probability for a pair of wedge
shaped regions with an angle of order $|\theta_{1}-\theta_{2}|$,
the last condition in (a) implies Assumption M (ii). For Assumption
M (iii), there exists a positive constant $C^{\prime}$ such that
for any $\varepsilon>0$ 
\begin{eqnarray*}
	&  & P\sup_{\theta\in\Theta:|\theta-\vartheta|<\varepsilon}h_{n}|f_{n,\theta}-f_{n,\vartheta}|^{2}\\
	& = & \int\int K(s)^{2}\sup_{\theta\in\Theta:|\theta-\vartheta|<\varepsilon}|[\mathbb{I}\{x^{\prime}\theta\ge0\}-\mathbb{I}\{x^{\prime}\vartheta\ge0\}]|^{2}p(x,c+sb_{n})dxds\\
	& \le & C^{\prime}E\left[\left.\sup_{\theta\in\Theta:|\theta-\vartheta|<\varepsilon}|[\mathbb{I}\{x^{\prime}\theta\ge0\}-\mathbb{I}\{x^{\prime}\vartheta\ge0\}]|^{2}\right|w=c\right],
\end{eqnarray*}
for all $\vartheta$ in a neighborhood of $\theta_{0}$ and $n$ large
enough. Again, the right hand side is the conditional probability
for a pair of wedge shaped regions with an angle of order $\varepsilon$.
Thus the last condition in (a) also guarantees Assumption M (iii).
Since $\{f_{n,\theta}:\theta\in S\}$ satisfies Assumption M, Theorem
1 implies the limiting distribution of $(nh_{n})^{1/3}(\hat{\theta}-\theta_{0})$
for the random coefficient model.

\subsection{Minimum volume predictive region \label{subsec:mv}}

As an illustration of Theorem 2, we now consider the example in eq.
(1), the minimum volume predictor for a strictly stationary process
proposed by Polonik and Yao (2000). Suppose we are interested in predicting
$y\in\mathbb{R}$ from $x\in\mathbb{R}$ based on the observations
$\{y_{t},x_{t}\}$. The minimum volume predictor of $y$ at $x=c$
in the class $\mathcal{I}$ of intervals of $\mathbb{R}$ at level
$\alpha\in[0,1]$ is defined as 
\[
\hat{I}=\arg\min_{S\in\mathcal{I}}\mu(S)\quad\mbox{s.t. }\hat{P}(S)\ge\alpha,
\]
where $\mu$ is the Lebesgue measure and $\hat{P}(S)=\sum_{t=1}^{n}\mathbb{I}\{y_{t}\in S\}K\left(\frac{x_{t}-c}{h_{n}}\right)/\sum_{t=1}^{n}K\left(\frac{x_{t}-c}{h_{n}}\right)$
is the kernel estimator of the conditional probability $P\{y_{t}\in S|x_{t}=c\}$.
Since $\hat{I}$ is an interval, it can be written as $\hat{I}=[\hat{\theta}-\hat{\nu},\hat{\theta}+\hat{\nu}]$,
where 
\[
\hat{\theta}=\arg\max_{\theta\in\mathbb{R}}\hat{P}([\theta-\hat{\nu},\theta+\hat{\nu}]),\qquad\hat{\nu}=\inf\{\nu\in\mathbb{R}:\sup_{\theta\in\mathbb{R}}\hat{P}([\theta-\nu,\theta+\nu])\ge\alpha\}.
\]
To study the asymptotic property of $\hat{I}$, we impose the following
assumptions. 
\begin{description}
	\item [{(a)}] $\{y_{t},x_{t}\}$ satisfies Assumption D. $I_{0}=[\theta_{0}-\nu_{0},\theta_{0}+\nu_{0}]$
	is the unique shortest interval such that $P\{y_{t}\in I_{0}|x_{t}=c\}\ge\alpha$.
	The conditional density $\gamma_{y|x=c}$ of $y_{t}$ given $x_{t}=c$
	is bounded and strictly positive at $\theta_{0}\pm\nu_{0}$, and its
	derivative satisfies $\dot{\gamma}_{y|x=c}(\theta_{0}-\nu_{0})-\dot{\gamma}_{y|x=c}(\theta_{0}+\nu_{0})>0$. 
	\item [{(b)}] $K$ is bounded and symmetric, and satisfies $\lim_{a\to\infty}|a|K(a)=0$.
	As $n\to\infty$, $nh_{n}\to\infty$ and $nh_{n}^{4}\to0$. 
\end{description}
For notational convenience, assume $\theta_{0}=0$ and $\nu_{0}=1$.
We first derive the convergence rate for $\hat{\nu}$. Note that $\hat{\nu}=\inf\{\nu\in\mathbb{R}:\sup_{\theta\in\mathbb{R}}\hat{g}([\theta-\nu,\theta+\nu])\ge\alpha\hat{\gamma}(c)\}$,
where $\hat{g}(S)=\frac{1}{nh_{n}}\sum_{t=1}^{n}\mathbb{I}\{y_{t}\in S\}K\left(\frac{x_{t}-c}{h_{n}}\right)$
and $\hat{\gamma}(c)=\frac{1}{nh_{n}}\sum_{t=1}^{n}K\left(\frac{x_{t}-c}{h_{n}}\right)$.
By applying Lemma M' and a central limit theorem, we can obtain uniform
convergence rate
\begin{eqnarray*}
	&  & \max\left\{ |\hat{\gamma}(c)-\gamma(c)|,\sup_{(\theta,\nu)\in\mathbb{R}^{2}}|\hat{g}([\theta-\nu,\theta+\nu])-P\{y_{t}\in[\theta-\nu,\theta+\nu]|x_{t}=c\}\gamma(c)|\right\} \\
	& = & O_{p}((nh_{n})^{-1/2}+h_{n}^{2}).
\end{eqnarray*}
Thus the same argument to Kim and Pollard (1990, pp. 207-208) yields
$\hat{\nu}-1=O_{p}((nh_{n})^{-1/2}+h_{n}^{2})$. Let $\hat{\theta}=\arg\min_{\theta\in\mathbb{R}}\hat{g}([\theta-\hat{\nu},\theta+\hat{\nu}])$.
Consistency follows from uniqueness of $(\theta_{0},\nu_{0})$ in
(a) and the uniform convergence 
\[
\sup_{\theta\in\mathbb{R}}|\hat{g}([\theta-\hat{\nu},\theta+\hat{\nu}])-P\{y_{t}\in[\theta-1,\theta+1]|x_{t}=c\}\gamma(c)|\overset{p}{\to}0,
\]
which is obtained by applying Nobel and Dembo (1993, Theorem 1).

Now let $z=(y,x)^{\prime}$ and 
\[
f_{n,\theta,\nu}(z)=\frac{1}{h_{n}}K\left(\frac{x-c}{h_{n}}\right)[\mathbb{I}\{y\in[\theta-\nu,\theta+\nu]\}-\mathbb{I}\{y\in[-\nu,\nu]\}].
\]
Note that $\hat{\theta}=\arg\max_{\theta\in\mathbb{R}}\mathbb{P}_{n}f_{n,\theta,\hat{\nu}}$.
We apply Theorem 2 to obtain the convergence rate of $\hat{\theta}$.
For the condition in eq. (10) of the paper, observe that 
\begin{eqnarray*}
	&  & P(f_{n,\theta,\nu}-f_{n,0,1})=P(f_{n,\theta,\nu}-f_{n,0,\nu})+P(f_{n,0,\nu}-f_{n,0,1})\\
	& = & -\frac{1}{2}\{-\dot{\gamma}_{y|x}(1|c)+\dot{\gamma}_{y|x}(-1|c)\}\gamma_{x}(c)\theta^{2}+\{\dot{\gamma}_{y|x}(1|c)+\dot{\gamma}_{y|x}(-1|c)\}\gamma_{x}(c)\theta\nu\\
	&  & +o(\theta^{2}+|\nu-1|^{2})+O(h_{n}^{2}).
\end{eqnarray*}
The condition in eq. (10) holds with $V_{1}=\{\dot{\gamma}_{y|x}(1|c)-\dot{\gamma}_{y|x}(-1|c)\}\gamma_{x}(c)$.
Assumption M (iii) for $\{f_{n,\theta,\nu}:\theta\in\mathbb{R},\nu\in\mathbb{R}\}$
is verified in the same manner as in Section B.2 of the supplementary
material. It remains to verify Assumption M (ii) for the class $\{f_{n,\theta,1}:\theta\in\mathbb{R}\}$.
Pick any $\theta_{1}$ and $\theta_{2}$. Some expansions yield 
\begin{eqnarray*}
	&  & h_{n}\left\Vert f_{n,\theta_{1},1}-f_{n,\theta_{2},1}\right\Vert _{2}^{2}\\
	& = & \int K(a)^{2}\left|\begin{array}{c}
		\Gamma_{y|x}(\theta_{2}+1|x=c+ah_{n})-\Gamma_{y|x}(\theta_{1}+1|x=c+ah_{n})\\
		+\Gamma_{y|x}(\theta_{2}-1|x=c+ah_{n})-\Gamma_{y|x}(\theta_{1}-1|x=c+ah_{n})
	\end{array}\right|\gamma_{x}(c+ah_{n})da\\
	& \ge & \int K(a)^{2}\{\gamma_{y|x}(\dot{\theta}+1|x=c+ah_{n})+\gamma_{y|x}(\ddot{\theta}-1|x=c+ah_{n})\}\gamma_{x}(c+ah_{n})da|\theta_{1}-\theta_{2}|,
\end{eqnarray*}
where $\Gamma_{y|x}$ is the conditional distribution function of
$y$ given $x$, and $\dot{\theta}$ and $\ddot{\theta}$ are points
between $\theta_{1}$ and $\theta_{2}$. By (a), Assumption M (ii)
is satisfied. Therefore, we can conclude that $\hat{\nu}-\nu_{0}=O_{p}((nh_{n})^{-1/2}+h_{n}^{2})$
and $\hat{\theta}-\theta_{0}=O_{p}((nh_{n})^{-1/3}+h_{n})$. This
result confirms positively the conjecture of Polonik and Yao (2000,
Remark 3b) on the exact convergence rate of $\hat{I}$.

\subsection{Dynamic maximum score}

As a further application of Theorem 1, consider the maximum score
estimator (Manski, 1975) for the regression model $y_{t}=x_{t}^{\prime}\theta_{0}+u_{t}$,
that is 
\[
\hat{\theta}=\arg\max_{\theta\in S}\sum_{t=1}^{n}[\mathbb{I}\{y_{t}\ge0,x_{t}^{\prime}\theta\ge0\}+\mathbb{I}\{y_{t}<0,x_{t}^{\prime}\theta<0\}],
\]
where $S$ is the surface of the unit sphere in $\mathbb{R}^{d}$.
Since $\hat{\theta}$ is determined only up to scalar multiples, we
standardize it to be unit length. A key insight of this estimator
is to explore a median or quantile restriction in disturbances of
latent variable models to construct a population criterion that identifies
structural parameters of interest. 

We impose the following assumptions. Let $h(x,u)=\mathbb{I}\{x^{\prime}\theta_{0}+u\ge0\}-\mathbb{I}\{x^{\prime}\theta_{0}+u<0\}$. 
\begin{description}
	\item [{(a)}] $\{x_{t},u_{t}\}$ satisfies Assumption D. $x_{t}$ has compact
	support and a continuously differentiable density $p$. The angular
	component of $x_{t}$, considered as a random variable on $S$, has
	a bounded and continuous density, and the density for the orthogonal
	angle to $\theta_{0}$ is bounded away from zero. 
	\item [{(b)}] Assume that $|\theta_{0}|=1$, $\mathrm{median}(u|x)=0$,
	the function $\kappa(x)=E[h(x_{t},u_{t})|x_{t}=x]$ is non-negative
	for $x^{\prime}\theta_{0}\ge0$ and non-positive for $x^{\prime}\theta_{0}<0$
	and is continuously differentiable, and $P\{x^{\prime}\theta_{0}=0,\dot{\kappa}(x)^{\prime}\theta_{0}p(x)>0\}>0$. 
\end{description}
Except for Assumption D, which allows dependent observations, all
assumptions are similar to the ones in Kim and Pollard (1990, Section
6.4). First, note that the criterion function is written as 
\[
f_{\theta}(x,u)=h(x,u)[\mathbb{I}\{x^{\prime}\theta\ge0\}-\mathbb{I}\{x^{\prime}\theta_{0}\ge0\}].
\]
We can see that $\hat{\theta}=\arg\max_{\theta\in S}\mathbb{P}_{n}f_{\theta}$
and $\theta_{0}=\arg\max_{\theta\in S}Pf_{\theta}$. Existence and
uniqueness of $\theta_{0}$ are guaranteed by (b) (see, Manski, 1985).
Also the uniform law of large numbers for an absolutely regular process
by Nobel and Dembo (1993, Theorem 1) implies $\sup_{\theta\in S}|\mathbb{P}_{n}f_{\theta}-Pf_{\theta}|\overset{p}{\to}0$.
Therefore, $\hat{\theta}$ is consistent for $\theta_{0}$.

We now verify that $\{f_{\theta}:\theta\in S\}$ satisfy Assumption
M with $h_{n}=1$. Assumption M (i) is already verified. By Jensen's
inequality, 
\[
\left\Vert f_{\theta_{1}}-f_{\theta_{2}}\right\Vert _{2}=\sqrt{P|\mathbb{I}\{x^{\prime}\theta_{1}\ge0\}-\mathbb{I}\{x^{\prime}\theta_{2}\ge0\}|}\ge P\{x^{\prime}\theta_{1}\ge0>x^{\prime}\theta_{2}\mbox{ or }x^{\prime}\theta_{2}\ge0>x^{\prime}\theta_{1}\},
\]
for any $\theta_{1},\theta_{2}\in S$. Since the right hand side is
the probability for a pair of wedge shaped regions with an angle of
order $|\theta_{1}-\theta_{2}|$, the last condition in (a) implies
Assumption M (ii). For Assumption M (iii), pick any $\varepsilon>0$
and observe that 
\[
P\sup_{\theta\in\Theta:|\theta-\vartheta|<\varepsilon}|f_{\theta}-f_{\vartheta}|^{2}=P\sup_{\theta\in\Theta:|\theta-\vartheta|<\varepsilon}\mathbb{I}\{x^{\prime}\theta\ge0>x^{\prime}\vartheta\mbox{ or }x^{\prime}\vartheta\ge0>x^{\prime}\theta\},
\]
for all $\vartheta$ in a neighborhood of $\theta_{0}$. Again, the
right hand side is the probability for a pair of wedge shaped regions
with an angle of order $\varepsilon$. Thus the last condition in
(a) also guarantees Assumption M (iii). This yields the convergence
rate of $n^{-1/3}$ for $\hat{\theta}$ and the stochastic equicontinuity
of the empirical process of the rescaled and centered functions $g_{n,s}=n^{1/6}(f_{\theta_{0}+sn^{-1/3}}-f_{\theta_{0}})$.
For its finite dimensional convergence, we can check the Lindeberg
condition in Lemma C (i.e., eq. (5) in the paper) by Lemma 2 as in
Section B.2, see (\ref{eq:gns}) below.

We next compute the expected value and covariance kernel of the limit
process (i.e., $V$ and $H$ in Theorem 1). Due to strict stationarity
(in Assumption D), we can apply the same argument to Kim and Pollard
(1990, pp. 214-215) to derive the second derivative 
\[
V=\left.\frac{\partial^{2}Pf_{\theta}}{\partial\theta\partial\theta^{\prime}}\right|_{\theta=\theta_{0}}=-\int\mathbb{I}\{x^{\prime}\theta_{0}=0\}\dot{\kappa}(x)^{\prime}\theta_{0}p(x)xx^{\prime}d\sigma,
\]
where $\sigma$ is the surface measure on the boundary of the set
$\{x:x^{\prime}\theta_{0}\ge0\}$. The matrix $V$ is negative definite
under the last condition of (b). Now pick any $s_{1}$ and $s_{2}$,
and define $q_{n,t}=f_{\theta_{0}+n^{-1/3}s_{1}}(x_{t},u_{t})-f_{\theta_{0}+n^{-1/3}s_{2}}(x_{t},u_{t})$.
The covariance kernel is written as $H(s_{1},s_{2})=\frac{1}{2}\{L(s_{1},0)+L(0,s_{2})-L(s_{1},s_{2})\}$,
where 
\[
L(s_{1},s_{2})=\lim_{n\to\infty}n^{4/3}\mathrm{Var}(\mathbb{P}_{n}q_{n,t})=\lim_{n\to\infty}n^{1/3}\{\mathrm{Var}(q_{n,t})+\sum_{m=1}^{\infty}\mathrm{Cov}(q_{n,t},q_{n,t+m})\}.
\]
The limit of $n^{1/3}\mathrm{Var}(q_{n,t})$ is given in Kim and Pollard
(1990, p. 215). For $\mathrm{Cov}(q_{n,t},q_{n,t+m})$, we note that
$q_{n,t}$ takes only three values, $-1$, $0$, or $1$. The definition
of $\beta_{m}$ and Assumption D imply 
\[
|P\{q_{n,t}=j,q_{n,t+m}=k\}-P\{q_{n,t}=j\}P\{q_{n,t+m}=k\}|\leq n^{-2/3}\beta_{m},
\]
for all $n,m\ge1$ and $j,k=-1,0,1$, i.e., $\{q_{n,t}\}$ is a $\beta$-mixing
array and its mixing coefficients are bounded by $n^{-2/3}\beta_{m}$.
Then, $\{q_{n,t}\}$ is an $\alpha$-mixing array whose mixing coefficients
are bounded by $2n^{-2/3}\beta_{m}$ as well. By applying the $\alpha$-mixing
inequality, the covariance is bounded as 
\[
\mathrm{Cov}(q_{n,t},q_{n,t+m})\le Cn^{-2/3}\beta_{m}\left\Vert q_{n,t}\right\Vert _{p}^{2},
\]
for some $C>0$ and $p>2$. Note that 
\[
\left\Vert q_{n,t}\right\Vert _{p}^{2}\leq[P|\mathbb{I}\{x^{\prime}(\theta_{0}+s_{1}n^{-1/3})>0\}-\mathbb{I}\{x^{\prime}(\theta_{0}+s_{2}n^{-1/3})>0\}|]^{2/p}=O(n^{-2/(3p)}).
\]
Combining these results, we get $n^{1/3}\sum_{m=1}^{\infty}\mathrm{Cov}(q_{n,t},q_{n,t+m})\to0$
as $n\to\infty$. Therefore, the covariance kernel $H$ is same as
the independent case in Kim and Pollard (1990, p. 215).

Since $\{f_{\theta}:\theta\in S\}$ satisfies Assumption M, Theorem
1 implies that even if the data obey a dependence process specified
in Assumption D, the maximum score estimator possesses the same limiting
distribution as the independent sampling case.

\subsection{Dynamic least median of squares}

As another application of Theorem 2, consider the least median of
squares estimator for the regression model $y_{t}=x_{t}^{\prime}\beta_{0}+u_{t}$,
that is 
\[
\hat{\beta}=\arg\min_{\beta}\mathrm{median}\{(y_{1}-x_{1}^{\prime}\beta)^{2},\ldots,(y_{n}-x_{n}^{\prime}\beta)^{2}\}.
\]
We impose the following assumptions. 
\begin{description}
	\item [{(a)}] $\{x_{t},u_{t}\}$ satisfies Assumption D. $\{x_{t}\}$ and
	$\{u_{t}\}$ are independent. $P|x_{t}|^{2}<\infty$, $Px_{t}x_{t}^{\prime}$
	is positive definite, and the distribution of $x_{t}$ puts zero mass
	on each hyperplane. 
	\item [{(b)}] The density $\gamma$ of $u_{t}$ is bounded, differentiable,
	and symmetric around zero, and decreases away from zero. $|u_{t}|$
	has the unique median $\nu_{0}$ and $\dot{\gamma}(\nu_{0})<0$, where
	$\dot{\gamma}$ is the first derivative of $\gamma$.
\end{description}
Except for Assumption D, which allows dependent observations, all
assumptions are similar to the ones in Kim and Pollard (1990, Section
6.3). 

It is known that $\hat{\theta}=\hat{\beta}-\beta_{0}$ is written
as $\hat{\theta}=\arg\max_{\theta}\mathbb{P}_{n}f_{\theta,\hat{\nu}}$,
where 
\[
f_{\theta,\nu}(x,u)=\mathbb{I}\{x^{\prime}\theta-\nu\le u\le x^{\prime}\theta+\nu\},
\]
and $\hat{\nu}=\inf\{\nu:\sup_{\theta}\mathbb{P}_{n}f_{\theta,\nu}\ge\frac{1}{2}\}$.
Let $\nu_{0}=1$ to simplify the notation. Since $\{f_{\theta,\nu}:\theta\in\mathbb{R}^{d},\nu\in\mathbb{R}\}$
is a VC subgraph class, Arcones and Yu (1994, Theorem 1) implies the
uniform convergence $\sup_{\theta,\nu}|\mathbb{P}_{n}f_{\theta,\nu}-Pf_{\theta,\nu}|=O_{p}(n^{-1/2})$.
Thus, the same argument to Kim and Pollard (1990, pp. 207-208) yields
the convergence rate $\hat{\nu}-1=O_{p}(n^{-1/2})$.

Now, we verify the conditions in Theorem 2. By expansions, the condition
in eq. (10) of the paper is verified as 
\begin{eqnarray}
	P(f_{\theta,\nu}-f_{0,1}) & = & P|\{\Gamma(x^{\prime}\theta+\nu)-\Gamma(\nu)\}-\{\Gamma(x^{\prime}\theta-\nu)-\Gamma(-\nu)\}|\nonumber \\
	&  & +P|\{\Gamma(\nu)-\Gamma(1)\}-\{\Gamma(-\nu)-\Gamma(-1)\}|\nonumber \\
	& = & \dot{\gamma}(1)\theta^{\prime}Pxx^{\prime}\theta+o(|\theta|^{2}+|\nu-1|^{2}).\label{eq:expLMS}
\end{eqnarray}
To check Assumption M (iii) for $\{f_{\theta,\nu}:\theta\in\mathbb{R}^{d},\nu\in\mathbb{R}\}$,
pick any $\varepsilon>0$ and decompose 
\[
P\sup_{(\theta,\nu):|(\theta,\nu)-(\theta^{\prime},\nu^{\prime})|<\varepsilon}|f_{\theta,\nu}-f_{(\theta^{\prime},\nu^{\prime})}|^{2}\le P\sup_{(\theta,\nu):|(\theta,\nu)-(\theta^{\prime},\nu^{\prime})|<\varepsilon}|f_{\theta,\nu}-f_{\theta,\nu^{\prime}}|^{2}+P\sup_{\theta:|\theta-\theta^{\prime}|<\varepsilon}|f_{\theta,\nu^{\prime}}-f_{\theta^{\prime},\nu^{\prime}}|^{2},
\]
for $(\theta^{\prime},\nu^{\prime})$ in a neighborhood of $(0,\ldots,0,1)$.
By similar arguments to (\ref{eq:expLMS}), these terms are of order
$|\nu-\nu^{\prime}|^{2}$ and $|\theta-\theta^{\prime}|^{2}$, respectively,
which are bounded by $C\varepsilon$ with some $C>0$ independent
of $\varepsilon$.

We now verify that $\{f_{\theta,1}:\theta\in\mathbb{R}^{d}\}$ satisfies
Assumption M with $h_{n}=1$. By (b), $Pf_{\theta,1}$ is uniquely
maximized at $\theta_{0}=0$. So Assumption M (i) is satisfied. Since
Assumption M (iii) is already shown, it remains to verify Assumption
M (ii). Some expansions (using symmetry of $\gamma(\cdot)$) yield
\begin{eqnarray*}
	\left\Vert f_{\theta_{1},1}-f_{\theta_{2},1}\right\Vert _{2}^{2} & = & P|\Gamma(x^{\prime}\theta_{1}+1)-\Gamma(x^{\prime}\theta_{2}+1)+\Gamma(x^{\prime}\theta_{1}-1)-\Gamma(x^{\prime}\theta_{2}-1)|\\
	& \ge & (\theta_{2}-\theta_{1})^{\prime}P\dot{\gamma}(-1)xx^{\prime}(\theta_{2}-\theta_{1})+o(|\theta_{2}-\theta_{1}|^{2}),
\end{eqnarray*}
i.e., Assumption M (ii) is satisfied under (b). Therefore, $\{f_{\theta,1}:\theta\in\mathbb{R}^{d}\}$
satisfies Assumption M.

We finally derive the finite dimensional convergence through Lemma
2. Let
\[
g_{n,s}(z_{t})=n^{1/6}(\mathbb{I}\{|x_{t}^{\prime}(\theta_{0}+sn^{-1/3})-u_{t}|\le1\}-\mathbb{I}\{|x_{t}^{\prime}\theta_{0}-u_{t}|\le1\}).
\]
Then for any $c>0$, 
\begin{eqnarray}
	&  & P\{|g_{n,s}(z_{t})|>c\}\label{eq:gns}\\
	& \leq & P\{x_{t}^{\prime}(\theta_{0}+sn^{-1/3})\leq u_{t}-1\leq x_{t}^{\prime}\theta_{0}\}+P\{x_{t}^{\prime}(\theta_{0}+sn^{-1/3})\leq u_{t}+1\leq x_{t}^{\prime}\theta_{0}\}\nonumber \\
	&  & +P\{x_{t}^{\prime}\theta_{0}\leq u_{t}-1\leq x_{t}^{\prime}(\theta_{0}+sn^{-1/3})\}+P\{x_{t}^{\prime}\theta_{0}\leq u_{t}+1\leq x_{t}^{\prime}(\theta_{0}+sn^{-1/3})\}.\nonumber 
\end{eqnarray}
However, each of these terms are bounded by $O(n^{-1/3}E|x_{t}^{\prime}s|)$
due to the bounded ness of the density of $u_{t}$, the independence
between $x_{t}$ and $u_{t}$, and the law of iterated expectations
that $P\{h(x_{t},u_{t})\in A\}=EP\{h(x_{t},u_{t})\in A|x_{t}\}$ for
any $h$ and $A$. Thus, by Lemma 2, the central limit theorem in
Lemma C applies. 

It remains to characterize the covariance kernel $H(s_{1},s_{2})$
for any $s_{1}$ and $s_{2}$. Define $q_{n,t}=f_{\theta_{0}+n^{-1/3}s_{1}}(x_{t},u_{t})-f_{\theta_{0}+n^{-1/3}s_{2}}(x_{t},u_{t})$.
Then, the standard algebra yields that
\[
H(s_{1},s_{2})=\frac{1}{2}\{L(s_{1},0)+L(0,s_{2})-L(s_{1},s_{2})\},
\]
where $L(s_{1},s_{2})=\lim_{n\to\infty}n^{4/3}\mathrm{Var}(\mathbb{P}_{n}q_{n,t})=\lim_{n\to\infty}n^{1/3}\{\mathrm{Var}(q_{n,t})+\sum_{m=1}^{\infty}\mathrm{Cov}(q_{n,t},q_{n,t+m})\}.$
The limit of $n^{1/3}\mathrm{Var}(q_{n,t})$ is given by $2\gamma\left(1\right)P|x^{\prime}(s_{2}-s_{1})|$
by direct algebra as in Kim and Pollard (1990, p. 213). For the covariance
$\mathrm{Cov}(q_{n,t},q_{n,t+m})$, note that $q_{n,t}$ can take
only three values, $-1$, $0$, or $1$. By the definition of $\beta_{m}$,
Assumption D implies 
\[
|P\{q_{n,t}=j,q_{n,t+m}=k\}-P\{q_{n,t}=j\}P\{q_{n,t+m}=k\}|\leq n^{-2/3}\beta_{m},
\]
for all $n,m\ge1$ and $j,k=-1,0,1$, i.e., $\{q_{n,t}\}$ is a $\beta$-mixing
array whose mixing coefficients are bounded by $n^{-2/3}\beta_{m}$.
In turn, this implies that $\{q_{n,t}\}$ is an $\alpha$-mixing array
whose mixing coefficients are bounded by $2n^{-2/3}\beta_{m}$. Thus,
by applying the $\alpha$-mixing inequality, the covariance is bounded
as 
\[
\mathrm{Cov}(q_{n,t},q_{n,t+m})\le Cn^{-2/3}\beta_{m}\left\Vert q_{n,t}\right\Vert _{p}^{2},
\]
for some $C>0$ and $p>2$. Note that proceeding as in the bound for
(\ref{eq:gns}) we can show that $\left\Vert q_{n,t}\right\Vert _{p}^{2}=O(n^{-2/(3p)}).$
Combining these results, $n^{1/3}\sum_{m=1}^{\infty}\mathrm{Cov}(q_{n,t},q_{n,t+m})\to0$
as $n\to\infty$. 

Therefore, by Theorem 2, we conclude that $n^{1/3}(\hat{\beta}-\beta_{0})$
converges in distribution to the argmax of $Z(s)$, which is a Gaussian
process with expected value $\dot{\gamma}(1)s^{\prime}Pxx^{\prime}s$
and the covariance kernel $H$, for which $L(s_{1},s_{2})=2\gamma(1)P|x^{\prime}(s_{1}-s_{2})|.$

\subsection{Monotone density}

Preliminary results (Lemmas M, M', C, and 1) to show Theorem 1 may
be applied to establish weak convergence of certain processes. As
an example, consider estimation of a decreasing marginal density function
of $z_{t}$ with support $[0,\infty)$. We impose Assumption D for
$\{z_{t}\}$. The nonparametric maximum likelihood estimator $\hat{\gamma}(c)$
of the density $\gamma(c)$ at a fixed $c>0$ is given by the left
derivative of the concave majorant of the empirical distribution function
$\hat{\Gamma}$. It is known that $n^{1/3}(\hat{\gamma}(c)-\gamma(c))$
can be written as the left derivative of the concave majorant of the
process $W_{n}(s)=n^{2/3}\{\hat{\Gamma}(c+sn^{-1/3})-\hat{\Gamma}(c)-\gamma(c)sn^{-1/3}\}$
(Prakasa Rao, 1969). Let $f_{\theta}(z)=\mathbb{I}\{z\le c+\theta\}$
and $\Gamma$ be the distribution function of $\gamma$. Decompose
\[
W_{n}(s)=n^{1/6}\mathbb{G}_{n}(f_{sn^{-1/3}}-f_{0})+n^{2/3}\{\Gamma(c+sn^{-1/3})-\Gamma(c)-\gamma(c)sn^{-1/3}\}.
\]
A Taylor expansion implies convergence of the second term to $\frac{1}{2}\dot{\gamma}(c)s^{2}<0$.
For the first term $Z_{n}(s)=n^{1/6}\mathbb{G}_{n}(f_{sn^{-1/3}}-f_{0})$,
we can apply Lemmas C and M' to establish the weak convergence. Lemma
C (setting $g_{n}$ as any finite dimensional projection of the process
$\{n^{1/6}(f_{sn^{-1/3}}-f_{0}):s\}$) implies finite dimensional
convergence of $Z_{n}$ to projections of a centered Gaussian process
with the covariance kernel 
\[
H(s_{1},s_{2})=\lim_{n\to\infty}n^{1/3}\sum_{t=-n}^{n}\{\Gamma_{0t}(c+s_{1}n^{-1/3},c+s_{2}n^{-1/3})-\Gamma(c+s_{1}n^{-1/3})\Gamma(c+s_{2}n^{-1/3})\},
\]
where $\Gamma_{0t}$ is the joint distribution function of $(z_{0},z_{t})$.
For stochastic asymptotic equicontinuity of $Z_{n}$, we apply Lemma
M' by setting $g_{n,s}=n^{1/6}(f_{sn^{-1/3}}-f_{0})$. The envelope
condition is clearly satisfied. The condition in eq. (8) of the paper
is verified as 
\begin{eqnarray*}
	&  & P\sup_{s:|s-s^{\prime}|<\varepsilon}|g_{n,s}-g_{n,s^{\prime}}|^{2}\\
	& = & n^{1/3}P\sup_{s:|s-s^{\prime}|<\varepsilon}|\mathbb{I}\{z\le c+sn^{-1/3}\}-\mathbb{I}\{z\le c+s^{\prime}n^{-1/3}\}|\\
	& \le & n^{1/3}\max\{\Gamma(c+sn^{-1/3})-\Gamma(c+(s-\varepsilon)n^{-1/3}),\Gamma(c+(s+\varepsilon)n^{-1/3})-\Gamma(c+sn^{-1/3})\}\\
	& \le & \gamma(0)\varepsilon.
\end{eqnarray*}
Therefore, by applying Lemmas C and M', $W_{n}$ weakly converges
to $Z$, a Gaussian process with expected value $\frac{1}{2}\dot{\gamma}(c)s^{2}$
and covariance kernel $H$.

The remaining part follows by the same argument to Kim and Pollard
(1990, pp. 216-218) (by replacing their Lemma 4.1 with our Lemma 1).
Then we can conclude that $n^{1/3}(\hat{\gamma}(c)-\gamma(c))$ converges
in distribution to the derivative of the concave majorant of $Z$
evaluated at $0$.

\subsection{Binary choice with interval regressor}

we consider a binary choice model with an interval regressor studied
by Manski and Tamer (2002). More precisely, let $y=\mathbb{I}\{x^{\prime}\theta_{0}+w+u\ge0\}$,
where $x$ is a vector of observable regressors, $w$ is an unobservable
regressor, and $u$ is an unobservable error term satisfying $P\{u\le0|x,w\}=\alpha$
(we set $\alpha=.5$ to simplify the notation). Instead of $w$, we
observe the interval $[w_{l},w_{u}]$ such that $P\{w_{l}\le w\le w_{u}\}=1$.
Here we normalize that the coefficient of $w$ to determine $y$ equals
one. In this setup, the parameter $\theta_{0}$ is partially identified
and its identified set is written as (Manski and Tamer 2002, Proposition
2) 
\[
\Theta_{I}=\{\theta\in\Theta:P\{x^{\prime}\theta+w_{u}\le0<x^{\prime}\theta_{0}+w_{l}\mbox{ or }x^{\prime}\theta_{0}+w_{u}\le0<x^{\prime}\theta+w_{l}\}=0\}.
\]
Let $\tilde{x}=(x^{\prime},w_{l},w_{u})^{\prime}$ and $q_{\hat{\nu}}(\tilde{x})$
be an estimator for $P\{y=1|\tilde{x}\}$ with the estimated parameters
$\hat{\nu}$. Suppose $P\{y=1|\tilde{x}\}=q_{\nu_{0}}(\tilde{x})$.
By exploring the maximum score approach, Manski and Tamer (2002) developed
the set estimator for $\Theta_{I}$, that is
\begin{equation}
	\hat{\Theta}=\{\theta\in\Theta:\max_{\theta\in\Theta}S_{n}(\theta)-S_{n}(\theta)\leq\epsilon_{n}\},\label{eq:MTE}
\end{equation}
where
\[
S_{n}(\theta)=\mathbb{P}_{n}(y-.5)[\mathbb{I}\{q_{\hat{\nu}}(\tilde{x})>.5\}\mathrm{sgn}(x^{\prime}\theta+w_{u})+\mathbb{I}\{q_{\hat{\nu}}(\tilde{x})\le.5\}\mathrm{sgn}(x^{\prime}\theta+w_{l})].
\]
Manski and Tamer (2002) established the consistency of $\hat{\Theta}$
to $\Theta_{I}$ under the Hausdorff distance. To establish the consistency,
they assumed the cutoff value $\epsilon_{n}$ is bounded from below
by the (almost sure) decay rate of $\sup_{\theta\in\Theta}|S_{n}(\theta)-S(\theta)|$,
where $S(\theta)$ is the limiting object of $S_{n}(\theta)$. As
Manski and Tamer (2002, Footnote 3) argued, characterization of the
decay rate is a complex task because $S_{n}(\theta)$ is a step function
and $\mathbb{I}\{q_{\hat{\nu}}(\tilde{x})>.5\}$ is a step function
transform of the nonparametric estimate of $P\{y=1|\tilde{x}\}$.
Therefore, it has been an open question. Obtaining the lower bound
rate of $\epsilon_{n}$ is important because we wish to minimize the
volume of the estimator $\hat{\Theta}$ without losing the asymptotic
validity. By applying Theorem 4, we can explicitly characterize the
decay rate for the lower bound of $\epsilon_{n}$ and establish the
convergence rate of this estimator.

A little algebra shows that the set estimator in (\ref{eq:MTE}) is
written as 
\[
\hat{\Theta}=\{\theta\in\Theta:\max_{\theta\in\Theta}\mathbb{P}_{n}f_{\theta,\hat{\nu}}-\mathbb{P}_{n}f_{\theta,\hat{\nu}}\leq\hat{c}n^{-1/2}\},
\]
where $z=(x^{\prime},w,w_{l},w_{u},u)^{\prime}$, $h(x,w,u)=\mathbb{I}\{x^{\prime}\theta_{0}+w+u\ge0\}-\mathbb{I}\{x^{\prime}\theta_{0}+w+u<0\}$,
and 
\begin{equation}
	f_{\theta,\nu}(z)=h(x,w,u)[\mathbb{I}\{x^{\prime}\theta+w_{u}\ge0,q_{\nu}(\tilde{x})>.5\}-\mathbb{I}\{x^{\prime}\theta+w_{l}<0,q_{\nu}(\tilde{x})\le.5\}].\label{eq:MTf}
\end{equation}
We impose the following assumptions. Let $\partial\Theta_{I}$ be
the boundary of $\Theta_{I}$, $\kappa_{u}(\tilde{x})=(2q_{\nu_{0}}(\tilde{x})-1)\mathbb{I}\{q_{\nu_{0}}(\tilde{x})>.5\}$,
and $\kappa_{l}(\tilde{x})=(1-2q_{\nu_{0}}(\tilde{x}))\mathbb{I}\{q_{\nu_{0}}(\tilde{x})\le.5\}$.
\begin{description}
	\item [{(a)}] $\{x_{t},w_{t},w_{lt},w_{ut},u_{t}\}$ satisfies Assumption
	D. $x|w_{u}$ has a bounded and continuous conditional density $p(\cdot|w_{u})$
	for almost every $w_{u}$. There exists an element $x_{j}$ of $x$
	whose conditional density $p(x_{j}|w_{u})$ is bounded away from zero
	over the support of of $w_{u}$. The same condition holds for $x|w_{l}$.
	The conditional densities of $w_{u}|w_{l},x$ and $w_{l}|w_{u},x$
	are bounded. $q_{\nu}(\cdot)$ is continuously differentiable at $\nu_{0}$
	a.s. and the derivative is bounded for almost every $\tilde{x}$. 
	\item [{(b)}] For each $\theta\in\partial\Theta_{I}$, $\kappa_{u}(\tilde{x})$
	is non-negative for $x^{\prime}\theta+w_{u}\ge0$, $\kappa_{l}(\tilde{x})$
	is non-positive for $x^{\prime}\theta+w_{l}\le0$, $\kappa_{u}(\tilde{x})$
	and $\kappa_{l}(\tilde{x})$ are continuously differentiable, and
	it holds 
	\begin{eqnarray*}
		&  & P\{x^{\prime}\theta+w_{u}=0,q_{\nu_{0}}(\tilde{x})>.5,(\theta^{\prime}\partial\kappa_{u}(\tilde{x})/\partial x)p(x|w_{l},w_{u})>0\}>0,\mbox{ or}\\
		&  & P\{x^{\prime}\theta+w_{l}=0,q_{\nu_{0}}(\tilde{x})\le.5,(\theta^{\prime}\partial\kappa_{l}(\tilde{x})/\partial x)p(x|w_{l},w_{u})>0\}>0.
	\end{eqnarray*}
	\item [{(c)}] There exist some $\varepsilon,C>0$ such that for any $|\nu-\nu_{0}|<\varepsilon$,
	$|q_{\nu}(\tilde{x})-q_{\nu_{0}}(\tilde{x})|\leq c|\nu-\nu_{0}|\sqrt{k_{n}}$.
\end{description}
To apply Theorem 4, we verify that $\{f_{\theta,\nu_{0}}:\theta\in\Theta\}$
satisfy Assumption S with $h_{n}=1$. We first check Assumption S
(i). This class is clearly bounded. From Manski and Tamer (2002, Lemma
1 and Corollary (a)), $Pf_{\theta,\nu_{0}}$ is maximized at any $\theta\in\Theta_{I}$
and $\Theta_{I}$ is a bounded convex set. By applying the argument
in Kim and Pollard (1990, pp. 214-215), the second directional derivative
at $\theta\in\partial\Theta_{I}$ with the orthogonal direction outward
from $\Theta_{I}$ is
\begin{eqnarray*}
	&  & -2P\int\mathbb{I}\{x^{\prime}\theta=-w_{u}\}\theta^{\prime}\frac{\partial\kappa_{u}(\tilde{x})}{\partial x}p(x|w_{l},w_{u})(x^{\prime}\theta)^{2}d\sigma_{u}\\
	&  & -2P\int\mathbb{I}\{x^{\prime}\theta=-w_{l}\}\theta^{\prime}\frac{\partial\kappa_{u}(\tilde{x})}{\partial x}p(x|w_{l},w_{u})(x^{\prime}\theta)^{2}d\sigma_{l},
\end{eqnarray*}
where $\sigma_{u}$ and $\sigma_{l}$ are the surface measures on
the boundaries of the sets $\{x:x^{\prime}\pi_{\theta}+w_{u}\ge0\}$
and $\{x:x^{\prime}\pi_{\theta}+w_{l}\ge0\}$, respectively. Since
this matrix is negative definite by (b), Assumption S (i) is verified.
We next check Assumption S (ii). By $h(x,w,u)^{2}=1$, observe that
\[
\left\Vert f_{\theta,\nu_{0}}-f_{\pi_{\theta},\nu_{0}}\right\Vert _{2}\ge\sqrt{2}\min\left\{ \begin{array}{c}
P\{x^{\prime}\theta\ge-w_{u}\ge x^{\prime}\pi_{\theta}\mbox{ or }x^{\prime}\theta<-w_{u}<x^{\prime}\pi_{\theta}\}\mathbb{I}\{q_{\nu_{0}}(\tilde{x})>.5\},\\
P\{x^{\prime}\theta\ge-w_{l}\ge x^{\prime}\pi_{\theta}\mbox{ or }x^{\prime}\theta<-w_{l}<x^{\prime}\pi_{\theta}\}\mathbb{I}\{q_{\nu_{0}}(\tilde{x})\le.5\}
\end{array}\right\} .
\]
for any $\theta\in\Theta$. Since the right hand side is the minimum
of probabilities for pairs of wedge shaped regions with angles of
order $|\theta-\pi_{\theta}|$, (a) implies Assumption S (ii). We
now check Assumption S (iii). By $h(x,w,u)^{2}=1$, the triangle inequality,
and $|\mathbb{I}\{q_{\nu_{0}}(\tilde{x})>0.5\}|\le1$, we obtain 
\begin{eqnarray}
	&  & P\sup_{\theta\in\Theta:0<|\theta-\pi_{\theta}|<\varepsilon}|f_{\theta,\nu_{0}}-f_{\pi_{\theta},\nu_{0}}|^{2}\nonumber \\
	& \le & P\sup_{\theta\in\Theta:0<|\theta-\pi_{\theta}|<\varepsilon}\mathbb{I}\{x^{\prime}\theta\ge-w_{u}\ge x^{\prime}\pi_{\theta}\mbox{ or }x^{\prime}\theta<-w_{u}<x^{\prime}\pi_{\theta}\}\nonumber \\
	&  & +P\sup_{\theta\in\Theta:0<|\theta-\pi_{\theta}|<\varepsilon}\mathbb{I}\{x^{\prime}\theta\ge-w_{l}\ge x^{\prime}\pi_{\theta}\mbox{ or }x^{\prime}\theta<-w_{l}<x^{\prime}\pi_{\theta}\},\label{eq:MT}
\end{eqnarray}
for any $\varepsilon>0$. Again, the right hand side is the sum of
the probabilities for pairs of wedge shaped regions with angles of
order $\varepsilon$. Thus, (a) also guarantees Assumption S (iii).

Next, we verify eqs. (16) and (17). Let $I_{\nu}(\tilde{x})=\mathbb{I}\{q_{\nu}(\tilde{x})>.5\ge q_{\nu_{0}}(\tilde{x})\mbox{ or }q_{\nu}(\tilde{x})\le.5<q_{\nu_{0}}(\tilde{x})\}$
and note that $|f_{\theta,\nu}-f_{\theta,\nu_{0}}|^{2}\leq\mathbb{I}\{x^{\prime}\theta\ge-w_{u}\ge x^{\prime}\pi_{\theta}\mbox{ or }x^{\prime}\theta<-w_{u}<x^{\prime}\pi_{\theta}\}I_{\nu}(\tilde{x})\leq I_{\nu}(\tilde{x})$.
Furthermore, 
\begin{eqnarray*}
	P\sup_{\nu\in\Lambda:|\nu-\nu_{0}|<\varepsilon}\mathbb{I}\{q_{\nu}(\tilde{x})>.5\ge q_{\nu_{0}}(\tilde{x})\} & = & P\sup_{\nu\in\Lambda:|\nu-\nu_{0}|<\varepsilon}\mathbb{I}\{q_{\nu}(\tilde{x})-q_{\nu_{0}}(\tilde{x})>.5-q_{\nu_{0}}(\tilde{x})\ge0\}\\
	& \leq & CP\sup_{\nu\in\Lambda:|\nu-\nu_{0}|<\varepsilon}|q_{\nu}(\tilde{x})-q_{\nu_{0}}(\tilde{x})|\\
	& \leq & C\sqrt{k_{n}}\varepsilon,
\end{eqnarray*}
where the first inequality is due to the boundedness of the conditional
density of $q_{\nu_{0}}(\tilde{x})$ and the second to Condition (c).
This verifies eq. (16), and eq. (17) is verified in the same manner
as Assumption S (ii) in the preceding paragraph.

Finally, for eq. (18), note that 
\begin{eqnarray}
	&  & |P(f_{\theta,\nu}-f_{\theta,\nu_{0}})-P(f_{\pi_{\theta},\nu}-f_{\pi_{\theta},\nu_{0}})|\nonumber \\
	& \le & P\mathbb{I}\{x^{\prime}\theta\ge-w_{u}\ge x^{\prime}\pi_{\theta}\mbox{ or }x^{\prime}\theta<-w_{u}<x^{\prime}\pi_{\theta}\}I_{\nu}(\tilde{x})\nonumber \\
	&  & +P\mathbb{I}\{x^{\prime}\theta\ge-w_{l}\ge x^{\prime}\pi_{\theta}\mbox{ or }x^{\prime}\theta<-w_{l}<x^{\prime}\pi_{\theta}\}I_{\nu}(\tilde{x}),\label{eq:MT1}
\end{eqnarray}
for each $\theta\in\{\theta\in\Theta:|\theta-\pi_{\theta}|<\varepsilon\}$
and $\nu$ in a neighborhood of $\nu_{0}$. For the first term of
(\ref{eq:MT1}), the law of iterated expectation and an expansion
of $q_{\nu}(\tilde{x})$ around $\nu_{0}$ based on (a) imply 
\begin{eqnarray*}
	&  & P\mathbb{I}\{x^{\prime}\theta\ge-w_{u}\ge x^{\prime}\pi_{\theta}\mbox{ or }x^{\prime}\theta<-w_{u}<x^{\prime}\pi_{\theta}\}I_{\nu}(\tilde{x})\\
	& \le & P\mathbb{I}\{x^{\prime}\theta\ge-w_{u}\ge x^{\prime}\pi_{\theta}\mbox{ or }x^{\prime}\theta<-w_{u}<x^{\prime}\pi_{\theta}\}A(w_{u},x)|v-\nu_{0}|,
\end{eqnarray*}
for some bounded function $A$. The second term of (\ref{eq:MT1})
is bounded in the same manner. Therefore, $|P(f_{\theta,\nu}-f_{\theta,\nu_{0}})-P(f_{\pi_{\theta},\nu}-f_{\theta,\nu_{0}})|=O(|\theta-\pi_{\theta}||v-\nu_{0}|)$
and eq. (18) is verified. Since all conditions of Theorem 4 are verified,
we conclude that the convergence rate of Manski and Tamer's (2002)
set estimator $\hat{\Theta}$ in (\ref{eq:MTE}) is characterized
by eqs. (19) and (20).

\subsection{Hough transform estimator}

In the statistics literature on computer vision algorithm, Goldenshluger
and Zeevi (2004) investigated the so-called Hough transform estimator
for the regression model 
\begin{equation}
	\hat{\beta}=\arg\max_{\beta\in B}\sum_{t=1}^{n}\mathbb{I}\{|y_{t}-x_{t}^{\prime}\beta|\le h|x_{t}|\},\label{eq:hough}
\end{equation}
where $B$ is some parameter space, $x_{t}=(1,\tilde{x}_{t})^{\prime}$
for a scalar $\tilde{x}_{t}$, and $h$ is a fixed tuning constant.
Goldenshluger and Zeevi (2004) derived the cube root asymptotics for
$\hat{\beta}$ with fixed $h$, and discussed carefully about the
practical choice of $h$. However, for this estimator, $h$ plays
a role of the bandwidth and the analysis for the case of $h_{n}\to0$
is a substantial open question (see, Goldenshluger and Zeevi, 2004,
pp. 1915-6). Here we focus on the Hough transform estimator in (\ref{eq:hough})
with $h=h_{n}\to0$ and study its asymptotic property. The estimators
by Chernoff (1964) and Lee (1989) with varying $h$ can be analyzed
in the same manner.

Let us impose the following assumptions. 
\begin{description}
	\item [{(a)}] $\{x_{t},u_{t}\}$ satisfies Assumption D. $x_{t}$ and $u_{t}$
	are independent. $P|x_{t}|^{3}<\infty$, $Px_{t}x_{t}^{\prime}$ is
	positive definite, and the distribution of $x_{t}$ puts zero mass
	on each hyperplane. The density $\gamma$ of $u_{t}$ is bounded,
	continuously differentiable in a neighborhood of zero, symmetric around
	zero, and strictly unimodal at zero. 
	\item [{(b)}] As $n\to\infty$, $h_{n}\to0$ and $nh_{n}^{5}\to\infty$. 
\end{description}
Let $z=(x,u)$. Note that $\hat{\theta}=\hat{\beta}-\beta_{0}$ is
written as $\hat{\theta}=\arg\max_{\theta\in\Theta}\mathbb{P}_{n}f_{n,\theta}$,
where 
\[
f_{n,\theta}(z)=h_{n}^{-1}\mathbb{I}\{|u-x^{\prime}\theta|\le h_{n}|x|\}.
\]
The consistency of $\hat{\theta}$ follows from the uniform convergence
$\sup_{\theta\in\Theta}|\mathbb{P}_{n}f_{n,\theta}-Pf_{n,\theta}|\overset{p}{\to}0$
by applying Nobel and Dembo (1993, Theorem 1).

In order to apply Theorem 5, we verify that $\{f_{n,\theta}\}$ satisfies
Assumption M (i), (ii), and (iii'). Obviously $h_{n}f_{n,\theta}$
is bounded. Since $\lim_{n\rightarrow\infty}Pf_{n,\theta}=2P\gamma(x^{\prime}\theta)|x|$
and $\gamma$ is uniquely maximized at zero (by (a)), $\lim_{n\rightarrow\infty}Pf_{n,\theta}$
is uniquely maximized at $\theta=0$. Since $\gamma$ is continuously
differentiable in a neighborhood of zero, $Pf_{n,\theta}$ is twice
continuously differentiable at $\theta=0$ for all $n$ large enough.
Let $\Gamma$ be the distribution function of $\gamma$. An expansion
yields 
\begin{align*}
	P(f_{n,\theta}-f_{n,0}) & =h_{n}^{-1}P\{\Gamma(x^{\prime}\theta+h_{n}|x|)-\Gamma(h_{n}|x|)\}-h_{n}^{-1}P\{\Gamma(x^{\prime}\theta-h_{n}|x|)-\Gamma(-h_{n}|x|)\}\\
	& =\ddot{\gamma}(0)\theta^{\prime}P(|x|xx^{\prime})\theta\{1+O(h_{n})\}+o(|\theta|^{2}),
\end{align*}
i.e., the condition in eq. (2) holds with $V=\ddot{\gamma}(0)P(|x|xx^{\prime})$.
Note that $\ddot{\gamma}(0)<0$ by (a). Therefore, Assumption M (i)
is satisfied.

For Assumption M (ii), pick any $\theta_{1}$ and $\theta_{2}$ and
note that 
\begin{eqnarray*}
	h_{n}\left\Vert f_{n,\theta_{1}}-f_{n,\theta_{2}}\right\Vert _{2}^{2} & = & 2P\{\gamma(x^{\prime}\theta_{1})+\gamma(x^{\prime}\theta_{2})\}|x|\\
	&  & -2h_{n}^{-1}P\{x^{\prime}\theta_{1}-h_{n}|x|<u<x^{\prime}\theta_{2}+h_{n}|x|,\mbox{ }-2h_{n}|x|<x^{\prime}(\theta_{2}-\theta_{1})<0\}\\
	&  & -2h_{n}^{-1}P\{x^{\prime}\theta_{2}-h_{n}|x|<u<x^{\prime}\theta_{1}+h_{n}|x|,\mbox{ }-2h_{n}|x|<x^{\prime}(\theta_{1}-\theta_{2})<0\}.
\end{eqnarray*}
Since the second and third terms converge to zero (by a change of
variable), Assumption M (ii) holds true.

We now check Assumption M (iii'). Observe that 
\begin{eqnarray*}
	P\sup_{\theta\in\Theta:|\theta-\vartheta|<\varepsilon}h_{n}^{2}|f_{n,\theta}-f_{n,\vartheta}|^{2} & \le & P\sup_{\theta\in\Theta:|\theta-\vartheta|<\varepsilon}\mathbb{I}\{|u-x^{\prime}\vartheta|\le h_{n}|x|,\mbox{ }|u-x^{\prime}\theta|>h_{n}|x|\}\\
	&  & +P\sup_{\theta\in\Theta:|\theta-\vartheta|<\varepsilon}\mathbb{I}\{|u-x^{\prime}\theta|\le h_{n}|x|,\mbox{ }|u-x^{\prime}\vartheta|>h_{n}|x|\},
\end{eqnarray*}
for all $\vartheta$ in a neighborhood of $0$. Since the same argument
applies to the second term, we focus on the first term (say, $T$).
If $\varepsilon\le2h_{n}$, then an expansion around $\varepsilon=0$
implies 
\[
T\le P\{(h_{n}-\varepsilon)|x|\le u\le h_{n}|x|\}=P\gamma(h_{n}|x|)|x|\varepsilon+o(\varepsilon).
\]
Also, if $\varepsilon>2h_{n}$, then an expansion around $h_{n}=0$
implies 
\[
T\le P\{-h_{n}|x|\le u\le h_{n}|x|\}\le P\gamma(0)|x|\varepsilon+o(h_{n}).
\]
Therefore, Assumption M (iii') is satisfied.

Finally, the covariance kernel is obtained by a similar way as Section
B.1. Let $r_{n}=(nh_{n}^{2})^{1/3}$ be the convergence rate in this
example. The covariance kernel is written by $H(s_{1},s_{2})=\frac{1}{2}\{L(s_{1},0)+L(0,s_{2})-L(s_{1},s_{2})\}$,
where $L(s_{1},s_{2})=\lim_{n\rightarrow\infty}\mathrm{Var}(r_{n}^{2}\mathbb{P}_{n}g_{n,t})$
with $g_{n,t}=f_{n,s_{1}/r_{n}}-f_{n,s_{2}/r_{n}}$. An expansion
implies $n^{-1}\mathrm{Var}(r_{n}^{2}g_{n,t})\rightarrow2\gamma(0)P|x^{\prime}(s_{1}-s_{2})|$.
We can also see that the covariance term $n^{-1}\sum_{m=1}^{\infty}\mathrm{Cov}(r_{n}^{2}g_{n,t},r_{n}^{2}g_{n,t+m})$
is negligible. Therefore, by Theorem 5, the limiting distribution
of the Hough transform estimator with the bandwidth $h_{n}$ is obtained
as 
\[
(nh_{n}^{2})^{1/3}(\hat{\beta}-\beta_{0})\overset{d}{\rightarrow}\arg\max_{s\in\mathbb{R}^{d}}Z(s),
\]
where $Z(s)$ is a Gaussian process with continuous sample paths,
expected value $\ddot{\gamma}(0)s^{\prime}P(|x|xx^{\prime})s/2$,
and covariance kernel $H(s_{1},s_{2})=2\gamma(0)P|x^{\prime}(s_{1}-s_{2})|$.


\begin{thebibliography}{10}
\bibitem{key-1} Abrevaya, J. and J. Huang (2005) On the bootstrap of the maximum score estimator, \emph{Econometrica}, 73, 1175-1204.

\bibitem{key-2}Adamczak, R. (2008) A tail inequality for suprema of unbounded empirical processes with applications to Markov chains, \emph{Electronic Journal of Probability}, 13, 1000-1034.

\bibitem{key-3}Andrews, D. F., Bickel, P. J., Hampel, F. R., Huber, P. J., Rogers, W. H. and J. W. Tukey (1972) \emph{Robust Estimates of Location}, Princeton University Press, Princeton.

\bibitem{key-5}Anevski, D. and O. H\"{o}ssjer (2006) A general asymptotic scheme for inference under order restrictions, \emph{Annals of Statistics}, 34, 1874-1930.

\bibitem{key-7}Banerjee, M. and I. W. McKeague (2007) Confidence sets for split points in decision trees, \emph{Annals of Statistics}, 35, 543-574.

\bibitem{key-8}Baraud, Y. (2010) A Bernstein-type inequality for suprema of random processes with applications to model selection in non-Gaussian regression, \emph{Bernoulli}, 16, 1064-1085.

\bibitem{key-9}Belloni, A., Chen, D., Chernozhukov, V. and C. Hansen (2012) Sparse models and methods for optimal instruments with application to eminent domain, \emph{Econometrica}, 80, 2369-2429.

\bibitem{key-10}B\"{u}hlmann, P. and B. Yu (2002) Analyzing bagging, \emph{Annals of Statistics}, 30, 927-961.

\bibitem{key-11}Carrasco, M. and X. Chen (2002) Mixing and moment properties of various GARCH and stochastic volatility models, \emph{Econometric Theory}, 18, 17-39.

\bibitem{key-12}Chan, K. S. (1993) Consistency and Limiting Distribution of the Least Squares Estimator of a Threshold Autoregressive Model. \emph{Ann. Statist.} 21, 520-533.

\bibitem{key-13}Chen, X. (2007) Large sample sieve estimation of semi-nonparametric models, \emph{Handbook of Econometrics}, vol. 6B, ch. 76, Elsevier.

\bibitem{key-14}Chen, X., Hansen, L. P. and M. Carrasco (2010) Nonlinearity and temporal dependence, \emph{Journal of Econometrics}, 155, 155-169.

\bibitem{key-15}Chernoff, H. (1964) Estimation of the mode, \emph{Annals of the Institute of Statistical Mathematics}, 16, 31-41.

\bibitem{key-16}Chernozhukov, V., Hong, H. and E. Tamer (2007) Estimation and confidence regions for parameter sets in econometric models, \emph{Econometrica}, 75, 1243-1284.

\bibitem{key-17}de Jong, R. M. and T. Woutersen (2011) Dynamic time series binary choice, \emph{Econometric Theory}, 27, 673-702.

\bibitem{key-19}Doukhan, P., Massart, P. and E. Rio (1995) Invariance principles for absolutely regular empirical processes, \emph{Annales de l'Institut Henri Poincaré, Probability and Statistics}, 31, 393-427.

\bibitem{key-20}Gautier, E. and Y. Kitamura (2013) Nonparametric estimation in random coefficients binary choice models, \emph{Econometrica}, 81, 581-607.

\bibitem{key-21}Goldenshluger, A. and A. Zeevi (2004) The Hough transform estimator, \emph{Annals of Statistics}, 32, 1908-1932.

\bibitem{key-22}Hirano, K. and J. R. Porter, (2003) Asymptotic efficiency in parametric structural models with parameter-dependent support, \emph{Econometrica}, 71, 1307-1338.

\bibitem{key-23}Honor\'{e}, B. E. and E. Kyriazidou (2000) Panel data discrete choice models with lagged dependent variables, \emph{Econometrica}, 68, 839-874.

\bibitem{KP90}Kim, J. and D. Pollard (1990) Cube root asymptotics, \emph{Annals of Statistics}, 18, 191-219.

\bibitem{key-24}Koo, B. and M. H. Seo (2015) Structural Break Models under Misspecification: Implication for Forecasting, \emph{Journal of Econometrics}, 188, 166-181.

\bibitem{key-25}Kosorok, M. R. (2008) \emph{Introduction to Empirical Processes and Semiparametric Inference}, Springer.

\bibitem{key-26}Lee, M.-J. (1989) Mode regression, \emph{Journal of Econometrics}, 42, 337-349.

\bibitem{key-27}Manski, C. F. (1975) Maximum score estimation of the stochastic utility model of choice, \emph{Journal of Econometrics}, 3, 205-228.

\bibitem{key-29}Manski, C. F. and E. Tamer (2002) Inference on regressions with interval data on a regressor or outcome, \emph{Econometrica}, 70, 519-546.

\bibitem{key-30}Merlev\`{e}de, F., Peligrad, M. and E. Rio (2009) Bernstein inequality and moderate deviations under strong mixing conditions, \emph{IMS Collections: High Dimensional Probability V}, 5, 273-292.

\bibitem{key-31}Merlev\`{e}de, F., Peligrad, M. and E. Rio (2011) A Bernstein type inequality and moderate deviations for weakly dependent sequences, \emph{Probability Theory and Related Fields}, 151, 435-474.

\bibitem{key-33}Nickl, R. and J. S\"{o}hl (2016) Nonparametric Bayesian posterior contraction rates for discretely observed scalar diffusions, Working paper, arXiv:1510.05526v2.

\bibitem{key-35}Paulin, D. (2015) Concentration inequalities for Markov chains by Marton couplings and spectral methods, \emph{Electronic Journal of Probability}, 20, 1-32.

\bibitem{key-36}Prakasa Rao, B. L. S. (1969) Estimation of a unimodal density, \emph{Sankhy\={a}}, A, 31, 23-36.

\bibitem{key-37}Politis, D. N., Romano, J. P. and M. Wolf (1999) \emph{Subsampling}, New York: Springer-Verlag.

\bibitem{key-38}Pollard, D. (1989) Asymptotics via empirical processes, Statistica Sinica, 4, 341-354.

\bibitem{key-40}Polonik, W. and Q. Yao (2000) Conditional minimum volume predictive regions for stochastic processes, \emph{Journal of the American Statistical Association}, 95, 509-519.

\bibitem{key-41}Rio, E. (1997) About the Lindeberg method for strongly mixing sequences, \emph{ESAIM: Probability and Statistics}, 1, 35-61.

\bibitem{key-42}Romano, J. P. and A. M. Shaikh (2008) Inference for identifiable parameters in partially identified econometric models, \emph{Journal of Statistical Planning and Inference}, 138, 2786-2807.

\bibitem{key-43}Rousseeuw, P. J. (1984) Least median of squares regression, \emph{Journal of the American Statistical Association}, 79, 871-880.

\bibitem{key-44}Sen, B., Banerjee, M. and M. Woodroofe (2010) Inconsistency of bootstrap: the Grenander estimator, \emph{Annals of Statistics}, 38, 1953-1977.

\bibitem{key-45}Talagrand, M. (2005) \emph{The Generic Chaining}, Springer.

\bibitem{key-46}van der Vaart, A. W. and J. A. Wellner (1996) \emph{Weak Convergence and Empirical Processes}, Springer, New York.

\bibitem{key-47} van der Vaart, A. W. and J. A. Wellner (2007) Empirical processes indexed by estimated functions, \emph{IMS Lecture Notes: Asymptotics: Particles, Processes and Inverse Problems,} 55, 234-252.

\bibitem{key-48}Yao, W., Lindsay, B. G. and R. Li (2012) Local modal regression, \emph{Journal of Nonparametric Statistics}, 24, 647-663.

\bibitem{key-49}Zinde-Walsh, V. (2002) Asymptotic theory for some high breakdown point estimators, \emph{Econometric Theory}, 18, 1172-1196. 
\end{thebibliography}

\begin{thebibliography}{10}
	\bibitem{key-1}Andrews, D. W. K. (1993) An introduction to econometric
	applications of empirical process theory for dependent random variables,
	\emph{Econometric Reviews}, 12, 183-216.
	
	\bibitem{key-8}Arcones, M. A. and B. Yu (1994) Central limit theorems
	for empirical and U-processes of stationary mixing sequences, \emph{Journal
		of Theoretical Probability}, 7, 47-71.
	
	\bibitem{key-10}Chernoff, H. (1964) Estimation of the mode, \emph{Annals
		of the Institute of Statistical Mathematics}, 16, 31-41.
	
	\bibitem{key-2}Doukhan, P., Massart, P. and E. Rio (1995) Invariance
	principles for absolutely regular empirical processes, \emph{Annales
		de l'Institut Henri Poincaré, Probability and Statistics}, 31, 393-427.
	
	\bibitem{key-4}Goldenshluger, A. and A. Zeevi (2004) The Hough transform
	estimator, \emph{Annals of Statistics}, 32, 1908-1932.
	
	\bibitem{key-1}Honoré, B. E. and E. Kyriazidou (2000) Panel data
	discrete choice models with lagged dependent variables, \emph{Econometrica},
	68, 839-874.
	
	\bibitem{key-4}Kim, J. and D. Pollard (1990) Cube root asymptotics,
	\emph{Annals of Statistics}, 18 , 191-219.
	
	\bibitem{key-11}Lee, M.-J. (1989) Mode regression, \emph{Journal
		of Econometrics}, 42, 337-349.
	
	\bibitem{key-5}Manski, C. F. (1975) Maximum score estimation of the
	stochastic utility model of choice, \emph{Journal of Econometrics},
	3, 205-228.
	
	\bibitem{key-7}Manski, C. F. (1985) Semiparametric analysis of discrete
	response: Asymptotic properties of the maximum score estimator, \emph{Journal
		of Econometrics}, 27, 313-333.
	
	\bibitem{key-3}Manski, C. F. and E. Tamer (2002) Inference on regressions
	with interval data on a regressor or outcome, \emph{Econometrica},
	70, 519-546.
	
	\bibitem{key-6}Nobel, A. and A. Dembo (1993) A note on uniform laws
	of averages for dependent processes, \emph{Statistics and Probability
		Letters}, 17, 169-172.
	
	\bibitem{key-2}Polonik, W. and Q. Yao (2000) Conditional minimum
	volume predictive regions for stochastic processes, \emph{Journal
		of the American Statistical Association}, 95, 509-519.
	
	\bibitem{key-9}Prakasa Rao, B. L. S. (1969) Estimation of a unimodal
	density, \emph{Sankhy\={a}}, A, 31, 23-36.
	
	\bibitem{key-3}Rio, E. (1997) About the Lindeberg method for strongly
	mixing sequences, \emph{ESAIM: Probability and Statistics}, 1, 35-61.
\end{thebibliography}
\end{document}